\documentstyle[12pt]{amsart}

\input{amssymb.sty}
\input xy
\xyoption{all}
\theoremstyle{plain}
  \newtheorem{thm}{Theorem}[section]
  \newtheorem{lem}[thm]{Lemma}
  \newtheorem{cor}[thm]{Corollary}
  \newtheorem{prop}[thm]{Proposition}
  \newtheorem{fact}[thm]{Fact}
\theoremstyle{definition}
  \newtheorem{defn}[thm]{Definition}
  
  \newtheorem{clm}[thm]{Claim}

  \newtheorem{notation}{Notation\!\!}

\theoremstyle{remark}
  \newtheorem{rem}[thm]{Remark}
  \newtheorem{acknowledgment}{Acknowledgment}

\numberwithin{equation}{section}

\newcommand{\alp}{\alpha}
\newcommand{\Amp}{{\operatorname{Amp}}}
\newcommand{\bL}{{\bold L}}
\newcommand{\bR}{{\bold R}}
\newcommand{\CC}{{\mathbb{C}}}
\newcommand{\Coh}{{\operatorname{Coh}}}
\newcommand{\Cok}{\operatorname{Cok}}
\newcommand{\DD}{{\bold D}}
\newcommand{\ext}{\operatorname{ext}}
\newcommand{\Ext}{\operatorname{Ext}}
\newcommand{\ff}{{\bf f}}
\newcommand{\Hilb}{\operatorname{Hilb}}
\newcommand{\Hom}{\operatorname{Hom}}
\newcommand{\id}{\operatorname{id}}
\newcommand{\IIm}{\operatorname{Im}}
\newcommand{\Ker}{\operatorname{Ker}}
\newcommand{\KK}{{\cal K}}

\newcommand{\modm}{M_-}
\newcommand{\Modm}{M_-(c_1,c_2)}
\newcommand{\modms}{M_-^s}

\newcommand{\modp}{M_+}
\newcommand{\Modp}{M_+(c_1,c_2)}
\newcommand{\modps}{M_+^s}

\newcommand{\NN}{{\mathbb{N}}}
\newcommand{\NS}{\operatorname{NS}}
\newcommand{\Num}{{\operatorname{Num}}}
\newcommand{\Ox}{{\cal O}_X}

\newcommand{\Pf}{P^{\bf f}}
\newcommand{\Pic}{\operatorname{Pic}}
\newcommand{\PP}{{\mathbb{P}}}
\newcommand{\pr}{\operatorname{pr}}
\newcommand{\Proj}{\operatorname{Proj}}
\newcommand{\Qf}{Q^{\bf f}}
\newcommand{\Qm}{Q_-^{ss}}
\newcommand{\Qms}{Q_-^{s}}
\newcommand{\Qp}{Q_+^{ss}}
\newcommand{\Qps}{Q_+^{s}}
\newcommand{\QQ}{{\mathbb{Q}}}
\newcommand{\Quot}{{\operatorname{Quot}}}
\newcommand{\rk}{{\operatorname{rk}}}

\newcommand{\Spec}{\operatorname{Spec}}
\newcommand{\supp}{\operatorname{supp}}
\newcommand{\Sym}{\operatorname{Sym}}
\newcommand{\tilmodm}{\tilde{M}_-}
\newcommand{\tilmodms}{\tilde{M}_-^s}
\newcommand{\tilmodp}{\tilde{M}_+}

\newcommand{\tilQm}{\tilde{Q}_-^{ss}}
\newcommand{\tilQms}{\tilde{Q}_-^s}
\newcommand{\tilQp}{\tilde{Q}_+^{ss}}

\newcommand{\ZZ}{{\mathbb{Z}}}

\begin{document}
\bibliographystyle{amsplain}

\title[moduli of sheaves under change of polarization]
{A sequence of blowing-ups connecting moduli of\\
sheaves and the Donaldson polynomial\\
 under change of polarization}
\author{Kimiko Yamada}
\address{Department of Mathematics, Faculty of science and technology,
Sophia University, 7-1, Kioi-cho, Chiyoda-ku, Tokyo, Japan}
\email{yamada@@mm.sophia.ac.jp}
\thanks{Partially supported by JSPS Research Fellowship for
Young Scientists. This article has been published as the reference
\cite{Yam:Dthesis} says, while Introduction is revised.
The author submits it to arXiv for the sake of accessibility.}
\maketitle
%
%
\begin{abstract}
Let $H$ and $H'$ be two ample line bundles over a nonsingular projective
surface $X$, and $M(H)$ (resp. $M(H')$) the coarse moduli scheme
of $H$-semistable (resp. $H'$-semistable) sheaves of fixed type
$(r=2,c_1,c_2)$.
In a moduli-theoretic way that comes from elementary transforms, we
connect $M(H)$ and $M(H')$ by a sequence of blowing-ups when walls
separating $H$ and $H'$ are not necessarily good.
As an application, we also consider the polarization change problem
of Donaldson polynomials.
\end{abstract}
%
%
\section*{Introduction}
Let $X$ be a nonsingular projective surface over $\CC$, $H$ an ample line
bundle on
$X$, and $M_H(c_1,c_2)$ the moduli scheme of S-equivalence classes of
rank-two $H$-semistable sheaves on $X$ with fixed Chern classes
$(c_1,c_2)\in \Pic(X)\times\ZZ$. It is projective over $\CC$.
(See Section \ref{ss:background} for walls and chambers.)  \par
Fix two ample line bundles $H_1$ and $H_2$ on $X$ which lie in
neighboring chambers of type $(c_1,c_2)$. In this article,
we first connect $M_{H_1}(c_1,c_2)$ with $M_{H_2}(c_1,c_2)$ by a sequence
of blowing-ups and blowing-downs
\begin{equation}\label{flipintro}
\xymatrix{
& \tilde{M}_1 \ar[dl]_{\phi_1} \ar[dr]^{\psi_1} & & \ar[dl]_{\phi_2}
\tilde{M}_2 \cdots &
\cdots \tilde{M}_l \ar[dr]^{\psi_l} & \\
M_{H_1}(c_1,c_2) & & M_2 & & & M_{H_2}(c_1,c_2) }
\end{equation}
in a natural, moduli-theoretic way which comes from elementary transforms.
As an application, we study the exceptional divisor $E_i$ of $\phi_i$ in
\eqref{flipintro} to observe the following fact, which originates in
differential geometry, from an algebro-geometric view:
Donaldson polynomials of $X$ are independent of the choice of its Riemannian metric
when $b_2^+(X)= 2p_g(X)+1 >1$. \par
We explain the content of this paper.
In Section \ref{ss:background}, we remember some background materials
including the notion of $a$-stability ($0<a<1$), which originally was
proposed by Ellingsrud-G\"{o}ttsche \cite{EG:variation}
and others,
and the coarse moduli scheme $M_a(c_1,c_2)$ of $a$-semistable sheaves
of type $(c_1,c_2)$.
Lemma 3.10 and 3.11 in \cite{EG:variation}, that we shall recall at
Lemma \ref{lem:HNF}, say that when one wants to compare $M_{H_1}$ with
$M_{H_2}$, he suffices to see how $a$-stability of sheaves changes
as parameters $a$ do.
In Section \ref{ss:subscheme} we pick parameters $a_-<a_+$ which are
in adjacent minichambers (Definition \ref{defn:wall}), and endow the
subset
\[ M_- :=M_{a_-}(c_1,c_2) \supset P_- =
   \left\{ [E] \bigm| E\, \text{ is not $a_+$-semistable} \right\} \] 
with a natural subscheme structure.
In Section \ref{ss:morphisms}, we shall connect $M_-$ and 
$M_+= M_{a_+}(c_1,c_2)$.
Suppose for simplicity that $M_-$ has a universal family ${\cal U}_-$.
Then there is a relative family of Harder-Narasimhan filtrations for
$a_+$-stability
\begin{equation}\label{eq:a+HNF}
0 \longrightarrow {\cal F} \longrightarrow 
{\cal U}_1 |_{X\times P_-} \longrightarrow {\cal G}\longrightarrow 0,
\end{equation}
which is an exact sequence of $P_-$-flat sheaves over $X_{P_-}$.
For the blowing-up $\phi_-: \tilde{M_-}\rightarrow M_-$ along $P_-$,
we modify 
$\tilde{{\cal U}}_- := (\operatorname{id}\times \phi_-)^* {\cal U}_-$
to a new $\tilde{M_-}$-flat sheaf ${\cal W}_+$ over $X_{\tilde{M_-}}$
by using elementary transforms.
This ${\cal W}_+$ is accompanied by an exact sequence
\begin{equation}\label{eq:eqseqW}
0 \longrightarrow \tilde{{\cal G}}\otimes {\cal O}_E(-E) \longrightarrow
{\cal W}_+ |_E \longrightarrow \tilde{{\cal F}} \longrightarrow 0,
\end{equation}
where $E$ is the exceptional divisor of $\phi_-$,
$\tilde{{\cal F}}=(\operatorname{id}\times \phi_-)^*{\cal F}$, and so on.
In Lemma \ref{lem:nontrivial} and Lemma \ref{lem:BisWp}, we verify
that this exact sequence \eqref{eq:eqseqW} relates with the
infinitesimal relative obstruction theory of $P_-\subset M_-$.
As a corollary ${\cal W}_+$ is a flat family of $a_+$-semistable
sheaves, so it leads to a morphism
$\overline{\phi}_+: \tilde{M_-}\rightarrow M_+$.
In Section \ref{ss:blowing-up} we prove that $\overline{\phi}_+$ is in
fact the blowing-up of $M_+$ along
\[ M_+ \supset P_+ =
   \left\{ [E] \bigm| E\, \text{ is not $a_-$-semistable} \right\}. \] 
We therefore arrive at blowing-ups
\[  M_- \overset{\phi_-}{\longleftarrow} \tilde{M_-}
    \overset{\overline{\phi}_+}{\longrightarrow} M_+ .  \]
In Section \ref{ss:overPicHilb} we pay attention to the morphism
$E\rightarrow \Pic(X)\times \Hilb(X)\times\Hilb(X)$ deduced from 
${\cal F}$ and ${\cal G}$ at \eqref{eq:a+HNF}. \par
We begin Section \ref{ss:AGanalogy} with reviewing algebro-geometric
analogy of Donaldson polynomials (\cite{Li:AGinterpret})
$ \gamma_H(c_2): \operatorname{Sym}^{d(c_2)} \NS(X)\rightarrow \ZZ$,
where $H$ is an ample line bundle on $H$, $c_2$ an integer, and $d(c_2)$
a certain constant depending on $c_2$. $\gamma_H$ is defined from $M_H(0,c_2)$.
As an application of arguments in the preceding sections,
we prove the following result
in Section \ref{ss:AGanalogy}, \ref{ss:PtimesP} and
\ref{ss:incidence}:
suppose that ample line bundles $H_1$ and $H_2$ are in neighboring
chambers of type $(0,c_2)$ separated by a wall of type $(0,c_2)$, say $W$.
Now denote by $A^+(W)$ the set of all the triples $\ff=(f,m,n)\in
\Num(X)\times \NN^{\times 2}$ which satisfy $f\in 2\Num(X)$,
$H_1\cdot f>0$, $m+n=c_2+(f^2/4)$, and the set
\[ W^f= \left\{ x\in\Num(X) \bigm| x\cdot f=0 \right\} \]
is equal to $W$.
Then, for $\ff\in A^+(W)$ one can define a homomorphism
$C(c_2,\ff): \operatorname{Sym}^{d(c_2)} \NS(X)\rightarrow\ZZ$ such that
\[ \gamma_{H_1}- \gamma_{H_2} =\sideset{}{_{\ff\in A^+(W)}}\sum
 C(c_2,\ff).\]
In Section \ref{ss:subscheme} we shall divide $P_1$ into $\sideset{}{_{\ff\in A^+(a)}}
\coprod P_1^{\ff}$ as a disjoint union of components in a natural way,
and $C(c_2,\ff)$ is the contribution of $P_1^{\ff}$ to $\gamma_{H_1}- \gamma_{H_2}$. 
Let $\Pic^{f/2}(X)$ designate an open subset of $\Pic(X)$
\[ \bigl\{ L\in\Pic(X) \bigm| [2L]=f\; \text{in}\; \Num(X) \bigr\}. \]
\begin{prop}\label{thm:MainIntro}
Suppose that $q(X)=0$ and that some global section $\kappa\in\Gamma(K_X)$
gives a nonsingular curve ${\cal K}\subset X$. Let ${\cal S}$ be any compact
subset of the ample cone $\Amp(X)$.
Then there are constants $d_0({\cal S})$, $d_1(X)$ and $d_2(X)$ depending on
${\cal S}$ such that the following hold:\par
Assume that $\ff=(f,m,n)\in A^+(W)$ satisfies that
\begin{enumerate}
\item the functions $T^{\ff}=\Pic^{f/2}(X)\times \Hilb^m(X)\times
\Hilb^n(X) \rightarrow \ZZ$ defined by 
\begin{align*}
(L,Z_1,Z_2)\mapsto & \dim\Ext^1_X({\cal O}(L)\otimes I_{Z_1}, 
{\cal O}(-L)\otimes I_{Z_2}) \quad\text{and} \\
(L,Z_1,Z_2)\mapsto & \dim\Ext^1_X({\cal O}(-L)\otimes I_{Z_2},
{\cal O}(L)\otimes I_{Z_1})
\end{align*}
are locally-constant, and that
\item $-f^2 > (4/3) c_2 +d_1({\cal S})\sqrt{c_2} +d_2({\cal S})$.
\end{enumerate}
Then $C(c_2,\ff)$ is zero if $c_2\geq d_0({\cal S})$.
\end{prop}
How strong are these conditions (i) and (ii)?
As to (ii), recall that $f\in\NS(X)$ defines a wall of type $(0,c_2)$
if $W^f \cap \Amp(X)\neq\emptyset$, $f\equiv 0 \mod 2\Num(X)$ and
$0 < -f^2 \leq 4c_2$. 
Thus the condition (ii) is relatively weak when $c_2$ is sufficiently
large with respect to ${\cal S}$.
The condition (i) is more strict, while this is always valid when $X$ is
a $K3$ surface. 
We prove Proposition \ref{thm:MainIntro} in 
Section \ref{ss:AGanalogy}, \ref{ss:PtimesP} and \ref{ss:incidence}. 
Roughly speaking, it is important to grasp the structure of the
exceptional divisor of $\phi_i$ at the sequence \eqref{flipintro}. \par
Next, let us explain the background and characteristics of this article.
At least two methods have been developed for the purpose of connecting
$M_{H_1}$ and $M_{H_2}$ by a sequence of morphisms, which are often 
birational.
First, Matsuki-Wentworth \cite{MW:Mumford} succeeded to reduce this
problem to a subject concerning Thaddeus principle \cite{Tha:GITflip},
which considers how the GIT quotient $R^{ss}(L)//G$ varies as a
$G$-linearized polarization $L$ does, where $R$ is a quasi-projective
scheme on which a reductive algebraic group $G$ acts.
Second, by using elementary transforms Ellingsrud-G\"{o}ttsche 
\cite{EG:variation} and Friedman-Qin \cite{FQ:flips} constructed a
diagram of blowing-ups \eqref{flipintro} when walls of type $(c_1,c_2)$
separating $H_1$ and $H_2$ are good. This implies that
\begin{equation}\label{eq:P1}
M_{H_1}(c_1,c_2) \supset P_1 =
   \left\{ [E] \bigm| E\, \text{ is not $H_2$-semistable} \right\} 
\end{equation}
is non-singular.
One can say that the good point of using elementary transforms for this
problem is definiteness;
the centers of blowing-ups in \eqref{flipintro} are
directly described in terms of moduli theory.
One can also relate universal sheaves of moduli spaces
in \eqref{flipintro} very concretely.
Thanks to such definiteness, it should be possible to derive interesting
properties of this flip with the help of moduli theory.
Consequently we can in this article investigate exceptional divisors via
the obstruction theory of universal families.
Above-mentioned papers \cite{EG:variation} and \cite{FQ:flips} 
stimulated the author to write this article. However the hypothesis that walls
are good is rather strong, especially in case where $\kappa(X)>0$.
In this article we make no assumption on walls, so we have to proceed
more carefully.
For example, in obtaining the sequence \eqref{flipintro}, we have to
watch not only tangent spaces but also much more about infinitesimal behavior
of  $P_1\subset M_{H_1}$ at \eqref{eq:P1}; see Section 
\ref{ss:morphisms}, \ref{ss:blowing-up} and \ref{ss:overPicHilb}.
As to the wall-crossing formula of Donaldson polynomials, after
completing this work the author realized that Mochizuki
\cite{Mo:invariants} shown the independence of $\gamma_H(c_2)$ from $H$
when $p_g>0$. Mochizuki's method uses Thaddeus' master spaces and the
localization theorem by Graber-Pandharipande
\cite{Gra-Pan:localization}, and quite differs from ours.
\begin{acknowledgment}
The author is grateful
to Prof.\ Akira Ishii for informing the author of Mochizuki's work, and
giving useful advice especially to Section 3. 
Deep appreciation also goes to Prof. Zhenbo Qin, who gave
valuable advice especially to Lemma \ref{lem:Ti}.
\end{acknowledgment}
\begin{notation}
\begin{enumerate}
\item A scheme is algebraic over $\CC$. For a surface $X$, $\Num(X)$ is
the quotient of $\Pic(X)$ modulo the numerically equivalence.
$\Amp(X)\subset \Num(X) \sideset{}{_{\ZZ}}\otimes \mathbb{R}$ is the 
ample cone of $X$. For a closed subscheme $D$ of $S$, $I_D= I_{D,S}$ means 
its ideal sheaf. The stability of coherent torsion-free sheaves is in the
sense of Gieseker-Maruyama.
\item For $T$-schemes $f: X\rightarrow T$ and $g: S\rightarrow T$, let
$X_S$ denote $X \sideset{}{_T}\times S$. Let ${\cal F}$ be a sheaf on $X$,
and $D\subset T$ a subscheme. We often shorten a sheaf $(\id_X \times g)^*
{\cal F}$ on $X_S$ to $g^* {\cal F}$, and shorten ${\cal F}|_{X_D}$ to
${\cal F}|_D$. $\hom$ and $\ext^i$ indicate, respectively, $\dim\Hom$
and $\dim\Ext^i$.
\end{enumerate}
\end{notation}
%
%
\section{Background materials}\label{ss:background}

In this section let us review some background materials introduced in
\cite{EG:variation} and \cite{Qi:equivalence}. Let $X$ be a nonsingular
surface, and fix a line bundle $c_1$ on $X$ and an integer $c_2$ such
that $4c_2-c_1^2>0$.

\begin{defn}
\begin{enumerate}

\item For $f\in\Num(X)$ we define $W^{f}\subset\Amp(X)$ by
\[ W^{f}=\bigl\{ x\in\Amp(X) \bigm| x\cdot f =0 \bigr\}. \]
$f$ is said to be {\it define a wall of type} $(c_1,c_2)$ if
$W^{f}$ is nonempty, $0<-f^2 \leq 4c_2-c_1^2$ and $f-c_1$ is divisible
by 2 in $\Num(X)$. Then $W^f$ is called {\it a wall of type} $(c_1,c_2)$.

\item {\it A chamber of type $(c_1,c_2)$} is a connected component of
the complement of the union of all walls of type $(c_1,c_2)$.
Two different chambers are said to be {\it neighboring} if the
intersection of their closures contains a nonempty open subset of a wall.

\end{enumerate}
\end{defn}
For an ample line bundle $H$ on $X$ we denote by $M_H(c_1,c_2)$ the
coarse moduli scheme of $H$-semistable rank-two sheaves with Chern classes
$(c_1,c_2)$.
\begin{lem}
\begin{enumerate}

\item For $H$ not contained in any wall of type $(c_1,c_2)$, $M_H(c_1,c_2)$
depends only on the chamber containing $H$.

\item The set of walls of type $(c_1,c_2)$ is locally finite.
\end{enumerate}
\end{lem}
\begin{pf}
(i) is \cite[Proposition 2.7]{EG:variation}. (ii) is \cite[Proposition 2.1.6]
{Qi:equivalence}.
\end{pf}\par
Let $H_+$ and $H_-$ be ample line bundles lying in neighboring chambers
${\cal C}_+$ and ${\cal C}_-$ respectively, and $H$ an ample line bundle
contained in the wall $W$ separating ${\cal C}_+$ and ${\cal C}_-$, and
not contained in any wall but $W$.
Such a setting is natural because of the lemma above.
We can assume that $M=H_+-H_-$ is effective by replacing $H_+$ by its
high multiple if necessary.
\begin{lem}\label{lem:transfertoparab}
There is an integer $n_0$ such that if $E$ is a rank-two sheaf with
Chern classes $(c_1,c_2)$ on $X$ then the following holds for any
integer $l\geq n_0$:
\begin{enumerate}
\item $E$ is $H_-$-stable $($resp. semistable$)$ if and only if $E(-lM)$ is
$H$-stable $($resp. semistable$)$.
\item $E$ is $H_+$-stable $($resp. semistable$)$ if and only if $E(lM)$ is
$H$-stable $($resp. semistable$)$.
\end{enumerate}
\end{lem}
\begin{pf} 
\cite[Page 6, Lemma 3.1]{EG:variation}.
\end{pf}\par
Let $C$ denote $(n_0+1)M$ in this section, where $n_0$ is that in the
lemma above.
\begin{defn}
Let $a$ be a real number between 0  and 1.
\begin{enumerate}
\item For a torsion-free sheaf $E$, we define $P_a(E)$ by
$P_a(E)=[(1-a)\chi(E(-C))+a\chi(E(C))]/\rk(E)$.

\item A torsion-free sheaf $E$ on $X$ is said to be {\it $a$-stable}
(resp. {\it $a$-semistable}) if every subsheaf $F\subsetneq E$ satisfies
$P_a(F(lH))\leq P_a(E(lH))$ (resp. $P_a(F(lH))< P_a(E(lH))$) for
sufficiently large integer $l$.

\item $E$ is $a$-semistable if and only if parabolic sheaf $(E(C), E(-C),
a)$ is parabolic semistable with respect to $H$.
Hence from \cite{Yk:moduli}, there is a coarse moduli scheme of S-equivalence
classes of $a$-semistable rank-two sheaves with Chern classes $(c_1,c_2)$
on $X$, denoted by $M_a(c_1,c_2)$. This is projective over $\CC$.
$M_a^s(c_1,c_2)\subset M_a(c_1,c_2)$ denotes the open subscheme of
$a$-stable sheaves.

\end{enumerate}
\end{defn}
Remark that $M_0(c_1,c_2)$ (resp. $M_1(c_1,c_2)$)is  naturally isomorphic
to $M_{H_-}(c_1,c_2)$ (resp. $M_{H_+}(c_1,c_2)$) by Lemma 
\ref{lem:transfertoparab}. So we would like to study how $M_a(c_1,c_2)$
changes as $a$ varies.
\begin{defn}\label{defn:wall}
For a real number $0\leq a\leq 1$, $A^+(a)$ is the set of $(f,m,n)\in
\Num(X)\times \ZZ_{\geq 0}^2$ satisfying that $W^{f}$ is equal to
the wall $W$ dividing $H_+$ and $H_-$, $H_+\cdot f>0$,
$m+n=c_2-(c_1^2-{f}^2)/4$, and $m-n=\langle f\cdot(c_1-K_X)\rangle/2 + 
(2a-1)\langle f\cdot C\rangle$.
$a$ is called {\it a miniwall} if $A^+(a)$ is nonempty. Remark that the
number of miniwalls is finite.
{\it A minichamber} is a connected component of the complement of the set 
of all miniwalls in $[0,1]$.
Two minichambers are said to be {\it neighboring} if their closures
intersect.
\end{defn}
\begin{lem}\label{lem:HNF}
Let $a_-<a_+$ be in neighboring minichambers separated by a miniwall $a$.
For torsion-free rank-two sheaf $E$ with Chern classes $(c_1,c_2)$, the
following holds.
\begin{enumerate}

\item If $E$ is $a_-$-semistable and not $a_+$-semistable, then $E$ is
given by a nontrivial extension
\begin{equation}\label{HNF}
0 \longrightarrow \Ox(F)\otimes I_{Z_1} \longrightarrow E
\longrightarrow \Ox(c_1-F)\otimes I_{Z_2} \longrightarrow 0,
\end{equation}
where $Z_1$ and $Z_2$ are zero-dimensional subschemes of $X$ such that
\begin{equation}\label{HNF2}
(2F-c_1, l(Z_1), l(Z_2))\in A^+(a).
\end{equation}

\item Conversely suppose that $E$ is given by a nontrivial extension
\eqref{HNF} satisfying \eqref{HNF2}.
Then $E$ is $a_-$-stable, strictly $a$-semistable, and not $b$-semistable
for any $b>a$.
\end{enumerate}
\end{lem}
\begin{pf}
\cite[Lemma 3.10 and 3.11]{EG:variation}.
\end{pf} \par
We fix ample line bundles $H_{\pm}$ and $H$, and neighboring minichambers
$a_-<a_+$ separated by a miniwall $a$.
We shorten $M_{a_{\pm}}(c_1,c_2)$ to $M_{\pm}(c_1,c_2)$ for simplicity.
%
%
\section{Subscheme consisting of not $a_+$-semistable sheaves}
\label{ss:subscheme}

In this section we shall give a natural subscheme structure to a
well-defined subset 
\begin{equation}\label{setforcenter}
 \modm \supset \bigl\{ [E] \bigm| 
   \text{ $E$ is not $a_+$-semistable} \bigr\} 
\end{equation}
contained in $\modms$.
This closed subscheme shall be the center of a blowing-up later.

We begin with a quick review of the construction of $M_{\pm}(c_1,c_2)=
M_{\pm}$ referring to \cite{Yk:moduli}.
Let ${\cal F}_-(c_1,c_2)$ (or ${\cal F}_+(c_1,c_2)$, resp.) denote the
family of all $a_-$-semistable ($a_+$-semistable, resp.) rank-two sheaves
with Chern classes $(c_1,c_2)$ on $X$.
By the boundedness of $a_{\pm}$-semistablity, there is an integer $N_0$
such that the following conditions are satisfied for any $E\in
{\cal F}_-(c_1,c_2)\cup {\cal F}_+(c_1,c_2)$. 

\begin{enumerate}
\item If $m\geq N_0$, then both $E(C)(mH)|_{2C}$ and $E(-C)(mH)$ are
generated by its global sections.
\item If $m\geq N_0$, then $h^i(E(C)(mH)|_{2C})=0$ and
$h^i(E(-C)(mH))=0$ for $i>0$.
\end{enumerate}

We fix an integer $m\geq N_0$. Then $h^0(E(C)(mH))=R$ is
independent of $E\in {\cal F}_+(c_1,c_2)\cup {\cal F}_-(c_1,c_2)$.
We denote by $Q$ the Quot-scheme $\Quot^{P(l)}_{{\cal O}(-C-mH)^{\oplus R}
/X}$, where $P(l)$ is the Hilbert polynomial $\chi(E(lH))$ of $E\in
{\cal F}_{\pm}(c_1,c_2)$.
On $X_Q$ there is the universal quotient sheaf $\tau_0: {\cal O}_{X_Q}
(-C-mH)^{\oplus R} \rightarrow {\cal U}$.
Now let $Q_{\pm}^s$ (or $Q_{\pm}^{ss}$, resp.) be the maximal open
subset of $Q$ such that, for every $t\in Q^s_{\pm}$ ($Q^{ss}_{\pm}$, resp.),
\[ H^0(\tau_0(C+mH)\otimes k(t)): k(t)^{\oplus R} \rightarrow
   H^0({\cal U}(C+mH)\otimes k(t))\]
is isomorphic, ${\cal U}\otimes k(t)$ satisfies the hypothesis (i) and (ii)
above, and ${\cal U}\otimes k(t)$ is $a_{\pm}$-stable ($a_{\pm}$-semistable,
resp.).
Let us denote the universal quotient sheaf of $Q_{\pm}^{ss}$ by
${\cal U}_{\pm}\in\Coh(X_{Q_{\pm}^{ss}})$.
$\overline{G}=\operatorname{PGL}(R,\CC)$ naturally acts on $Q_{\pm}^{ss}$
and $Q_{\pm}^s$.
By \cite{Yk:moduli} we can construct a good quotient of $Q_{\pm}^{ss}$
(or $Q_{\pm}^s$, resp.) by $\overline{G}$ when $m$ is sufficiently large.
This quotient turns out to be the moduli scheme $M_{\pm}(c_1,c_2)$
($M^s_{\pm}(c_1,c_2)$, resp.). 
Moreover, because a $a_{\pm}$-stable sheaf is simple, one can prove that
the quotient map $\pi_{\pm}: Q^s_{\pm} \rightarrow M^s_{\pm}(c_1,c_2)$ is a
principal fiber bundle with group $\overline{G}$ \cite{Mu:GIT} in a
similar fashion to the proof of \cite[Proposition 6.4]{Ma:moduli2}.

Now we try to give a closed-subscheme structure to the subset
\eqref{setforcenter}.
For $\ff=(f,m,n)\in A^+(a)$, we can define a functor
\[ {\cal Q}^{\ff}: (\operatorname{Sch}/\Qm)^{\circ} \rightarrow
(\operatorname{Sets})\]
as follows:
${\cal Q}^{\ff}(S\rightarrow \Qm)$ is the set of all $S$-flat quotient sheaves
${\cal U}_-\sideset{}{_{\Qm}}\otimes {\cal O}_S \rightarrow {\cal G}'$
such that, for
every geometric point $t\in S$, the induced exact sequence
\[ 0 \longrightarrow \Ker \longrightarrow {\cal U}_-\otimes k(t)
   \longrightarrow {\cal G}'\otimes k(t) \longrightarrow 0 \]
satisfies that $(c_1-2c_1({\cal G}'\otimes k(t)), c_2(\Ker), c_2({\cal G}'
\otimes k(t)))=(f,m,n)$.
This functor ${\cal Q}^{\ff}$ is represented by a relative Quot-scheme
$\Qf$, that is projective over $\Qm$. On $X_{\Qf}$ there is the universal
quotient $\tau_{\ff}: {\cal U}_-\otimes {\cal O}_{\Qf} \rightarrow {\cal G}$.
\begin{lem}\label{lem:torsionfree}
${\cal G}\otimes k(s)$ is torsion-free for every closed point $s\in\Qf$.
\end{lem}
\begin{pf}
The proof is by contradiction. Assume that ${\cal G}\otimes k(s)$ is
not torsion-free, and denote its torsion part by $T\neq 0$.
Then we have a new quotient sheaf
\[ {\cal U}_-\otimes k(s) \rightarrow {\cal G}\otimes k(s) \rightarrow
G'={\cal G}\otimes k(s)/T.\]
Then $P_a({\cal G}\otimes k(s)(lH))> P_a(G'(lH))$
if $l$ is sufficiently large. From the definition of $\ff$ and $\Qf$
one can show that
$P_a({\cal G}\otimes k(s)(lH))=P_a({\cal U}_-\otimes k(s)(lH))$
for all $l$.
So the quotient sheaf ${\cal U}_-\otimes k(s)\rightarrow G'$ satisfies
that
\begin{equation}\label{G'fora}
P_a({\cal U}_-\otimes k(s)(lH)) > P_a(G'(lH))
\end{equation}
if $l$ is sufficiently large.
On the other hand
\begin{equation}\label{G'foram}
P_{a_-}({\cal U}_-\otimes k(s)(lH)) \leq P_{a_-}(G'(lH))
\end{equation} 
if $l$ is sufficiently large since ${\cal U}_-\otimes k(s)$ is
$a_-$-semistable.
From \eqref{G'fora}, \eqref{G'foram} and the Riemann-Roch theorem,
there should be an integer $a_-\leq b< a$ such that
$P_b({\cal U}_-\otimes k(s)(lH))=P_b(G'(lH))$
for all $l$.
We can easily prove that $b$ is a miniwall, which contradicts the
choice of $a_-$ and $a$.
\end{pf}\par
\begin{lem}\label{lem:subschofQm}
The structural morphism $i=i^{\ff}: \Qf \rightarrow \Qm$ is a closed
immersion.
\end{lem}
\begin{pf}
For $s\in\Qf$ we put $t=i(s)$. First we claim that their residue fields
satisfy $k(s)=k(t)$. Indeed, any member $\lambda\in\operatorname{Gal}
(k(s)/k(t))$ induces another $k(s)$-valued point
\[ \Spec(k(s)) \overset{\lambda}{\longrightarrow} \Spec(k(s))
   \rightarrow \Qf \]
of $\Qf$. We denote this $k(s)$-valued point by $s'$.
$s$ and $s'$ respectively give exact sequences
\begin{equation}\label{twoHNF}
\xymatrix{0 \ar[r] & K \ar[r] & ({\cal U}_-\otimes k(t))\otimes k(s)
\ar@{=}[d] \ar[r] & {\cal G}\otimes k(s) \ar[r] & 0 \\
0 \ar[r] & K' \ar[r] & ({\cal U}_-\otimes k(t))\otimes k(s') \ar[r] &
{\cal G}\otimes k(s') \ar[r] & 0. }
\end{equation}
Because of the definition of $\ff$ and $\Qf$, it holds that
\[0< \{c_1(K)-c_1({\cal G}\otimes k(s))\}\cdot H_+\quad\text{and that}\quad
0< \{c_1(K')-c_1({\cal G}'\otimes k(s))\}\cdot H_+.\]
Besides, the lemma above tells us that both ${\cal G}\otimes k(s)$ and
${\cal G}\otimes k(s')$ are torsion-free and rank-one.
Thus two horizontal rows in \eqref{twoHNF} respectively give the
Harder-Narasimhan filtration of ${\cal U}_-\otimes k(s)$ with respect
to $H_+$-stability.
Because of the uniqueness of the Harder-Narasimhan filtration, two
quotient sheaves in \eqref{twoHNF} are isomorphic, that is, $s=s'$.
Accordingly $\operatorname{Gal}(k(s)/k(t))=\{ 1\}$, and hence
$k(s)=k(t)$ since $\operatorname{ch}(k(t))=0$.

Next, $i$ is injective and hence finite.
Indeed, suppose that two points $s$ and $s'$ in $\Qf$ satisfy that
$i(s)=i(s')=t$. Then $k(s)=k(s')=k(t)$ as mentioned above, and we have
two exact sequences
\begin{equation*}
\xymatrix{0 \ar[r] & K \ar[r] & {\cal U}_-\otimes k(t)
\ar@{=}[d] \ar[r] & {\cal G}\otimes k(s) \ar[r] & 0 \\
0 \ar[r] & K' \ar[r] & {\cal U}_-\otimes k(t)\otimes k(s') \ar[r] &
{\cal G}\otimes k(s') \ar[r] & 0.}
\end{equation*}
Then one can prove that $s=s'$ in $\Qf$, in the same way as the
preceding paragraph.

Next, $i$ is unramified.
To prove this, we only need to show that the tangent map
$T_t\, i:\, T_t\,\Qf \rightarrow T_s\,\Qm $ is injective.
$t\in\Qf$ gives an exact sequence
\begin{equation}\label{HNFagain}
0 \longrightarrow K \longrightarrow {\cal U}_-\otimes k(s) 
\longrightarrow {\cal G}\otimes k(t) \longrightarrow 0
\end{equation}
on $X_{k(s)}=X_{k(t)}$. By \cite[Page 43]{HL:text}
$\Ker(T_t\, i)$ can be identified with
$\Hom_{X_k(t)}(K, {\cal G}\otimes k(t))$, which is equal to zero
because \eqref{HNFagain} gives the Harder-Narasimhan filtration of
${\cal U}_-\otimes k(s)$.

Last, $i$ is a closed immersion.
Since $i$ is injective and unramified, the fiber $i^{-1}(t)$ is
naturally isomorphic to $\Spec(k(s))$ for $s\in\Qf$.
%
Since $i$ is finite, $i^{-1}(t)$ is isomorphic to $\Spec( i_*{\cal O}_{\Qf}
\otimes k(t))$.
These facts tell us that the natural homomorphism ${\cal O}_{\Qm}\otimes
k(t) \rightarrow i_*{\cal O}_{\Qf}\otimes k(t)$ is surjective since
$k(t)=k(s)$. So ${\cal O}_{\Qf} \rightarrow i_*{\cal O}_{\Qf}$ itself
should be surjective.
This means that a finite morphism $i$ is a closed immersion.
\end{pf} \par
We therefore obtain a closed subscheme $\Qf$ of $\Qm$, which is
contained in $\Qms$ by virtue of Lemma \ref{lem:HNF}.
Remembering the way to define the natural action $\bar{\sigma}: \bar{G}
\times\Qms \rightarrow \Qms$, one can verify the following:
\begin{lem}\label{lem:QfGstable}
Denote by $\overline{\sigma}_-: \overline{G}\times\Qms \rightarrow \Qms$
the natural action of $\overline{G}$ on $\Qms$.
Then the morphism $\id\times \overline{\sigma}_-: \overline{G}\times\Qms
\rightarrow \overline{G}\times\Qms$ satisfies that
$(\id\times\overline{\sigma}_-)(\overline{G}\times\Qf)=\overline{G}
\times\Qf $
as subschemes of $\overline{G}\times\Qf$.
\end{lem}
This lemma means that
\begin{equation}\label{QfdescentdataGQ}
\pr_2^*{\cal O}_{\Qf}= {\cal O}_{\overline{G}\times\Qf}=
(\id_{\overline{G}}\times \overline{\sigma}_-)^*{\cal O}_{\overline{G}\times
\Qf} = \overline{\sigma}_-^*{\cal O}_{\Qf}
\end{equation}
as quotient sheaves of ${\cal O}_{\overline{G}\times\Qms}$.
Since $\pi_-:\Qms \rightarrow \modms$ is a principal fiber bundle with
group $\overline{G}$,
$(\overline{\sigma}_-, \pr_2): \overline{G}\times\Qms \rightarrow
\Qms \sideset{}{_{\modms}}\times \Qms$ is isomorphic.
Thus, the identification \eqref{QfdescentdataGQ} corresponds to an 
isomorphism
\begin{equation}\label{Qfdescentdata}
\alpha_2: \pr_2^*{\cal O}_{\Qf} \rightarrow \pr_1^*{\cal O}_{\Qf}
\end{equation}
of quotient sheaves of ${\cal O}_{\Qms\sideset{}{_M}\times\Qms}$,
where $\pr_i: \Qms \sideset{}{_{\modms}}\times \Qms \rightarrow \Qms$
is the $i$-th projection for $i=1,2$.
Since \eqref{QfdescentdataGQ} results from Lemma \ref{lem:QfGstable},
one can check that the isomorphism \eqref{Qfdescentdata} satisfies that
$\pr_{12}^*(\alpha_2) \circ \pr_{23}^*(\alpha_2) = 
 \pr_{13}^*(\alpha_2)$, where $\pr_{ij}: \Qms\sideset{}{_{\modms}}\times
\Qms \sideset{}{_{\modms}}\times \Qms \rightarrow 
\Qms \sideset{}{_{\modms}}\times \Qms$ is the $(i,j)$-th projection.

By faithfully-flat quasi-compact descent theory, we get
a coherent sheaf ${\cal F}$ on $\modm$
and a homomorphism $p':{\cal O}_{\modms} \rightarrow {\cal F}$
such that $\pi_-^*{\cal F}={\cal O}_{\Qf}$ and that $\pi_-^*(p')=p^{\ff}$.
This $p': {\cal O}_{\modms} \rightarrow {\cal F}$ should be surjective
since $\pi_-: \Qms\rightarrow\modms$ is faithfully-flat, and hence $p'$
gives a closed subscheme $\Pf$ of $\modms$ such that
$\pi_-^{-1}(\Pf)=\Qf$.
On the other hand $\Qf$ is a closed subscheme of $\Qm$ fixed by
$\overline{G}$, and so $\Pf=\pi_-(\Qf)$ is closed not only in $\modms$ 
but also in
$\modm$ by the property of good quotient.
Summarizing:
\begin{lem}\label{lem:setforcenter}
The closed subscheme $\Qf$ of $\Qm$ obtained in Lemma \ref{lem:subschofQm}
descends to a closed subscheme $\Pf$ of $\modm$ such that $\pi_-^{-1}
(\Pf)=\Qf$, where $\pi_-: \Qm\rightarrow\modm$ is the quotient map.
$\Pf$ is contained in $\modms$.
Set-theoretically, $\sideset{}{_{\ff\in A^+(a)}}\coprod \Pf$ coincides
with the subset \eqref{setforcenter}.
Both $\Qf \cap Q^{\ff'}$ and $\Pf \cap P^{\ff'}$ are empty if $\ff$ and
$\ff'$ are mutually different member of $A^+(a)$.
\end{lem}
%
At the end of this section, we define a closed subset
\begin{equation}\label{setforcenterp}
 \modp \supset \bigl\{ [E] \bigm| 
   \text{ $E$ is not $a_-$-semistable} \bigr\}
\end{equation}
similarly to the above $\modm$. First we define $-\ff$.
\begin{defn}\label{defn:-ff}
For $\ff=(f,m,n)\in A^+(a)$, we define $-\ff\in\Num(X)\times
\ZZ_{\geq 0}^{\times 2}$ by $-\ff=(-f,n,m)$.
\end{defn}
In the same way as the case of $i^{\ff}: \Qf\rightarrow\Qm$ and
$\Pf\subset\Qm$, we can show that
a projective $\Qp$-scheme $Q^{-\ff}$ can be defined; the structural
morphism $i^{-\ff}: Q^{-\ff}\rightarrow \Qp$ is a closed immersion
which factors through $\Qps$;
by using faithfully-flat quasi-compact descent theory, we can obtain
a closed subscheme $P^{-\ff}\subset\modps$ such that $\pi_+^{-1}
(P^{-\ff})=Q^{-\ff}$;
for different members  $\ff$ and $\ff'$ of $A^+(a)$, we see
that $P^{-\ff}\cap P^{-\ff'}$ is empty;
set-theoretically, $\sideset{}{_{\ff\in A^+(a)}}\coprod P^{-\ff}$
coincides with the subset \eqref{setforcenterp} of $\modp$.
%
%
\section{A sequence of morphisms connecting $\Modm$ with $\Modp$}
\label{ss:morphisms}
Let $V_-$ be a closed subscheme $\sideset{}{_{\ff\in A^+(a)}}\coprod \Qf$
of $\Qm$, and $\varphi_-: \tilQm \rightarrow \Qm$ the blowing-up of $\Qm$
along $V_-$, with exceptional divisor $D_-$.
Similarly, let $P_-$ be a closed subscheme $\sideset{}{_{\ff\in A^+(a)}}
\coprod \Pf$ of $\modm$, and $\phi_-: \tilmodm \rightarrow \modm$ the
blowing-up of $\modm$ along $P_-$, with exceptional divisor $E_-$.
\begin{equation*}
\xymatrix{
V_-= \sideset{}{_{\ff}}\coprod \Qf\, \ar@{^{(}-}[r]  &
\Qm \ar[d]^{\pi_-} & \tilQm \ar[l]_{\varphi_-} \ar[d]^{\tilde{\pi}_-}
& \,D_- \ar@{_{(}-}[l]  \\
P_-= \sideset{}{_{\ff}}\coprod \Pf\, \ar@{^{(}-}[r] & \modm &
\tilmodm \ar[l]_{\phi_-} & \,E_- \ar@{_{(}-}[l] }
\end{equation*}
Because $\varphi_-^{-1}\pi_-^{-1}(P_-)=\varphi_-(V_-)=D_-$ is an
effective Cartier divisor on $\tilQm$, a morphism $\tilde{\pi}_-$
is induced.
In this section, we begin with constructing a morphism $\tilde{\varphi}_+:
\tilQm \rightarrow \modp$  using the method of elementary
transformation.
%
Joining the universal quotient sheaf
${\cal U}_-|_{X_{\Qf}} \twoheadrightarrow {\cal G}={\cal G}^{\ff}$
of $\Qf$, we have a quotient sheaf
${\cal U}_-|_{X_{V_-}} \twoheadrightarrow {\cal G}$
on $X_{V_-}=\sideset{}{_{\ff}}\coprod X_{\Qf}$.
This results in an exact sequence
\begin{equation}\label{HNFm}
0 \longrightarrow {\cal F} \longrightarrow {\cal U}_-|_{V_-}
\longrightarrow {\cal G} \longrightarrow 0
\end{equation}
of $V_-$-flat $X_{V_-}$-modules.
Pulling back this by $\id_X\times\varphi_-: X_{D_-}\rightarrow X_{V_-}$,
we get an exact sequence of $D_-$-flat sheaves
\begin{equation}\label{HNFtilQm}
0 \longrightarrow \tilde{\cal F} \longrightarrow \tilde{\cal U}_-|_{D_-}
\longrightarrow \tilde{\cal G} \longrightarrow 0
\end{equation}
on $X_{D_-}$.
Now let ${\cal W}_+$ denote
$\Ker( \tilde{\cal U}_- \twoheadrightarrow \tilde{\cal U}_-|_{D_-}
\rightarrow \tilde{\cal G})$, that is,
\begin{equation}\label{elemtransm}
0 \longrightarrow {\cal W}_+ \longrightarrow \tilde{\cal U}_-
\longrightarrow \tilde{\cal G} \longrightarrow 0
\end{equation}
is exact.
From \cite[Lemma A.3]{Fr:elliptic} ${\cal W}_+$ is flat over $\tilQm$.
\eqref{HNFtilQm} and
\eqref{elemtransm} induce a commutative diagram on $X_{\tilQm}$
\begin{equation}\label{Elemtransm1}
\xymatrix{
 & & 0 \ar[d] & 0 \ar[d] &  \\
0 \ar[r] & \tilde{\cal U}_-(-D_-) \ar[r] \ar@{=}[d] &
{\cal W}_+ \ar[r]^f \ar[d] & \tilde{\cal F} \ar[r] \ar[d] & 0 \\
0 \ar[r] & \tilde{\cal U}_-(-D_-) \ar[r] & \tilde{\cal U}_- \ar[r]
\ar[d]^h & \tilde{\cal U}_-|_{D_-} \ar[r] \ar[d] & 0 \\
 & & \tilde{\cal G} \ar[d] \ar@{=}[r] & \tilde{\cal G} \ar[d] & \\
 & & 0 & 0 & }
\end{equation}
whose rows and columns are exact.
The second column of \eqref{Elemtransm1} gives rise to an exact sequence
\begin{equation*}
0 \longrightarrow Tor^{X_{\tilQm}}_1(\tilde{\cal G}, {\cal O}_{X_{D_-}})
=\tilde{\cal G}(-D_-) \longrightarrow {\cal W}_+|_{X_{D_-}} \longrightarrow
\tilde{\cal U}_-|_{X_{D_-}} \longrightarrow \tilde{\cal G} \longrightarrow 0.
\end{equation*}
From \eqref{Elemtransm1}, this results in an exact sequence
\begin{equation}\label{inducedHNFforWp}
0 \longrightarrow \tilde{\cal G}(-D_-) \longrightarrow {\cal W}_+|_{X_{D_-}}
\overset{f|_{D_-}}{\longrightarrow} \tilde{\cal F} \longrightarrow 0.
\end{equation}
\eqref{inducedHNFforWp}
and the first row of \eqref{Elemtransm1} induce the following commutative
diagram on $X_{\tilQm}$:
\begin{equation}\label{Elemtransm2}
\xymatrix{
 & & 0 \ar[d] & 0 \ar[d] & \\
0 \ar[r] & {\cal W}_+(-D_-) \ar[r] \ar@{=}[d] & 
\tilde{\cal U}_-(-D_-) \ar[d] \ar[r]^{\bar{h}} & 
\tilde{\cal G}(-D_-) \ar[d] \ar[r] & 0 \\
0 \ar[r] & {\cal W}_+(-D_-) \ar[r] & {\cal W}_+ \ar[r] \ar[d]^f &
{\cal W}_+|_{D_-} \ar[r] \ar[d]^{f|_{D_-}} & 0 \\
 & & \tilde{\cal F} \ar@{=}[r] \ar[d] & \tilde{\cal F} \ar[d] & \\
 & & 0 & 0 & }
\end{equation}
such that its second column is equal to the first row of \eqref{Elemtransm1},
and that all rows and columns are exact.
%
For homomorphisms $h$ in \eqref{Elemtransm1} and $\bar{h}$ in
\eqref{Elemtransm2}, one can find an isomorphism
$j_g: \tilde{\cal G}(-D_-)\rightarrow \tilde{\cal G}(-D_-)$ such that
\begin{equation}\label{twoquotientsofUm}
\xymatrix{
\tilde{\cal U}_-(-D_-) \ar[r]^{\bar{h}} \ar@{=}[d] &
\tilde{\cal G}(-D_-) \ar[d]^{j_g} \\
\tilde{\cal U}_-(-D_-) \ar[r]^{h(-D_-)} & \tilde{\cal G}(-D_-) }
\end{equation}
is commutative, in view of the uniqueness of the Harder-Narasimhan
filtration and the simplicity of torsion-free rank-one sheaf. \par
%
Now we recall some obstruction theory. 
By the exact sequence
\begin{equation}\label{2D}
0 \longrightarrow {\cal O}_{D_-}(-D_-) \longrightarrow
{\cal O}_{2D_-} \longrightarrow {\cal O}_{D_-} \longrightarrow 0.
\end{equation}
and \eqref{HNFtilQm}, we have the following commutative
diagram on $X_{2D_-}$ whose rows and columns are exact:
\begin{equation}\label{obstdiag}
\xymatrix{
 & & 0 \ar[d] & & \\
0 \ar[r] & \tilde{\cal F}(-D_-) \ar[r] & \tilde{\cal U}_-\otimes 
{\cal O}_{D_-}(-D_-) \ar[r] \ar[d] & \tilde{\cal G}(-D_-) \ar[r] &
0 \\
 & & \tilde{\cal U}_-|_{X_{2D_-}} \ar[d] & & \\
0 \ar[r] & \tilde{\cal F} \ar[r] & \tilde{\cal U}_-|_{X_{D_-}} \ar[d]
\ar[r] & \tilde{\cal G} \ar[r] & 0 \\
 & & 0 & & }
\end{equation}
From this we can get a complex
$\tilde{\cal F}(-D_-) \overset{F}{\longrightarrow}
\tilde{\cal U}_-|_{X_{2D_-}} \overset{G}{\longrightarrow}
\tilde{\cal G}$,
and check that its middle cohomology $B=\Ker G/ \IIm F$ is a
${\cal O}_{X_{D_-}}$-module.
Then, again from \eqref{obstdiag} we can deduce an exact sequence
\begin{equation}\label{obstexseq}
0 \longrightarrow \tilde{\cal G}(-D_-) \overset{p}{\longrightarrow}
B \overset{q}{\longrightarrow} \tilde{\cal F} \longrightarrow 0
\end{equation}
of $D_-$-flat ${\cal O}_{X_{D_-}}$-modules.
\begin{lem}\label{lem:nontrivial}
The following conditions are equivalent for a closed point $t$ of
$D_-$:
\begin{enumerate}
\item The exact sequence
 \begin{equation}\label{obstexseqFibr}
 0 \longrightarrow \tilde{\cal G}(-D_-)\otimes k(t) \longrightarrow
 B\otimes k(t) \longrightarrow \tilde{\cal F}\otimes k(t) \longrightarrow 0
 \end{equation}
induced from \eqref{obstexseq} is trivial;
\item Let $\tilde{m}_t \subset {\cal O}_{\tilQm}$ be the maximal ideal
defining $t$ and $l$ the integer such that $I_{D_-,t}\subset \tilde{m}_t^l$
and that $I_{D_-,t}\not\subset \tilde{m}_t^{l+1}$.
Then there is a morphism $p_{l+1}: \Spec( {\cal O}_{\tilQm}/ 
\tilde{m}_t^{l+1}) \rightarrow V_-= \sideset{}{_{\ff}}\coprod \Qf$
such that
\begin{equation}\label{Pl+1}
\xymatrix{
\Spec({\cal O}_{\tilQm}/ {\cal O}(-D_-)+\tilde{m}_t^{l+1})
\ar@{_{(}->}[d] \ar@{^{(}->}[r] & D_- \ar[r]^{\varphi_-} &
V_- \ar@{_{(}->}[d] \\
\Spec({\cal O}_{\tilQm}/ \tilde{m}_t^{l+1}) \ar@{^{(}->}[r]
\ar[urr]_{p_{l+1}} & \tilQm \ar[r]_{\varphi_-} & \Qm }
\end{equation}
is commutative.
\end{enumerate}
\end{lem}
\begin{pf}
We put $A={\cal O}_{\tilQm} / \tilde{m}_t^{l+1}$ and
$A'={\cal O}_{\tilQm} / ({\cal O}(-D_-)+\tilde{m}_t^{l+1})$, which are
Artinian local rings.
Tensoring $A$ to \eqref{2D}, we have the following commutative diagram
whose rows are exact.
\begin{equation}\label{AtoAd}
\xymatrix{
 & {\cal O}_{D_-}(-D_-) \sideset{}{_{\tilQm}}\otimes A \ar[r]
\ar@{->>}[d]^q  & {\cal O}_{2D_-} \sideset{}{_Q}\otimes A \ar[r]
\ar@{=}[d] & {\cal O}_{D_-} \sideset{}{_Q}\otimes A \ar[r]
\ar@{=}[d] & 0 \\
0 \ar[r] & I={\cal O}(-D_-)+\tilde{m}_t^{l+1} / \tilde{m}_t^{l+1}
\ar[r] & A \ar[r] & A' \ar[r] & 0}
\end{equation}
Remark that $I$ is a $k(t)$-module because of the choice of $l$.
From its bottom row and \eqref{HNFtilQm} we get the following commutative
diagram on $X_A$ whose rows and columns are exact, similarly to
\eqref{obstdiag}:
\begin{equation}\label{obstdiagLocal}
\xymatrix{ 
 & & 0 \ar[d] & & \\
0 \ar[r] & \tilde{\cal F}_{k(t)}\otimes I \ar[r] & \tilde{\cal U}_{k(t)}
\otimes I \ar[r] \ar[d] & \tilde{\cal G}_{k(t)}\otimes I \ar[r] & 0 \\
 & & \tilde{\cal U}_-|_{X_A} \ar[d] & & \\
0 \ar[r] & \tilde{\cal F} \sideset{}{_{D_-}}\otimes A' \ar[r] &
\tilde{\cal U}_-|_{X_{A'}} \ar[r] \ar[d] & \tilde{\cal G} \sideset{}{_{D_-}}
\otimes A' \ar[r] & 0 \\
 & & 0 & & }
\end{equation}
Then one can deduce a complex $\tilde{\cal F}_{k(t)}\otimes I 
\overset{F'}{\longrightarrow} \tilde{\cal U}_-|_{X_{A'}}
\overset{G'}{\longrightarrow} \tilde{\cal G} \sideset{}{_{D_-}}\otimes
A'$ and an exact sequence of $X_{A'}$-modules
\begin{equation}\label{obstexseqLocal}
0 \longrightarrow \tilde{\cal G}_{k(t)} \otimes I \longrightarrow
B'=\Ker G'/ \IIm F' \longrightarrow \tilde{\cal F} 
\sideset{}{_{D_-}}\otimes A' \longrightarrow 0.
\end{equation}
Now recall that obstruction
theory shows the following fact \cite[Page 43]{HL:text}.
\begin{fact}\label{fact:exseq}
The exact sequence \eqref{obstexseqLocal} is trivial if and only if the
condition $(ii)$ in Lemma \ref{lem:nontrivial} is satisfied.
\end{fact}
From the commutativity of \eqref{AtoAd}, we can make a homomorphism
$B\sideset{}{_{D_-}}\otimes A' \rightarrow B'$ such that
\begin{equation}\label{Twoobstexseq}
\xymatrix{
0 \ar[r] & \tilde{\cal G}(-D_-) \sideset{}{_{D_-}}\otimes A' =
\tilde{\cal G}(-D_-) \sideset{}{_Q}\otimes A \ar[r] 
\ar@{->>}[d]^{\id\otimes q} & B \sideset{}{_{D_-}}\otimes A' \ar[r]
\ar@{->>}[d] & \tilde{\cal F} \sideset{}{_{D_-}}\otimes A' \ar[r]
\ar@{=}[d] & 0 \\
0 \ar[r] & \tilde{\cal G}_{k(t)}\otimes I = \tilde{\cal G} 
\sideset{}{_{\tilQm}}\otimes I \ar[r] & B' \ar[r] &
\tilde{\cal F} \sideset{}{_{D_-}}\otimes A' \ar[r] & 0 }
\end{equation}
is commutative, where the first row is obtained by tensoring $A$ to
\eqref{obstexseq}, and the second row is \eqref{obstexseqLocal}.
Further, the homomorphism $q$ in \eqref{AtoAd} gives a surjective
homomorphism 
$q\otimes k(t): {\cal O}_{D_-} \sideset{}{_{D_-}}\otimes k(t)
\twoheadrightarrow I$,
which should be isomorphic because $\rk_{k(t)} {\cal O}_{D_-}(-D_-)\otimes
k(t)=1$ and $I\neq 0$.
Accordingly we obtain a commutative diagram
\begin{equation}\label{extclasses}
\xymatrix{
\Ext^1_{X_{A'}}(\tilde{\cal F}\otimes A', \tilde{\cal G}(-D_-)\otimes A')
\ar[r]^{q_*} \ar[d]^{(\pi_t)_*} & \Ext^1_{X_{A'}}(\tilde{\cal F}\otimes A',
\tilde{\cal G}_{k(t)}\otimes I) \\
\Ext^1_{X_{A'}}(\tilde{\cal F}\otimes A', \tilde{\cal G}(-D_-)\otimes k(t))
\ar[ur]^{(q\otimes k(t))_*} & \Ext^1_{X_{k(t)}}(\tilde{\cal F}_{k(t)},
\tilde{\cal G}(-D_-)_{k(t)} \ar[l]^-{\pi_t^*}), }
\end{equation}
where $\pi_t$ is a natural homomorphism $A'\rightarrow k(t)$.
Remark that $\pi_t^*$ is isomorphic since $\tilde{\cal F}$ is $D_-$-flat.
Let $\lambda\in\Ext_{X_{A'}}^1(\tilde{\cal F}\otimes A', \tilde{\cal G}
(-D_-)\otimes A')$
be the extension class of the first row of \eqref{Twoobstexseq}.
Then one can prove that $(\pi_t^*)^{-1}(\pi_{t*}(\lambda))$ is the
extension class of \eqref{obstexseqFibr} and that $q_*(\lambda)$ is
the extension class of \eqref{obstexseqLocal} by using the commutativity
of \eqref{Twoobstexseq}.
Because $(q\otimes k(t))_*$ is isomorphic,
$(\pi_t^*)^{-1}(\pi_{t*}(\lambda))=0$ if and only if $q_*(\lambda)=0$.
This and Fact \ref{fact:exseq} complete the proof of this lemma.
\end{pf}\par
\begin{lem}\label{lem:BisWp}
There is an isomorphism $r_0: {\cal W}_+|_{X_{D_-}} \rightarrow B$
such that the following diagram is commutative:
\begin{equation}\label{BisWp}
\xymatrix{
0 \ar[r] & \tilde{\cal G}(-D_-) \ar[r] \ar[d]^{j_g} & {\cal W}_+|_{X_{D_-}}
\ar[r] \ar[d]^{r_0} & \tilde{\cal F} \ar[r] \ar@{=}[d] & 0 \\
0 \ar[r] & \tilde{\cal G}(-D_-) \ar[r]^p & B \ar[r]^q & \tilde{\cal F}
\ar[r] & 0 }
\end{equation}
Here the first row is the third column of \eqref{Elemtransm2}, the second
row is \eqref{obstexseq}, and $j_g$ is the isomorphism in
\eqref{twoquotientsofUm}.
\end{lem}
\begin{pf}
Tensoring ${\cal O}_{2D_-}$ to \eqref{Elemtransm1}, we have a
commutative diagram on $X_{2D_-}$
\begin{equation}\label{Elemtrans2D}
\xymatrix{
 & \tilde{\cal G}(-2D_-) \ar[r] \ar@{_{(}->}[d] & 0 \ar[d] & \\
\tilde{\cal U}_-(-D_-)|_{2D_-} \ar[r]^-{k'} \ar@{=}[d] &
{\cal W}_+|_{2D_-} \ar[r]^-{f|_{2D_-}} \ar[d]^r & \tilde{\cal F}
\ar[r] \ar[d] & 0 \\
\tilde{\cal U}_-(-D_-)|_{2D_-} \ar[r] & \tilde{\cal U}_-|_{2D_-}
\ar[r] \ar[d]_{h|_{2D_-}=G} & \tilde{\cal U}_-|_{D_-} \ar[r] \ar[d] &
0 \\
 & \tilde{\cal G} \ar[d] \ar@{=}[r] & \tilde{\cal G} \ar[d] & \\
 & 0 & 0 & }
\end{equation}
whose rows and columns are exact. In this diagram $h|_{2D_-}$ clearly is
equal to the homomorphism $G$ defined just below \eqref{obstdiag},
and so $r$ factors into
${\cal W}_+|_{2D_-} \overset{r_1}{\twoheadrightarrow} \IIm r=\Ker G
\rightarrow \tilde{\cal U}_-|_{2D_-}$.
One can readily check that
\begin{equation}\label{r0defn}
\xymatrix{ 
{\cal W}_+|_{2D_-} \ar@{->>}[r] \ar@{->>}[d]^{r_1} & {\cal W}_+|_{D_-}
\ar@{.>}[d]^{r_0} \ar@{->>}[rd]^{f|_{D_-}} & \\
\Ker G \ar@{->>}[r] & B=\Ker G/\IIm F \ar[r]_-{q} & \tilde{\cal F}}
\end{equation}
is commutative by the definition of $q$ in \eqref{obstexseq}.
Since $B$ is naturally regarded as an ${\cal O}_{X_{D_-}}$-module,
we can induce a homomorphism $r_0: {\cal W}_+|_{D_-} \rightarrow B$
such that the left side of \eqref{r0defn} becomes commutative.
Then one can also check the right side of \eqref{r0defn} is commutative,
since ${\cal W}_+|_{2D_-} \rightarrow {\cal W}_+|_{D_-}$ is surjective.
Therefore the right side of \eqref{BisWp} is surely commutative for
this $r_0$.

Next, by the definition of $p$ in \eqref{obstexseq} one can readily check that
\begin{equation*}
\xymatrix{
\tilde{\cal U}_-(-D_-)|_{2D_-} \ar@{->>}[rr]^-{h(-D_-)|_{2D_-}}
\ar[d]^{k'} && \tilde{\cal G}(-D_-) \ar[d]^p \\
{\cal W}_+|_{2D_-} \ar@{->>}[r]^{r_1} & \Ker G \ar@{->>}[r] & B}
\end{equation*}
is commutative, where $h|_{2D_-}$ and $k'$ are those of \eqref{Elemtrans2D}.
We have also the following commutative diagram:
\begin{equation*}
\xymatrix{
\tilde{\cal U}_-(-D_-)|_{2D_-} \ar[r]^{k'} \ar[d]^{\bar{h}|_{2D_-}} &
{\cal W}_+|_{2D_-} \ar[r]^{r_1} \ar@{->>}[d] & \Ker G \ar@{->>}[d] \\
\tilde{\cal G}(-D_-) \ar[r] & {\cal W}_+|_{D_-} \ar@{->>}[r]^{r_0} & B, }
\end{equation*}
where the left side is the upper-right side of \eqref{Elemtransm2},
and the right side is the left side of \eqref{r0defn}.
These two commutative diagrams gives rise to a commutative diagram
\begin{equation}\label{intoBandWp}
\xymatrix{
\tilde{\cal U}_-(-D_-)|_{2D_-} \ar@{->>}[r]^{\bar{h}|_{2D_-}} \ar@{=}[d] &
{\cal G}(-D_-) \ar[r] & {\cal W}_+|_{D_-} \ar[d]^{r_0} \\
\tilde{\cal U}_-(-D_-)|_{2D_-} \ar@{->>}[r]_-{h(-D_-)}
 & \tilde{\cal G}(-D_-) \ar[r]^p & B. }
\end{equation}
Then we can prove the right side of \eqref{BisWp} is commutative from
\eqref{twoquotientsofUm} and the surjectivity of $\bar{h}|_{2D_-}$.
\end{pf}\par
\begin{cor}\label{cor:WpisStables}
Let $t\in D_-$ be a closed point. Then the exact sequence
\[ 0 \longrightarrow \tilde{\cal G}(-D_-)\otimes k(t) \longrightarrow
{\cal W}_+\otimes k(t) \longrightarrow \tilde{\cal F}\otimes k(t)
\longrightarrow 0 \]
induced from the third column of \eqref{Elemtransm2} is nontrivial.
\end{cor}
\begin{pf}
Suppose not. Then Lemma \ref{lem:nontrivial} and Lemma \ref{lem:BisWp}
lead to a morphism $p_{l+1}: \Spec({\cal O}_{\tilQm}/ \tilde{m}_t^{l+1})
\rightarrow V_-$ such that \eqref{Pl+1} becomes commutative.
This $p_{l+1}$ induces a $\tilQm$-morphism
$q_{l+1}: \Spec({\cal O}_{\tilQm}/ \tilde{m}^{l+1}_t) \rightarrow
\tilQm \sideset{}{_{\Qm}}\times V_- = D_-$.
Thus $I_{D_-}$ is contained in $\tilde{m}_t^{l+1}$, which
contradicts the choice of $l$ in Lemma \ref{lem:nontrivial}.
\end{pf}\par
From the corollary above one can show that ${\cal W}_+\otimes k(t)\in
\Coh(X_{k(t)})$ is $a_+$-semistable for every point $t\in\tilQm$
in a similar fashion to the proof of Lemma \ref{lem:HNF} (ii).
This sheaf ${\cal W}_+$ accordingly gives a morphism $\tilde{\varphi}_+:
\tilQm \rightarrow \modp$.
Now we intend to construct a morphism $\bar{\phi}_+ : \tilmodm \rightarrow
\modp$ such that $\bar{\phi}_+\circ \tilde{\pi}_- : \tilQm \rightarrow
\modp$ is equal to $\tilde{\varphi}_+$.
\begin{lem}
The natural morphism $\tilQm \rightarrow \Qm \sideset{}{_{\modm}}\times
 \tilmodm$ is isomorphic.
\end{lem}
\begin{pf}
$\tilQms$ denotes the open subset
$(\phi_-\circ \tilde{\pi}_-)^{-1}(\modms)$ of $\tilQm$, and
$\tilmodms$ denotes $\phi_-^{-1}(\modms)$.
Because $E_-\subset\tilmodm$ is
contained in $\tilmodms$ it suffices to show that
$\tilQms \rightarrow \Qms \sideset{}{_{\modms}}\otimes \tilmodms$
is isomorphic.
Since $\pi_-: \Qms \rightarrow \modms$ is flat one can show that
$\pi_-^*( I_{P_-,\modms})=I_{V_-,\Qms}$, and hence that
$\pi_-^*( I_{P_-,\modms}^n)=I_{V_-,\Qms}^n$
for any $n$.
\end{pf}\par
Using this lemma one can induce an action
$\bar{\Sigma}_-: \bar{G}\times\tilQm = (\bar{G}\times\Qm) 
\sideset{}{_{\modm}}\times \tilmodm \rightarrow \tilQm=
\Qm \sideset{}{_{\modm}}\times \tilmodm$ from the action
$\bar{\sigma}_-: \bar{G}\times\Qm \rightarrow \Qm$.
\begin{lem}\label{lem:GactQm}
As to the morphism $\tilde{\varphi}_+$, the following is commutative:
\begin{equation*}
\xymatrix{
\bar{G}\times\tilQm \ar[r]_{\bar{\Sigma}-} \ar[d]^{\pr_2} &
\tilQm \ar[d]^{\tilde{\varphi}_+} \\
\tilQm \ar[r]^{\tilde{\varphi}_+} & \modp}
\end{equation*}
\end{lem}
%
%
%
One can prove this lemma easily.
$\pi_-: \Qm\rightarrow\modm$ is a good quotient by $\bar{\sigma}_-$,
so \cite[Page 8, Remark 5]{Mu:GIT} and \cite[Page 27, Theorem 1]{Mu:GIT}
imply that
$\tilde{\pi}_-: \tilQm=\Qm \sideset{}{_{\modm}}\times \tilmodm \rightarrow
\tilmodm$ is a categorical quotient by $\bar{\Sigma}_-$.
Therefore there is a unique morphism $\bar{\phi}_+: \tilmodm \rightarrow
\modp$ such that $\bar{\phi}_+\circ \tilde{\pi}_-: \tilQm \rightarrow
\modp$ is equal to $\tilde{\phi}_+$ because of the lemma above.

Consequently we can connect $\modm=\Modm$ with $\modp=\Modp$ by
\begin{equation}\label{flipm}
\xymatrix{
V_- \ar@{^{(}-}[r] & \Qm \ar[d]_{\pi_-} & \tilQm \ar[l]_{\varphi_-}
\ar[d]_{\tilde{\pi}_-} \ar[rd]^{\tilde{\varphi}_+} & \\
P_- \ar@{^{(}-}[r] & \modm & \tilmodm \ar[l]_{\phi_-} \ar[r]_{\bar{\phi}_+}
& \modp }
\end{equation}
when $P_-\subset\modm$ is nowhere dense. (Without this hypothesis $\tilmodm$
may be empty.)

We shall reverse $\modm$ and $\modp$ and follow a similar argument.
Let $V_+$ be a closed subscheme $\sideset{}{_{\ff\in A^+(a)}}\coprod Q^{-\ff}$,
and $P_+$ a closed subscheme $\sideset{}{_{\ff\in A^+(a)}}\coprod P^{-\ff}$ of
$\modp$, mentioned right after Definition \ref{defn:-ff}.
Let $\varphi_+: \tilQp \rightarrow \Qp$ be the blowing-up along $V_+$,
and $\phi_+:\tilmodp\rightarrow\modp$ the blowing-up along $P_+$.
Denote their exceptional divisors by $D_+\subset\tilQp$ and $E_+\subset
\tilmodp$ respectively.
Then we can construct a morphism $\bar{\varphi}_-: \tilQp\rightarrow\modm$
and make it descend to a morphism $\bar{\phi}_-: \tilmodp\rightarrow\modm$.
Thereby we get another sequence of morphisms connecting $\modm$ and
$\modp$ as follows:
\begin{equation}\label{flipp}
\xymatrix{
 & \tilQp \ar[dl]_{\tilde{\varphi}_-} \ar[d]^{\tilde{\pi}_+}
\ar[r]_{\varphi_+} & \Qp \ar[d]^{\pi_+} & V_+ \ar@{_{(}-}[l] \\
\modm & \tilmodp \ar[l]^{\bar{\phi}_-} \ar[r]_{\phi_+} & \modp &
P_+ \ar@{_{(}-}[l] }
\end{equation}
%
%
\section{$\bar{\phi}_+: \tilmodm\rightarrow\modp$ is blowing-up}
\label{ss:blowing-up}
In this section we would like to compare \eqref{flipm} with \eqref{flipp}
assuming that $P_-$ is nowhere dense in $\modm$.
The following lemma shall be needed later.
\begin{lem}\label{lem:locallift}
Let ${\cal U}_+$ be a universal quotient sheaf of $\Qp$ on $X_{\Qp}$,
and ${\cal W}_+$ the $X_{\tilQm}$-module defined at \eqref{elemtransm}.
There are an open covering $\sideset{}{_{\alp}}\bigcup U_{\alp}$ of
$\tilQm$, a morphism $\bar{\varphi}_+^{\alp}: U_{\alp}\rightarrow\Qp$
such that 
\begin{equation*}
\xymatrix{
U_{\alp} \ar[r]^{\bar{\varphi}_+^{\alp}} \ar@{_{(}-}[d] & \Qp
\ar[d]^{\pi_+} \\
\tilQm \ar[r]^{\tilde{\varphi}_+} & \modp}
\end{equation*}
is commutative, and an isomorphism
$\Phi_+^{\alp}: {\cal W}_+|_{U_{\alp}} \rightarrow (\bar{\varphi}_+^{\alp})^*
{\cal U}_+$ of $X_{U_{\alp}}$-modules.
Furthermore, we can assume that $U_{\alp}\cap U_{\beta} \subset \tilQms$ if
$\alp\neq\beta$.
\end{lem}
\begin{pf}
The proof of the first part is easy, so may be left to the
reader. Recall that both $\Qp$ and $\Qm$ are open subsets of a Quot-scheme
$Q$, and that ${\cal U}_+|_{\Qm\cap\Qp} = {\cal U}_-|_{\Qm\cap\Qp}$.
$U_0= \tilQm\setminus D_- = \Qm\setminus V_-$ is an open neighborhood
of $\tilQm\setminus\tilQms$, and is contained in $\Qm\cap\Qp$.
Let $\bar{\varphi}_+^0 : U_0 = \Qm\setminus V_- \rightarrow \Qp$ be
a natural open immersion, and $\Phi_+^0: {\cal W}_+|_{U_0}\rightarrow
{\cal U}_+|_{U_0}$ an isomorphism ${\cal W}_+|_{\tilQm\setminus D_-}
\rightarrow \tilde{\cal U}_-|_{\tilQm\setminus D_-} =
{\cal U}_-|_{\Qm\setminus V_-}= {\cal U}_+|_{\Qm\setminus V_-}$
induced from \eqref{elemtransm}. Then $\bar{\varphi}_+^0$ and
$\Phi_+^0$ satisfy the conditions in this lemma.
Thus we can assume that $U_{\alp}\cap U_{\beta} \subset \tilQms$ if
$\alp\neq\beta$.
\end{pf}\par
%
\begin{lem}\label{lem:inverseofPpisDm}
$\tilde{\varphi}_+^{-1}(P_+)$ is equal to $D_-= \varphi_-^{-1}(V_-)$
as closed subschemes in $\tilQm$.
\end{lem}
\begin{pf}
Clearly $D_-\subset \tilde{\varphi}_+^{-1}(P_+)$ from the construction
of $\tilde{\varphi}_+$.
We first consider the case where $D_-' := \tilde{\varphi}_+^{-1}(P_+)$
is a Cartier divisor of $\tilQm$.
By virtue of the definition of $V_+$, there is an exact sequence
\begin{equation}\label{HNFp}
0 \longrightarrow {\cal G}' \longrightarrow {\cal U}_+|_{X_{V_+}}
\longrightarrow {\cal F}' \longrightarrow 0
\end{equation}
of $V_+$-flat $X_{V_+}$-modules such that, for every closed point
$t$ of $Q^{-\ff}$, $(2c_1({\cal G}'_{k(t)})-c_1, c_2({\cal G}'_{k(t)}),
c_2({\cal F}'_{k(t)}))$ is equal to $-\ff$.
Similarly to Lemma \ref{lem:torsionfree}, ${\cal F}'$ and ${\cal G}'$
are flat family of torsion-free sheaves.
Pulling back this by $\bar{\varphi}_+^{\alp}$ of Lemma \ref{lem:locallift},
we have an exact sequence
\begin{equation}\label{localHNFp}
0 \longrightarrow (\bar{\varphi}_+^{\alp})^*\,{\cal G}'={\cal G}_{\alp}'
\longrightarrow (\bar{\varphi}_+^{\alp})^*\, {\cal U}_+|_{D_-'\cap U_{\alp}}
= \bar{\cal U}_+^{\alp}|_{D_-'\cap U_{\alp}} \longrightarrow
(\bar{\varphi}_+^{\alp})^*\, {\cal F}'= {\cal F}_{\alp}' \longrightarrow 0
\end{equation}
on $X\times (\bar{\varphi}_+^{\alp})^{-1}(V_+)= X_{D_-'\cap U_{\alp}}$,
where we put $(\bar{\varphi}_+^{\alp})^*{\cal U}_+=\bar{\cal U}_+^{\alp}$.
Let ${\cal V}_-$ denote $\Ker( \bar{\cal U}_+^{\alp} \twoheadrightarrow
\bar{\cal U}_+^{\alp}|_{D_-'\cap U_{\alp}} \twoheadrightarrow 
{\cal F}_{\alp}')$, that is,
\begin{equation}\label{localelemtransp}
0 \longrightarrow {\cal V}_- \longrightarrow \bar{\cal U}_+^{\alp}
\longrightarrow {\cal F}'_{\alp} \longrightarrow 0
\end{equation}
is exact. ${\cal V}_-$ is flat over $U_{\alp}$ since
$D'_-$ is a Cartier divisor of $\tilQm$.

Because $D_-'\supset D_-$, the isomorphism $\Phi_+^{\alp}$ in Lemma
\ref{lem:locallift} induces a surjection
$\bar{\cal U}_+^{\alp}|_{D'_-\cap U_{\alp}} \twoheadrightarrow
{\cal W}_+|_{D_-\cap U_{\alp}}$. Hence we have a diagram on $X_{D'_-\cap
U_{\alp}}$ 
\begin{equation}\label{UptoWp}
\xymatrix{
0 \ar[r] & {\cal G}'_{\alp} \ar[r] & \bar{\cal U}_+^{\alp}|_{D'_-\cap
U_{\alp}} \ar[r] \ar@{->>}[d] & {\cal F}'_{\alp} \ar[r] \ar@{.>}[d]^r &
0 \\
0 \ar[r] & \tilde{\cal G}(-D_-)|_{U_{\alp}} \ar[r] & {\cal W}_+|_{D_-
\cap U_{\alp}} \ar[r] & \tilde{\cal F}|_{U_{\alp}} \ar[r] & 0, }
\end{equation}
where the first row is \eqref{localHNFp} and the second row is the
restriction of the third column in \eqref{Elemtransm2} to
$X_{D_-\cap U_{\alp}}$.
One can check that
$Hom_{X_{D'_-\cap U_{\alp}}/ D'_-\cap U_{\alp}}({\cal G}'_{\alp},
\tilde{\cal F}|_{U_{\alp}})= 
Hom_{X_{D_-\cap U_{\alp}}/ D_-\cap U_{\alp}}({\cal G}'_{\alp}|_{D_-\cap
U_{\alp}}, \tilde{\cal F}|_{U_{\alp}})=0$ by base change theorem on
relative Ext sheaves, and so one can find $r: {\cal F}'_{\alp} \rightarrow
\tilde{\cal F}|_{U_{\alp}}$ such that \eqref{UptoWp} is commutative.
Then the following also is commutative:
\begin{equation}\label{VmtoUm}
\xymatrix{
0 \ar[r] & {\cal V}_- \ar[r] \ar@{.>}[d]^s & \bar{\cal U}_+^{\alp}
\ar[r] \ar[d]^{(\Phi_+^{\alp})^{-1}} & {\cal F}'_{\alp} \ar[r]
\ar[d]^r & 0 \\
0 \ar[r] & \tilde{\cal U}_-(-D_-)|_{U_{\alp}} \ar[r] & 
{\cal W}_+|_{U_{\alp}} \ar[r] & \tilde{\cal F}|_{U_{\alp}} \ar[r] &
0, }
\end{equation}
where the first row is \eqref{localelemtransp} and the second row
is the restriction of the second column in \eqref{Elemtransm2}
to $X_{U_{\alp}}$.
\begin{clm}\label{clm:assetDisDd}
Set-theoretically, $D_-\cap U_{\alp}$ coincides with $D'_-\cap U_{\alp}$.
\end{clm}
\begin{pf}
Suppose not. Then one can find a closed point $t\in D'_-$ that is not
contained in $D_-$.
Since $t\in D'_-$, \eqref{localHNFp} implies that $\bar{\cal U}_+^{\alp}
\otimes k(t)$ is not $a_-$-semistable. Since $t \not\in D_-$,
\eqref{elemtransm} implies that ${\cal W}_+\otimes k(t)$ is 
$a_-$-semistable.
This is a contradiction because $\bar{\cal U}_+^{\alp}\otimes k(t)$ is
isomorphic to ${\cal W}_+\otimes k(t)$.
\end{pf}\par
One can obtain the following commutative diagram by tensoring
${\cal O}_{D'_-\cap U_{\alp}}$ to the first row in \eqref{VmtoUm}
and ${\cal O}_{D_-\cap U_{\alp}}$ to the second row in \eqref{VmtoUm}
since $D'_-\supset D_-$:
\begin{equation}\label{VmtoUmD}
\xymatrix{
0 \ar[r] & {\cal F}'_{\alp}(-D'_-) \ar[r] \ar[d]^u & {\cal V}_-|_{D'_-
\cap U_{\alp}} \ar[r] \ar[d]^{s'} & \bar{\cal U}_+^{\alp}|_{D'_-\cap
U_{\alp}} \ar[r] \ar@{->>}[d]^{(\Phi^{-1})'} & {\cal F}'_{\alp} \ar[r]
\ar[d]^r & 0 \\
0 \ar[r] & \tilde{\cal F}(-D_-)|_{U_{\alp}} \ar[r] & 
\tilde{\cal U}(-D_-)|_{D_-\cap U_{\alp}} \ar[r] & {\cal W}_+|_{D_-\cap
U_{\alp}} \ar[r] & \tilde{\cal F}|_{U_{\alp}} \ar[r] & 0 }
\end{equation}
\begin{clm}
$s\otimes k(t): {\cal V}_-\otimes
k(t) \rightarrow \tilde{\cal U}_-(-D_-)\otimes k(t)$ in \eqref{VmtoUm}
is isomorphic for every closed point $t\in U_{\alp}$.
\end{clm}
\begin{pf}
We have to verify this only in case where $t$ is contained in $D'_-$.
By Claim \ref{clm:assetDisDd} $t$ is also contained in $D_-$.
Tensoring $k(t)$ to \eqref{VmtoUmD}, we obtain a commutative diagram
\begin{equation}\label{VmtoUmFibr}
\xymatrix{
0 \ar[r] & {\cal F}'_{\alp}(-D'_-)_{k(t)} \ar[r] \ar[d]^{u_t} &
{\cal V}_{-\: k(t)} \ar[r] \ar[d]^{s'_t} & \bar{\cal U}^{\alp}_{+\: k(t)}
\ar[r] \ar[d]^{(\Phi^{-1})'_t} & {\cal F}'_{\alp \: k(t)} \ar[r]
\ar[d]^{r_t} & 0 \\
0 \ar[r] & \tilde{\cal F}(-D_-)_{k(t)} \ar[r] & \tilde{\cal U}_-(-D_-)_{k(t)}
\ar[r] & {\cal W}_{+\: k(t)} \ar[r] & \tilde{\cal F}_{k(t)} \ar[r] &
0 }
\end{equation}
whose rows are exact.
$(\Phi^{-1})'_t$ is isomorphic by its definition. One can see that
also $r_t$ is isomorphic by the uniqueness of the Harder-Narasimhan
filtration with respect to $a_-$-stability. Thus $s'_t$ is nonzero map.
If $u_t$ is zero map, then $s'_t$ induces a nonzero homomorphism
\[\bar{s}'_t: {\cal G}'_{\alp}\otimes k(t) = \Cok( {\cal F}'_{\alp}
(-D'_-)_{k(t)} \rightarrow {\cal V}_{-\: k(t)}) \rightarrow
\tilde{\cal U}_-(-D_-)_{k(t)}\]
by \eqref{localHNFp}.
This $\bar{s}'_t$ should be injective because ${\cal G}'_{\alpha}\otimes
k(t)$ is torsion-free and rank-one.
This contradicts the $a_-$-semistability of 
$\tilde{\cal U}_-(-D_-)_{k(t)}$, and so $u_t$ should be nonzero, and hence
injective. Then one can see $s'_t$ is injective by diagram-chasing.
\eqref{VmtoUmFibr} implies the Chern classes of ${\cal V}_{-\: k(t)}$
are equal to those of $\tilde{\cal U}_-(-D_-)_{k(t)}$, we see that $s'_t=
s\otimes k(t)$ is isomorphic.
\end{pf}\par
Both ${\cal V}_-$ and $\tilde{\cal U}_-(-D_-)|_{U_{\alp}}$ are 
$U_{\alp}$-flat, and hence the claim above implies that $s$ in
\eqref{VmtoUm} is isomorphic. Then also $r$ in \eqref{VmtoUm} is
isomorphic. Because ${\cal F}'_{\alp}$ is $D'_-\cap U_{\alp}$-flat
and $\tilde{\cal F}|_{U_{\alp}}$ is $D_-\cap U_{\alp}$-flat, one can
verify that $D_-\cap U_{\alp}$ is equal to $D'_-\cap U_{\alp}$.
Since this holds good for every $U_{\alp}$, we conclude the proof of
this lemma in case where $\tilde{\varphi}_+^{-1}(P_+)$ is a Cartier
divisor.

Next, we consider the case where $\tilde{\varphi}_+^{-1}(P_+)=D'_-$
is not necessarily a Cartier divisor of $\tilQm$.
Let $\varphi_-^{(2)}: \tilde{Q}_-^{(2)} \rightarrow \tilQm$ be the blowing
up along $D'_-$. Let $\tilde{D}_-$ and $\tilde{D'}_-$ denote closed
subschemes $(\varphi_-^{(2)})^{-1}(D_-)$ and $(\varphi_-^{(2)})^{-1}(D'_-)$
of $\tilde{Q}_-^{(2)}$, respectively.
For a natural exact sequence
\[ 0 \longrightarrow {\cal O}_{\tilQm}(-D_-) \longrightarrow 
{\cal O}_{\tilQm} \rightarrow {\cal O}_{D_-} \longrightarrow 0 \]
on $\tilQm$, one can verify that also its pull-back by $\varphi_-^{(2)}$
\[ 0 \longrightarrow (\varphi_-^{(2)})^* {\cal O}_{\tilQm}(-D_-)
\longrightarrow {\cal O}_{\tilde{Q}_-^{(2)}} \longrightarrow
{\cal O}_{\tilde{D}_-} \longrightarrow 0 \]
is exact. 
In view of this, one can check that the pull-back of \eqref{Elemtransm2}
by $\id_X\times \varphi_-^{(2)}$
\begin{equation*}
\xymatrix{
 & & 0 \ar[d] & 0 \ar[d] & \\
0 \ar[r] & \tilde{\cal W}_+^{(2)}(-\tilde{D}_-) \ar[r] \ar@{=}[d] &
\tilde{\cal U}_-^{(2)}(-\tilde{D}_-) \ar[r] \ar[d] &
\tilde{\cal G}^{(2)}(-\tilde{D}_-) \ar[r] \ar[d] & 0 \\
0 \ar[r] & \tilde{\cal W}_+^{(2)}(-\tilde{D}_-) \ar[r] &
\tilde{\cal W}_+^{(2)} \ar[r] \ar[d] & \tilde{\cal W}_+^{(2)}|_{\tilde{D}_-}
\ar[r] \ar[d] & 0 \\
 & & \tilde{\cal F}^{(2)} \ar@{=}[r] \ar[d] & \tilde{\cal F}^{(2)} \ar[d]
 & \\
 & & 0 & 0 & }
\end{equation*}
satisfies that its rows and columns are exact, where $\tilde{\cal W}_+^{(2)}$
denotes $(\id_X\times \tilde{\varphi}_-^{(2)})^* {\cal W}_+$, and so on.
Now both $\tilde{D}'_-$ and $\tilde{D}_-$ are Cartier divisors, and we
can show that $\tilde{D}'_- =\tilde{D}_-$ as subschemes of $\tilde{Q}_-^{(2)}$
in the same way as the proof in the preceding case.
\begin{clm}\label{noetherian}
Let $R$ be a Noetherian ring, $t$ an element of $R$ which is not a
zero-divisor, and $tR\supset I$ an ideal of $R$.
Suppose that $\Proj_R(\oplus\,I^n / I^{n+1})=
\Proj_R(\oplus\, I^n/ tI^n)$ as subschemes in $\Proj_R(\oplus_{n\geq 0}\,
I^n)$.
Then $tR=I$ if $\Spec(R/I)$ is nowhere dense in $\Spec(R)$.
\end{clm}
Its proof is left to the reader. Now Lemma \ref{lem:inverseofPpisDm}
is immediate from Claim \ref{clm:assetDisDd} and Claim \ref{noetherian}.
\end{pf}\par
\begin{cor}\label{cor:inverseofPpisEm}
$(\bar{\phi}_+)^{-1}(P_+)$ coincides with $E_-= \phi_-^{-1}(P_-)$
as subschemes of $\tilmodm$.
\end{cor}
\begin{pf}
By Claim \ref{clm:assetDisDd}, closed subschemes $E_-$ and
$(\bar{\phi}_+)^{-1}(P_+)$ of $\tilmodm$ are contained in $\tilmodms$.
Thus $\tilde{\pi}_-: \tilde{\pi}_-^{-1}(E_-)=D_- \rightarrow E_-$ and
$\tilde{\pi}_-:
\tilde{\pi}_-^{-1} \bar{\phi}_+^{-1}(P_+)=\tilde{\varphi}_+^{-1}(P_+)
\rightarrow \bar{\phi}_+^{-1}(P_+)$ are faithfully-flat.
Hence this corollary is immediate from Lemma \ref{lem:inverseofPpisDm}.
\end{pf}\par
By the corollary above, there is a morphism $\Delta_+: \tilmodm \rightarrow
\tilmodp$ such that
$\phi_+\circ \Delta_+: \tilmodm \rightarrow \tilmodp \rightarrow
\modp$ is equal to $\bar{\phi}_+$.
Likewise, for $U_{\alp}\subset\tilQm$ and $\bar{\varphi}_+^{\alp}$ in
Lemma \ref{lem:locallift}, there is a morphism $\Delta_+^{\alp}:
U_{\alp} \rightarrow \tilQp$ such that
$\varphi_+\circ\Delta_+^{\alp}: U_{\alp}\rightarrow \tilQp \rightarrow
\Qp$ is equal to $\bar{\varphi}_+^{\alp}$
since $(\bar{\varphi}_+^{\alp})^{-1}(V_+)=
(\tilde{\varphi}_+)^{-1}(P_+)\cap U_{\alp}$ is a Cartier divisor of
$U_{\alp}$ by Lemma \ref{lem:inverseofPpisDm}.
\begin{lem}\label{lem:Delp&Delpalp}
\begin{equation*}
\xymatrix{
 U_{\alp} \ar@{^{(}-}_{i_{\alp}}[r] \ar^{\Delta_+^{\alp}}[d] &
 \tilQm \ar_{\tilde{\pi}_-}[r] & \tilmodm \ar[d]_{\Delta_+} \\
 \tilQp \ar^{\tilde{\pi}_+}[rr] && \tilmodp }
\end{equation*}
is commutative.
\end{lem}
\begin{pf}
One can check that both $\phi_+\circ (\Delta_+\circ \tilde{\pi}_-\circ
i_{\alp}): U_{\alp}\rightarrow\tilmodp\rightarrow\modp$ and
$\phi_+\circ (\tilde{\pi}_+\circ \Delta_+^{\alp})$ coincide with
$\pi_+\circ \bar{\varphi}_+^{\alp}: U_{\alp}\rightarrow \Qp\rightarrow
\modp$.
Then this lemma follows by the universal property of the blowing-up
$\phi_+: \tilmodp\rightarrow\modp$.
\end{pf}\par
\begin{prop}
The morphism $\bar{\phi}_-\circ \Delta_+: \tilmodm\rightarrow\tilmodp
\rightarrow\modm$ is equal to $\phi_-: \tilmodm \rightarrow\modm$.
\end{prop}
\begin{pf}
First, let us verify the commutativity of
\begin{equation}\label{Delpalp&phim}
\xymatrix{
 U_{\alp} \ar@{^{(}-}[r]_{i_{\alp}} \ar[d]_{\Delta_+^{\alp}} &
 \tilQm \ar[r]_{\varphi_-} & \Qm \ar[d]^{\pi_-} \\
 \tilQp \ar[rr]^{\tilde{\varphi}_-} && \modm. }
\end{equation}
Pulling back an exact sequence \eqref{HNFp} on $X_{V_+}$ by
$\id_X\times\varphi_+: X_{\tilQp} \rightarrow X_{\Qp}$, we obtain a
commutative diagram on $X_{\tilQp}$
\begin{equation*}
\xymatrix{
 0 \ar[r] & {\cal W}_- \ar[r] \ar@{->>}[d] & \varphi_+^*{\cal U}_+=
 \tilde{\cal U}_+ \ar[r] \ar@{->>}[d] & \varphi_+^*{\cal F}'=
 \tilde{\cal F}' \ar[r] \ar@{=}[d] & 0 \\
 0 \ar[r] & \varphi_+^*{\cal G}'=\tilde{\cal G}' \ar[r] &
 \tilde{\cal U}_+|_{D_+} \ar[r] & \tilde{\cal F}' \ar[r] & 0}
\end{equation*}
whose rows are exact. Remark that ${\cal W}_-$ is $\tilQp$-flat.
Pulling back this diagram by $\id_X\times \Delta_+^{\alp}$, we obtain
a commutative diagram on $X_{U_{\alp}}$
\begin{equation}\label{elemtranspDelta}
\xymatrix{
 & (\Delta_+^{\alp})^*{\cal W}_- \ar[r] \ar[d] & (\Delta_+^{\alp})^*
 \tilde{\cal U}_+ \ar[r] \ar@{->>}[d] & (\Delta_+^{\alp})^*\tilde{\cal F}'
 \ar[r] \ar@{=}[d] & 0 \\
 0 \ar[r] & (\Delta_+^{\alp})^* \tilde{\cal G}' \ar[r] &
 (\Delta_+^{\alp})^* \tilde{\cal U}_+|_{D_-\cap U_{\alp}} \ar[r] &
 (\Delta_+^{\alp})^* \tilde{\cal F}' \ar[r] & 0 }
\end{equation}
whose rows are exact, because $(\Delta_+^{\alp})^{-1}
(D_+)=D_-\cap U_{\alp}$ by Lemma \ref{lem:inverseofPpisDm}.
Compare this with a commutative diagram 
\begin{equation}\label{partofElemtransm2}
\xymatrix{
 0 \ar[r] & \tilde{\cal U}_-(-D_-) \ar[r] \ar@{->>}[d] & {\cal W}_+ \ar[r]
 \ar@{->>}[d] & \tilde{\cal F} \ar[r] \ar@{=}[d] & 0 \\
 0 \ar[r] & \tilde{\cal G}(-D_-) \ar[r] & {\cal W}_+|_{D_-} \ar[r] &
 \tilde{\cal F} \ar[r] & 0 }
\end{equation}
on $X_{\tilQm}$ in \eqref{Elemtransm2}.
Since $(\Delta_+^{\alp})^*\tilde{\cal U}_+ = (\bar{\varphi}_+^{\alp})^*
{\cal U}_+$, an isomorphism $\Phi_+^{\alp}$ in Lemma \ref{lem:locallift}
connects the second row of \eqref{elemtranspDelta} with that of
\eqref{partofElemtransm2}:
\begin{equation}\label{UpandWp}
\xymatrix{
 0 \ar[r] & (\Delta_+^{\alp})^* \tilde{\cal G}' \ar[r] &
 (\Delta_+^{\alp})^* \tilde{\cal U}_+ |_{D_-\cap U_{\alp}} \ar[r] 
 \ar[d]_{(\Phi_+^{\alp})^{-1}}&
 (\Delta_+^{\alp})^* \tilde{\cal F}' \ar[r] \ar@{.>}[d]^{\gamma} \ar[r]
 & 0 \\
 0 \ar[r] & \tilde{\cal G}(-D_-)|_{U_{\alp}} \ar[r] & {\cal W}_+|_{D_-\cap
 U_{\alp}} \ar[r] & \tilde{\cal F}|_{U_{\alp}}
 \ar[r] & 0}
\end{equation}
Remark that all sheaves in this diagram are flat over $D_-\cap U_{\alp}$.
One can check that two exact sequences in this diagram are relative
Harder-Narasimhan filtrations of $(\Delta_+^{\alp})^* \tilde{\cal U}_+
 |_{D_-\cap U_{\alp}} \simeq {\cal W}_+|_{D_-\cap U_{\alp}}$ with respect
to $a_-$-stability, and hence we get a homomorphism
$\gamma: (\Delta_+^{\alp})^* \tilde{\cal F}' \rightarrow \tilde{\cal F}
|_{U_{\alp}}$ which makes \eqref{UpandWp} commutative.
$\gamma\otimes k(t)$ is isomorphic for any $t\in D_-\cap U_{\alp}$ because
of the uniqueness of HNF, and so $\gamma$ should be isomorphic.
\eqref{elemtranspDelta}, \eqref{partofElemtransm2}, $\Phi_+^{\alp}$ and
$\gamma$ induce a surjective homomorphism
\begin{equation}\label{WmandUm}
s: (\Delta_+^{\alp})^*{\cal W}_- \rightarrow \tilde{\cal U}_-(-D_-)
|_{U_{\alp}}.
\end{equation}
In fact $s\otimes k(t)$ should be isomorphic for any closed point $t\in
U_{\alp}$, since $(\Delta_+^{\alp})^*{\cal W}_-\otimes k(t)$ and
$\tilde{\cal U}_-(-D_-)\otimes k(t)$ has the same Chern classes.
Thereby \eqref{WmandUm} is isomorphic, and hence \eqref{Delpalp&phim}
is commutative.
From Lemma \ref{lem:Delp&Delpalp} and \eqref{Delpalp&phim}, one can
verify
$\bar{\phi}_-\circ \Delta_+\circ \tilde{\pi}_-\circ i_{\alp}:
U_{\alp}\hookrightarrow \tilQm \rightarrow \modm$ equals
$\phi_-\circ \bar{\pi}_-\circ i_{\alp}: U_{\alp} \rightarrow \tilQm
\rightarrow \modm$ by diagram-chasing.
Hence
$(\bar{\phi}_-\circ \Delta_+)\circ \tilde{\pi}_-: \tilQm \rightarrow
\tilmodm \rightarrow \modm$ equals
$\phi_-\circ\tilde{\pi}_-: \tilQm \rightarrow \tilmodm \rightarrow \modm$.
As mentioned in the preceding section, $\tilde{\pi}_-: \tilQm\rightarrow
\tilmodm$ is a categorical quotient by $\bar{G}$.
Therefore we conclude that $\bar{\phi}_-\circ \Delta_+: \tilmodm
\rightarrow \tilmodp \rightarrow \modm$ coincides with $\phi_-$, thanks
to the property of categorical quotients.
\end{pf}\par
From the proposition above we get a morphism $\Delta_+$ such that
\begin{equation}\label{Delp}
\xymatrix{
 \tilmodm \ar[r]^{\bar{\phi}_+} \ar[d]_{\phi_-} \ar[dr]^{\Delta_+} &
 \modp \\
 \modm & \tilmodp \ar[l]^{\bar{\phi}_-} \ar[u]_{\phi_+} }
\end{equation}
is commutative. Quite similarly, there is a morphism $\Delta_-: \modp
\rightarrow \modm$ such that
\begin{equation}\label{Delm}
\xymatrix{
 \tilmodm \ar[r]^{\bar{\phi}_+} \ar[d]_{\phi_-} & \modp \\
 \modm & \tilmodp \ar[l]^{\bar{\phi}_-} \ar[u]_{\phi_+} \ar[ul]_{\Delta_-}}
\end{equation}
is commutative.
Thus $\phi_-\circ (\Delta_-\circ \Delta_+): \tilmodm \rightarrow \tilmodm
\rightarrow \modm$ is equal to $\phi_-: \tilmodm \rightarrow \modm$, and
so $\Delta_-\circ \Delta_+ : \tilmodm \rightarrow \tilmodp \rightarrow
\tilmodm$ should be $\id_{\tilmodm}$ because of the universal property
of blowing-up $\phi_-$.
Likewise $\Delta_+\circ\Delta_-: \tilmodp\rightarrow\tilmodp$ equals
$\id_{\tilmodp}$, and hence both $\Delta_+$ and $\Delta_-$ are
isomorphic. Summarizing:
\begin{prop}
As to \eqref{flipm} and \eqref{flipp}, there are isomorphisms $\Delta_+:
\tilmodm \rightarrow \tilmodp$ and $\Delta_-: \tilmodp\rightarrow
\tilmodm$ such that \eqref{Delp} and \eqref{Delm} are commutative.
In particular, the morphism $\bar{\phi}_+: \tilmodm
\rightarrow\modp$, which is constructed by the method of elementary
transform and descent theory, is the blowing-up of $\modp$ along
$W_+$.
\end{prop}
%
%
%
\section{Some structure of $P^{\ff}$ over $\Pic(X)\times\Hilb(X)\times
\Hilb(X)$}\label{ss:overPicHilb}
Let $\ff=(f,m,n)$ be a member of $A^+(a)$. \eqref{HNFm} gives an exact
sequence
\begin{equation}\label{HNFmf}
0 \longrightarrow {\cal F} \longrightarrow {\cal U}_-|_{\Qf} 
\longrightarrow {\cal G} \longrightarrow 0
\end{equation}
of $\Qf$-flat ${\cal O}_{X_{\Qf}}$-modules.
By Lemma \ref{lem:torsionfree} both ${\cal F}\otimes k(t)$ and
${\cal G}\otimes k(t)$ are torsion-free, rank-one, and hence $H$-stable
for any $t\in\Qf$.
Denote by $M_H(1,F,m)$ the coarse moduli scheme of $H$-stable rank-one
sheaves on $X$ with Chern classes $(F,m)\in\Num(X)\times\ZZ$.
Then ${\cal F}$ and ${\cal G}$ in \eqref{HNFmf} induce morphisms
$\tau_{\cal F}: \Qf \rightarrow M_H(1,(c_1+f)/2,m)$ and
$\tau_{\cal G}: \Qf \rightarrow M_H(1,(c_1-f)/2,n)$.
On the other hand $M_H(1,F,m)$ is isomorphic to $\Pic^F(X)\times
\Hilb^m(X)$, where $\Pic^F(X)$ is an open closed
subscheme $\bigl\{ L\in\Pic(X) \bigm| [L]=F \text{ in } \Num(X) \bigr\}$
of $\Pic(X)$.
Thereby, using $\tau_{\cal F}$ and $\tau_{\cal G}$ we obtain
a morphism
$\tau^Q:\Qf\rightarrow \Pic^{(c_1+f)/2}(X)\times \Hilb^m(X)\times
\Hilb^n(X)$ which has the following properties:
Let ${\cal P}\in\Coh(X_{\Pic})$ be a universal line bundle of $\Pic(X)$,
and let $I_{Z_1}\in\Coh (X_{\Hilb^m})$ (resp. $I_{Z_2}\in\Coh(X_{\Hilb^n})$)
be the ideal sheaf of a universal sheaf of $\Hilb^m(X)$ (resp. $\Hilb^n(X)$).
Define ${\cal F}_0$ and ${\cal G}_0\in\Coh(X_{\Pic\times \Hilb^m\times
\Hilb^n})$ by
\begin{equation}\label{F0G0}
{\cal F}_0 := \pr_{12}^*({\cal P})\otimes \pr_{13}^*(I_{Z_1})
\text{ and }
{\cal G}_0 := c_1\otimes \pr_{12}^*({\cal P}^{\vee})\otimes
\pr_{14}^*(I_{Z_2}).
\end{equation}
Then one can find line bundles $L_1$ and $L_2$ on $\Qf$ such that
\begin{equation}\label{FandF0}
{\cal F}\simeq (\tau^Q)^*{\cal F}_0 \otimes L_1 \text{ and }
{\cal G}\simeq (\tau^Q)^*{\cal G}_0 \otimes L_2.
\end{equation}
From now on, we shorten $\Pic^{(c_1+f)/2}(X)\times \Hilb^m(X)\times
\Hilb^n(X)$ to $T=T^{\ff}$.

One can show that $\tau^Q: \Qf\rightarrow T$ is $\bar{G}$-invariant in
a similar fashion to the proof of Lemma \ref{lem:GactQm}, and hence
$\tau^Q$ descends to a morphism $\tau_-: P^{\ff}\rightarrow T$, since
$\pi_-: \pi_-^{-1}(P^{\ff})=\Qf \rightarrow P^{\ff}$ is a categorical
quotient by $\bar{G}$.
In this section we would like to study some structure of $P^{\ff}$ as a
$T$-scheme.

One can find bounded complexes
$F^{\bullet}$ and $G^{\bullet}$ of locally-free ${\cal O}_T$-modules
of finite rank which allow quasi-isomorphisms $\tau_F: F^{\bullet}
\rightarrow {\cal F}_0$ and $\tau_G: G^{\bullet} \rightarrow {\cal G}_0$
of complexes. Let $q: X_T\rightarrow T$ be the projection.
The Serre duality \cite{Ha:residue} asserts a natural homomorphism
\begin{multline}\label{SerreDuality}
\Theta_q: \bR q_*\,  \bR Hom_{X_T}\left( Hom_{X_T}(F^{\bullet}, G^{\bullet}),
{\cal O}_{X_T}[2] \right)\rightarrow \\
\bR Hom_T\left( \bR q_*( Hom_{X_T}(F^{\bullet},G^{\bullet}(K_X)), {\cal O}_T
\right)
\end{multline}
in the derived category $D(T)$ is isomorphism.
Now we shall deduce the following from this.
\begin{prop}\label{prop:SerreDuality}
For any $T$-scheme $f: S\rightarrow T$, there is an isomorphism
\[ \Theta_{f^*q}: Ext^1_{X_S/S}(f^*{\cal G}_0, f^*{\cal F}_0) \rightarrow
   Ext^1_{X_S/S}(f^*{\cal F}_0, f^*{\cal G}_0(K_X))^{\vee} \]
of relative Ext sheaves.
\end{prop}
\begin{pf}
We prove this lemma only in case where $S=T$. It's easy to extend the proof
to general case.
As to the left side of \eqref{SerreDuality}, one can check that
\begin{equation}\label{RqRhom}
\left[ \bR q_*\, \bR Hom_{X_T}\left( Hom_{X_T}(F^{\bullet}, G^{\bullet}),
{\cal O}_{X_T}[2] \right) \right]_{-l} \simeq 
Ext^{2-l}_{X_T/T}({\cal G}_0, {\cal F}_0)
\end{equation}
for any integer $l$.
%
%
%
%
Now consider the right side of \eqref{SerreDuality}.
If we fix an affine open covering ${\bold U}=\{ U_i\}_i$ of $X_T$ such that
$q:U_i\hookrightarrow X_T \rightarrow T$ is affine, then we can construct
a quasi-isomorphism
\begin{equation*}
Hom_{X_T}(F^{\bullet}, G^{\bullet}(K_X))\longrightarrow Hom_{X_T}(F^{\bullet},
{\cal G}_0(K_X)) \longrightarrow 
{\cal C}^{\bullet}\left( Hom_{X_T}(F^{\bullet},
{\cal G}_0(K_X)), {\bold U} \right)
\end{equation*}
to the C\v{e}ch complex similarly to \cite[Lemma III.4.2]{Ha:text}. 
\[q_*\left( {\cal C}^{\bullet}( Hom_{X_T}(F^{\bullet},
{\cal G}_0(K_X)), {\bold U}) \right)\]
represents
$\bR q_*( Hom_{X_T}(F^{\bullet}, G^{\bullet}(K_X)))$
since
${\cal C}^p \left( Hom_{X_T}( F^q, {\cal G}_0(K_X) ), {\bold U}\right)$ is
$q_*$-acyclic.
Therefore, for an injective resolution $\iota_T: {\cal O}_T \rightarrow
K^{\bullet}$, a complex
\[Hom_T( q_*( {\cal C}^{\bullet}( Hom_{X_T}(F^{\bullet}, {\cal G}_0(K_X)),
{\bold U})), K^{\bullet}) \]
represents
$\bR Hom_T \left( \bR q_*( Hom_{X_T}(F^{\bullet}, G^{\bullet}(K_X))),
{\cal O}_T \right)$.
%
%
%
Furthermore, for any affine open subset $T_{\alp}$ of $T$, there is a
bounded complex $H_{\alp}^{\bullet}$
of free ${\cal O}_{T_{\alp}}$-modules of finite rank and with a
quasi-isomorphism
\begin{equation}\label{perfectCplx}
h_{\alp}: H_{\alp}^{\bullet} \rightarrow q_* {\cal C}^{\bullet}(
Hom_{X_T}(F^{\bullet}, {\cal G}_0(K_X)), {\bold U})|_{T_{\alp}}.
\end{equation}
by \cite[Page 47, Lemma 1.1]{Mu:Abelian}.
This $h_{\alp}$ and $\iota_T: {\cal O}_T \rightarrow K^{\bullet}$ give
rise to an isomorphism
\begin{multline}\label{CechtoPerfect}
\left[ Hom_T ( q_* {\cal C}^{\bullet}( Hom_{X_T}(F^{\bullet}, 
{\cal G}_0(K_X)), {\bold U}), K^{\bullet})|_{T_{\alp}} \right]_{-1} 
\simeq \\
\left[ Hom_{T_{\alp}}(H_{\alp}^{\bullet}, K^{\bullet}) \right]_{-1}\simeq
\left[ Hom_{T_{\alp}}(H_{\alp}^{\bullet}, {\cal O}_{T_{\alp}}) \right]_{-1}.
\end{multline}
\begin{clm}
This complex $H_{\alp}^{\bullet}$ induces an isomorphism
\begin{equation*}
i_{\alp}: Hom_{T_{\alp}}([H_{\alp}^{\bullet}]_1, {\cal O}_{T_{\alp}})
\rightarrow \left[Hom_{T_{\alp}}( H_{\alp}^{\bullet}, {\cal O}_{T_{\alp}}) 
\right]_{-1}.
\end{equation*}
\end{clm}
\begin{pf}
As a result of the base change theorem for relative Ext sheaves
\cite[Theorem 1.4]{La:extfamily}, $Ext^2_{X_T/T}({\cal F}_0, {\cal G}_0
(K_X))$ is equal to zero. Thus one can assume that $H_{\alp}^l=0$ if
$l\geq 2$. The remaining part of the proof is easy and left to the reader.
\end{pf}\par
From \eqref{perfectCplx}, \eqref{CechtoPerfect} and the claim above,
we obtain an isomorphism
\begin{multline*}
j_{\alp}: \left[ Hom_T( q_* {\cal C}^{\bullet}( Hom_{X_T}(F^{\bullet},
{\cal G}_0(K_X)), {\bold U}), K^{\bullet})\right]_{-1}|_{T_{\alp}} 
\rightarrow \\
Hom_T\left( \left[q_*{\cal C}^{\bullet}( Hom_{X_T}(F^{\bullet}, {\cal G}_0(K_X)
),
{\bold U}) \right]_1, {\cal O}_T \right)|_{T_{\alp}}
\end{multline*}
\begin{clm}
Let $T_{\alp}$ and $T_{\beta}$ be affine open subsets in $T$.
Then $j_{\alp}|_{T_{\alpha\beta}}= j_{\beta}|_{T_{\alpha\beta}}$.
\end{clm}
\begin{pf}
For $h_{\alp}$ and $h_{\beta}$ at \eqref{perfectCplx}, there are a bounded
complex $K_{\alp\beta}^{\bullet}$ of locally free 
${\cal O}_{T_{\alp\beta}}$-modules of finite rank, and quasi-isomorphisms
$k_{\alp}$ and $k_{\beta}$ such that
\begin{equation*}
\xymatrix{
 K^{\bullet}_{\alp\beta} \ar[r]_{k_{\beta}} \ar[d]_{k_{\alp}} &
 H^{\bullet}_{\beta} \ar[d]_{h_{\beta}|_{T_{\alp\beta}}} \\
 H^{\bullet}_{\alp} \ar[r]^-{h_{\alp}|_{T_{\alp\beta}}} &
 q_* {\cal C}^{\bullet} (Hom_{X_T}(F^{\bullet}, {\cal G}_0(K_X)),
 {\bold U})|_{T_{\alp\beta}} }
\end{equation*}
is commutative up to homotopy.
%
This $(K^{\bullet}_{\alp\beta}, k_{\alp}, k_{\beta})$ can be found by
using \cite[Page 47, Lemma 1.1]{Mu:Abelian} and the mapping cone
complex $Z^{\bullet}(f)$ \cite[Page 26]{Ha:residue}.
Then a quasi-isomorphism
$k_{\alp}|_{T_{\alp\beta}}\circ k_{\alp}: K_{\alp\beta}^{\bullet}
\rightarrow q_* {\cal C}^{\bullet}( Hom_{X_T}(F^{\bullet},
{\cal G}_0(K_X)), {\bold U})|_{T_{\alp\beta}}$ induces an isomorphism
$j_{\alp\beta}$ similarly to $j_{\alp}$. One can verify that both
$j_{\alp}|_{T_{\alp\beta}}$ and $j_{\beta}|_{T_{\alp\beta}}$ coincide
with $j_{\alp\beta}$.
\end{pf}\par
By this claim we can glue $\{ j_{\alp} \}_{\alp}$ to obtain an isomorphism
\begin{align*}
j: &\left[ Hom_T(q_*{\cal C}^{\bullet}( Hom_{X_T}( F^{\bullet}, {\cal G}_0
(K_X)), {\bold U}),\, K^{\bullet}) \right]_{-1}= \\
   & \left[ \bR Hom_T( \bR q_*( Hom_{X_T}(F^{\bullet}, G^{\bullet}(K_X))),
{\cal O}_T) \right]_{-1} \rightarrow \\
  & \qquad Hom_T\left( \left[ q_* {\cal C}^{\bullet}( Hom_{X_T}( F^{\bullet},
{\cal G}_0(K_X)), {\bold U}) \right]_1, {\cal O}_T \right) =
Ext^1_{X_T/T}( {\cal F}_0, {\cal G}_0(K_X))^{\vee}.
\end{align*}
Now this $j$, \eqref{SerreDuality} and \eqref{RqRhom} complete the
proof of this lemma.
\end{pf}\par
Remark that $Ext^1_{X_T/T}({\cal F}_0, {\cal G}_0(K_X))$ is not
isomorphic to $Ext^1_{X_T/T}({\cal G}_0, {\cal F}_0)^{\vee}$ in general.
\begin{lem}\label{lem:compatibleBchange}
A natural homomorphism
\[f^* Ext^1_{X_T/T}({\cal F}_0, {\cal G}_0(K_X)) \rightarrow
 Ext_{X_S/S}^1 (f^*{\cal F}_0, f^*{\cal G}_0(K_X))\]
is isomorphic for any $T$-scheme $f: S\rightarrow T$.
\end{lem}
\begin{pf}
This lemma is immediate from base change theorem
\cite[Page 104]{La:extfamily}.
\end{pf}\par
Now let us study a $T$-scheme $P^{\ff}$.
\begin{lem}\label{lem:PextPf}
There is a $T$-morphism $i_-: \PP( Ext_{X_T/T}^1({\cal F}_0, {\cal G}_0
(K_X))) \rightarrow P^{\ff}$.
\end{lem}
\begin{pf}
We shorten $\PP( Ext^1_{X_T/T}({\cal F}_0, {\cal G}_0(K_X)))$ to $\PP_-$,
and denote by $p_-: \PP_-\rightarrow T$ its structural morphism.
Proposition \ref{prop:SerreDuality} and Lemma \ref{lem:compatibleBchange}
lead to a natural isomorphism
\begin{multline*}
\Hom_{\PP_-}(p_-^* Ext_{X_T/T}^1( {\cal F}_0, {\cal G}_0(K_X)), {\cal O}(1))
\simeq \\
\Gamma( \PP_-, Ext_{X_{\PP_-}/{\PP_-}}^1({\cal G}_0, {\cal F}_0\otimes
{\cal O}_-(1))) \simeq \Ext^1_{X_{\PP_-}}({\cal G}_0, {\cal F}_0\otimes
{\cal O}_-(1)) 
\end{multline*}
since $Hom_{X_{\PP_-}/{\PP_-}}({\cal G}_0, {\cal F}_0\otimes {\cal O}_-(1))=0$
by base change theorem.
A tautological quotient line bundle
\begin{equation}\label{Om}
p_-^* Ext_{X_T/T}^1 ({\cal F}_0, {\cal G}_0(K_X)) \twoheadrightarrow
{\cal O}_-(1)
\end{equation}
on $\PP_-$ gives $\sigma\in\Ext^1_{X_{\PP_-}}({\cal G}_0,
{\cal F}_0\otimes {\cal O}_-(1))$ or an extension
\begin{equation}\label{ExtonPextm}
0 \longrightarrow {\cal F}_0\otimes {\cal O}_-(1) \longrightarrow {\cal V}_-
\longrightarrow {\cal G}_0 \longrightarrow 0.
\end{equation}
This ${\cal O}_{X_{\PP_-}}$-module ${\cal V}_-$ is $\PP_-$-flat.
For any point $t$ of $\PP_-$, Proposition \ref{prop:SerreDuality} and
Lemma \ref{lem:compatibleBchange} provide us with homomorphisms
\begin{multline*}
\kappa_1\circ (k(t)\otimes \Theta_q): \,
k(t)\otimes \Gamma\left(\PP_-, Ext_{X_{\PP_-}/{\PP_-}}^1({\cal G}_0,
{\cal F}_0\otimes {\cal O}-(1))\right) \rightarrow \\
k(t)\otimes \Gamma\left(\PP_-, Ext^1_{X_{\PP_-}/{\PP_-}}({\cal F}_0\otimes
{\cal O}_-(1), {\cal G}_0(K_X))^{\vee}\right) \rightarrow
\Ext^1_{X_{k(t)}}( {\cal F}_{0\, k(t)},\, {\cal G}_{0\, k(t)}(K_X))^{\vee}
\end{multline*}
and
\begin{multline*}
\Theta_{q\otimes k(t)}\circ \kappa_2: \,
k(t)\otimes \Gamma\left(\PP_-, Ext_{X_{\PP_-}/\PP_-}^1({\cal G}_0, 
{\cal F}_0 \otimes {\cal O}_-(1)) \right) \rightarrow \\
\Ext_{X_{k(t)}}^1( {\cal G}_{0\, k(t)}, {\cal F}_{0\, k(t)}) \rightarrow
\Ext_{X_{k(t)}}^1( {\cal F}_{0\, k(t)}, {\cal G}_{0\, k(t)}(K_X))^{\vee},
\end{multline*}
where $\kappa_i$ are natural maps.
In fact these homomorphisms are equal to each other because a trace map
$\operatorname{Tr}_q: R^2 q_*(K_X)\rightarrow {\cal O}_T$ is compatible
with base change by \cite[Page 172, Theorem 3.6.5]{Co:duality}.
The extension class of the exact sequence
\begin{equation}\label{ExtonPextmFibr}
0 \longrightarrow {\cal F}_{0\, k(t)} \longrightarrow {\cal V}_{-\, k(t)}
\longrightarrow {\cal G}_{0\, k(t)} \longrightarrow 0
\end{equation}
induced from \eqref{ExtonPextm} is equal to $\kappa_2(\sigma)\in
\Ext^1_{X_{k(t)}}( {\cal G}_{0\, k(t)}, {\cal F}_{0\, k(t)})$.
On the other hand $\kappa_1\circ (k(t)\otimes \Theta_q)(\sigma)\in
\Hom_{k(t)}(\Ext_{X_{k(t)}}^1({\cal F}_{0\, k(t)}, {\cal G}_{0\, k(t)}),
k(t))$ is nonzero since \eqref{Om} is surjective.
Therefore we see that \eqref{ExtonPextmFibr} is not trivial, which means
that ${\cal V}_-$ is a
flat family of $a_-$-stable sheaves by Lemma \ref{lem:HNF}.
${\cal V}_-$ gives a morphism $i_-: \PP_-\rightarrow \modm$.
It's easy to see that $i_-$ factors through $\PP_-\rightarrow P^{\ff}
\hookrightarrow \modm$ and that $i_-: \PP_-\rightarrow P^{\ff}$ is
a $T$-morphism.
\end{pf}\par
By \eqref{HNFmf} and \eqref{FandF0}, we have a natural exact sequence
\begin{equation}\label{rewriteHNFmf}
0 \longrightarrow (\tau^Q)^*{\cal F}_0 \otimes L_1 \longrightarrow
{\cal U}_-|_{Q^{\ff}} \longrightarrow (\tau^Q)^* {\cal G}_0\otimes L_2
\rightarrow 0
\end{equation}
on $X_{Q^{\ff}}$.
Similarly to the proof of the lemma above, one can show that
\[\Ext^1_{X_{\Qf}}( (\tau^Q)^* {\cal G}_0 \otimes L_2, (\tau^Q)^*
{\cal F}_0 \otimes L_1) \simeq
\Hom_{\Qf}((\tau^Q)^* Ext^1_{X_T/T}({\cal F}_0, {\cal G}_0(K_X)),\, L_1\otimes
L_2^{\vee}), \]
and that the homomorphism
$(\tau^Q)^* Ext_{X_T/T}^1( {\cal F}_0, {\cal G}_0(K_X)) \twoheadrightarrow
L_1\otimes L_2^{\vee}$
induced by \eqref{rewriteHNFmf} is surjective.
Thus this gives a morphism
$j^Q: \Qf\rightarrow \PP( Ext^1_{X_T/T}({\cal F}_0, {\cal G}_0(K_X)))$.
One can check that $j^Q$ is $\bar{G}$-invariant.
As a result, $j^Q$ descends to a morphism
\begin{equation}\label{jm}
j_-: P^{\ff}\rightarrow \PP(Ext_{X_T/T}({\cal F}_0, {\cal G}_0(K_X))).
\end{equation}

\begin{lem}
For morphisms $i_-$ in Lemma \ref{lem:PextPf} and $j_-$ at \eqref{jm},
it holds that $i_-\circ j_-= \id_{P^{\ff}}$ and that $j_-\circ i_-=
\id_{\PP_-}$.
\end{lem}
\begin{pf}
Since $\pi_-: \Qf\rightarrow P^{\ff}$ is a categorical quotient by $\bar{G}$,
$i_-\circ j_-$ is equal to $\id_{P^{\ff}}$ if and only if $(i_-\circ
j_-)\circ \pi_-= i_-\circ j^Q$ is equal to $\pi_-$. One can readily verify
this, and hence its proof is omitted.
$T$-morphism $i_-: \PP( Ext^1_{X_T/T}({\cal F}_0, {\cal G}_0(K_X)))=
\PP_-\rightarrow P^{\ff}\hookrightarrow \modm$ is induced from an
${\cal O}_{X_T}$-module ${\cal V}_-$ in \eqref{ExtonPextm}, and hence
one can find an affine open covering $\{ P_{\alp}\}_{\alp}$ of $\PP_-$
and a morphism $i_{\alp}: P_{\alp}\rightarrow\Qf$ such that
$\pi_-\circ i_{\alp}= i_-|_{U_{\alp}}$.
It's easy to show that $j^Q\circ i_{\alp}= j_-\circ i_-|_{P_{\alp}}:
P_{\alp}\rightarrow \PP_-$ is equal to $\id_{P_{\alp}}$, and hence
its proof is left to the reader.
\end{pf}\par
Summing up, we get the following:
\begin{prop}\label{prop:PfisPext}
Fix an element $\ff$ of $A^+(a)$. We define a scheme $T$,
${\cal O}_{X_T}$-modules ${\cal F}_0$ and ${\cal G}_0$, and line bundles
$L_1$ and $L_2$ over $Q^{\ff}$ as in \eqref{F0G0} and in \eqref{FandF0}.
\begin{enumerate}
\item $P^{\ff}$ can be regarded as a $T$-scheme.
\item There is an isomorphism $j_-: P^{\ff}\rightarrow \PP(
      Ext^1_{X_T/T}
( {\cal F}_0, {\cal G}_0(K_X)))$ over $T$ such that $L_1\otimes L_2^{\vee}
\in\Pic( \Qf)$ in \eqref{FandF0} is equal to $(j_-\circ\pi_-)^* {\cal O}_-(1)$,
where ${\cal O}_-(1)$ is the tautological line bundle of
$\PP ( Ext_{X_T/T}^1({\cal F}_0, {\cal G}_0(K_X)))$.
\end{enumerate}
\end{prop}
%
%
%
\section{Algebro-geometric analogy of $\mu$-map and \\
the Donaldson polynomial}\label{ss:AGanalogy}
From now on we shall consider the case of $c_1=0$. Hence $\modm$ stands for
$M_{a_-}(0,c_2)$, and so on.
We begin with reviewing the algebro-geometric analogy $\mu_-: \NS(X)
\rightarrow \NS(\modm)$ of $\mu$-map, which was introduced in
\cite{Li:AGinterpret}.
Let $C\subset X$ be a nonsingular curve, and $\theta_C$ a line bundle
on $C$ with $\deg(\theta_C)=g(C)-1$.
For a universal sheaf ${\cal U}_-$ of $\Qm$ on $X_{\Qm}$, one can show
that a complex $\bR \pr_{2\, *}({\cal U}_-|_C\otimes \theta_C)$ on
$\Qm$ locally is quasi-isomorphic to a bounded complex of free modules
of finite rank.
Thus its determinantal line bundle $\det \bR \pr_{2\, *}
({\cal U}_-|_C\otimes \theta_C)$ on $\Qm$ exists.
In fact this line bundle descends to a line bundle ${\cal L}_{\modm}
(C,\theta_C)^{\vee}$ on $\modm$, and its algebraic equivalence class
$[ {\cal L}_{\modm}(C,\theta_c)^{\vee}]\in\NS(\modm)$
is independent of the choice of $\theta_C$.
It can be checked that the correspondence $C\mapsto [{\cal L}_{\modm}
(C,\theta_C)]$ induces a homomorphism
$\mu_-: \NS(X)\rightarrow \NS(\modm)$.
One can also construct a homomorphism $\mu_+: \NS(X)\rightarrow
\NS(\modp)$ likewise. Let
$ \modm \overset{\phi_-}{\longleftarrow} \tilmodm 
\overset{\phi_+}{\longrightarrow}\modp $
be the sequence of morphisms \eqref{flipm}.
For $\ff\in A^+(a)$, we denote by $E^{\ff}$ the Cartier divisor
$\phi_-^{-1}( P^{\ff})$ on $\tilmodm$, which is equal to
$\tilde{\phi}^{-1}_+( P^{\ff})$ by Corollary \ref{cor:inverseofPpisEm}.
\begin{lem}\label{lem:differenceofmu}
For $\alp\in\NS(X)$, it holds that 
\[\phi_-^* \mu_-(\alp)-\bar{\phi}_+^* \mu_+(\alp)= 
\sideset{}{_{\ff\in A^+(a)}}\sum {\cal O}_{\tilmodm}
\left( \langle f\cdot\alp /2 \rangle\, E^{\ff}\right) \]
in $\NS(\tilmodm)$.
\end{lem}
\begin{pf}
For the simplicity of notation, we prove this lemma in case of
$\sharp A^+(a)=1$. It's easy to extend the proof to general case.
Let $C$ and $\theta_C$ be as explained above, and
\begin{equation}\label{elemtransmAgain}
0 \longrightarrow {\cal W}_+ \longrightarrow \tilde{\cal U}_- 
\longrightarrow \tilde{\cal G} \longrightarrow 0
\end{equation}
the exact sequence \eqref{elemtransm} on $X_{\tilQm}$.
Since ${\cal U}_-$ is a flat family of torsion free sheaves, one can
show that a ${\cal O}_{C_{\tilQm}}$-module $\tilde{\cal U}_-|_C$ is
$\tilQm$-flat.
By using the method of C\v{e}ch complex, one get a quasi-isomorphism
\begin{equation}\label{LRtoRL}
\bL f^* \bR \pr_{2\, *} (\tilde{\cal U}_-|_C \otimes \theta_C) \rightarrow
\bR \pr'_{2\, *} \bL f^{'*} (\tilde{\cal U}_-|_C \otimes \theta_C)=
\bR \pr'_{2\, *} f^{'*} (\tilde{\cal U}_-|_C \otimes \theta_C),
\end{equation}
where
\begin{equation*}
\xymatrix{
 C_S \ar[r]_{f'} \ar[d]_{\pr'_2} & C_{\tilQm} \ar[d]^{\pr_2} \\
 S \ar[r]^f & \tilQm}
\end{equation*}
is a fiber product. The analogy to these result about $(\tilde{\cal U}_-,
\tilQm)$ also holds to $({\cal W}_+, \tilQm)$ and $(\tilde{\cal G}, D_-)$.
\eqref{elemtransmAgain} gives a triangle
\[ \bR \pr_{2\, *}({\cal W}_+|_C \otimes \theta_C) \longrightarrow
   \bR \pr_{2\, *}(\tilde{\cal U}_-|_C \otimes \theta_C) \longrightarrow
   \bR \pr_{2\, *}(\tilde{\cal G}_-|_C \otimes \theta_C) \]
in $D(\tilQm)$, and hence an isomorphism
\begin{equation}\label{det&elemtrans}
\det \bR\pr_{2\, *}(\tilde{\cal U}_-|_C \otimes \theta_C) \simeq
\det \bR\pr_{2\, *}({\cal W}_+|_C \otimes \theta_C)\cdot
\det \bR\pr_{2\, *}(\tilde{\cal G}|_C \otimes \theta_C)
\end{equation}
in $\Pic(\tilQm)$ is induced.

$\det \bR \pr_{2\, *}( {\cal W}_+|_C \otimes \theta_C)$ naturally is
isomorphic to $\tilde{\pi}_+^* \bar{\phi}_+^* {\cal L}_{\modp}(C,
\theta_C)^{\vee}$.
Indeed, Lemma \ref{lem:locallift} gives an open covering
$\bigcup_{\alp}\, U_{\alp}$ of $\tilQm$, a morphism
$\bar{\varphi}_+^{\alp}: U_{\alp} \rightarrow \Qp$, and an isomorphism
of $X_{U_{\alp}}$-modules
$\Phi_+^{\alp}: {\cal W}_+|_{U_{\alp}} \rightarrow
(\bar{\varphi}^{\alp}_+)^* {\cal U}_+$.
By \eqref{LRtoRL} $\Phi_+^{\alp}$ naturally induces an isomorphism
\begin{multline*}
\det( \Phi_+^{\alp}): \det \bR \pr_{2\, *}({\cal W}_+|_C \otimes
\theta_C)|_{U_{\alp}} \rightarrow
\det \bR \pr_{2\, *} (\bar{\varphi}_+^{\alp})^* ({\cal U}_+|_C \otimes
\theta_C) \rightarrow \\
(\bar{\varphi}_+^{\alp})^* \det \bR \pr_{2\, *}( {\cal U}_+|_C \otimes
\theta_C)=
(\bar{\varphi}^{\alp}_+)^* \pi_+^* {\cal L}_{\modp}(C, \theta_C)^{\vee}=
\tilde{\varphi}_+^* {\cal L}_{\modp}(C, \theta_C)^{\vee}|_{U_{\alp}}.
\end{multline*}
In addition, if $\alp\neq\beta$ then the isomorphism
\begin{equation}\label{PhibPhia}
(\Phi_+^{\beta})^{-1}\circ \Phi_+^{\alp}: {\cal W}_+|_{U_{\alp\beta}}
\rightarrow {\cal W}_+|_{U_{\alp\beta}} 
\end{equation}
on $X_{U_{\alp\beta}}= X_{U_{\alp}\cap U_{\beta}}$ is given by
$\lambda_{\alp\beta}\in \Gamma({\cal O}^{\times}_{U_{\alp\beta}})$
since ${\cal W}_+|_{U_{\alp\beta}}$ is a flat family of simple sheaves
as mentioned right after Corollary \ref{cor:WpisStables}.
One can define the rank $R$ of a perfect complex $\bR \pr_{2\, *}
({\cal W}_+|_C\otimes \theta_C)$, and then
\[ \det(\Phi_+^{\beta})^{-1} \circ \det(\Phi_+^{\alp}) :
   \det \bR \pr_{2\, *}( {\cal W}_+|_C\otimes \theta_C)|_{U_{\alp\beta}}
\rightarrow  \det \bR \pr_{2\, *}( {\cal W}_+|_C \otimes \theta_C)
|_{U_{\alp\beta}} \]
is given by $\lambda_{\alp\beta}^{\times R}$.
This $R$ turns out to be zero because the Riemann-Roch theorem implies that
$\chi( C_{k(t)}, {\cal W}_+|_C \otimes {\theta}_C \otimes k(t))=0$
for every $t\in U_{\alp\beta}$.
Hence we can glue $\det (\Phi_+^{\alp})$ to obtain an isomorphism
\begin{equation}\label{detWp}
 \det \bR \pr_{2\, *} ({\cal W}_+|_C \otimes {\theta}_C) \simeq
 \tilde{\varphi}_+^* {\cal L}_{\modp} (C,\theta_C)^{\vee}=
 \tilde{\pi}_-^* \bar{\phi}_+^* {\cal L}_{\modp} (C, \theta_C)^{\vee}.
\end{equation}
From \eqref{LRtoRL}, we can get a natural isomorphisms
\[ \det \bR \pr_{2\, *}(\tilde{\cal U}_-|_C \otimes \theta_C) \rightarrow
   \varphi_-^* \det \bR \pr_{2\, *}( {\cal U}_-|_C \otimes \theta_C)=
   \tilde{\pi}_-^* \phi_-^* {\cal L}_{\modm}(C, \theta_C)^{\vee}. \]
Hence by \eqref{det&elemtrans} and \eqref{detWp}
$\tilde{\pi}_-^*( \tilde{\phi}_+^* {\cal L}_{\modp}(C, \theta_C)-
 \phi_-^* {\cal L}_{\modm}(C,\theta_C)) \simeq
 \det \bR \pr_{2\, *}(\tilde{\cal G}|_C \otimes \theta_C)$.
$\tilde{\cal G}|_C\otimes \theta_C$ is a sheaf on $C_{D_-}\subset
C_{\tilQm}$, and so $\det \bR \pr_{2\, *}(\tilde{\cal G}|_C\otimes
\theta_C)$ can be regarded as
$\det \bR\pr_{2\, *} j^C_* (\tilde{\cal G}|_C \otimes \theta_C)=
\det j_* \bR \pr'_{2\, *}(\tilde{\cal G}|_C \otimes \theta_C)$,
where
\begin{equation*}
\xymatrix{
 C_{D_-} \ar@{^{(}-}[r]^{j^C} \ar[d]^{\pr'_2} & C_{\tilQm} \ar[d]^{\pr_2} \\
 D_- \ar@{^{(}-}[r]^j & \tilQm }
\end{equation*}
is a fiber product.
By the Riemann-Roch theorem $\chi(C_{k(t)}, \tilde{\cal G}|_C\otimes 
\theta_C \otimes k(t))= -f\cdot C/2 $ for every $t\in D_-$.
Thus the rank of a complex
$\bR \pr'_{2\, *}(\tilde{\cal G}|_C \otimes \theta_C)$ on $D_-$ is equal
to $-f\cdot C/2$.
In view of this we can prove that
$\det \bR \pr_{2, *}(\tilde{\cal G}|_C \otimes {\theta}_C)=
 -\langle f\cdot C/2\rangle\, D_-$ in $\Pic( \tilQm)$; its proof is omitted.
Summing up, we obtain an isomorphism
\begin{multline}\label{isomofGlb}
\tilde{\pi}_-^* ( \tilde{\phi}_+^* {\cal L}_{\modp}(C,\theta_C)-
 \phi_-^* {\cal L}_{\modm}(C, \theta_C)) \simeq \\
-\sideset{}{_{\ff\in A^+(a)}}\sum {\cal O}_{\tilQm} (\langle f\cdot C/2\rangle
\, D^{\ff})=
-\sideset{}{_{\ff\in A^+(a)}}\sum \tilde{\pi}_-^* {\cal O}_{\tilmodm}
( \langle f\cdot C/2\rangle\, E^{\ff})
\end{multline}
in $\Pic(\tilQm)$.
Moreover, both sides in \eqref{isomofGlb} respectively have a natural
$\bar{G}$-linearized structure. One can check that \eqref{isomofGlb} is an
isomorphism of $\bar{G}$-linearized line bundles.
By \cite[Theorem 4]{Se:geomred} and \cite[Page 87, Theorem 4.2.16]{HL:text}
the natural homomorphism
\begin{equation}\label{Picinj}
 \tilde{\pi}_-^*: \Pic( \tilmodm) \rightarrow \Pic_{\bar{G}}(\tilQm)
\end{equation}
is injective, where $\Pic_{\bar{G}}(\tilQm)$ is the group of
$\bar{G}$-linearized line bundles on $\tilQm$.
%
Thereby \eqref{isomofGlb} and \eqref{Picinj} complete the proof of this lemma.
\end{pf}\par
Now we assume that
\begin{equation}\label{expecteddim}
\dim M_{H_{\pm}}(0,c_2)=4c_2-3\chi(\Ox)=d(c_2)
\end{equation}
and that
\begin{equation}\label{codimgeq2}
\operatorname{codim}(M_{\pm}, P_{\pm})\geq 2.
\end{equation}
These assumptions can be considered to be reasonably weak by the
following lemma.
\begin{lem}\label{lem:condtionisWeak}
Let $\Amp(X)$ be the ample cone of $X$, and $S\subset \Amp(X)$ a compact
subset containing $H_{\pm}$.
If $c_2$ is sufficiently large with respect to $S$, then assumptions
\eqref{expecteddim} and \eqref{codimgeq2} hold good.
\end{lem}
\begin{pf}
Refer to \cite[Theorem 2]{Zu:genericsmooth}, \cite{GL:irreducible},
and the proof of \cite[Theorem 2.3.]{Qi:birational}.
\end{pf}\par
By \eqref{expecteddim} we can define a multilinear map
\begin{equation*}
\gamma_{\pm}=\gamma_{\pm}(c_2): \Sym^{d(c_2)}\, \NS(X) \rightarrow \ZZ
\end{equation*}
by $\gamma_{\pm}(\alp_1, \ldots , \alp_{d(c_2)})= \mu_{\pm}(\alp_1)
\pmb{\cdot}\cdots \pmb{\cdot} \mu_{\pm}(\alp_{d(c_2)})$ using 
the intersection number of
line bundles on $M_{\pm}=M_{a_{\pm}}(0,c_2)$.
Similarly, a multilinear map
\begin{equation*}
\gamma_{H_{\pm}}=\gamma_{H_{\pm}}(c_2): \Sym^{d(c_2)}\, \NS(X)\rightarrow
\ZZ
\end{equation*}
can be defined by the intersection number of line bundles on $M_{H_{\pm}}
(0,c_2)$.

\begin{prop}\cite[Page 456]{Li:AGinterpret}\label{prop:AGinterpretIntro}
Suppose that $X$ is simply connected and that $p_g(X)>0$.
Then there is such a constant $A({\cal S})$ depending on a compact subset
${\cal S}\subset\Amp(X)$ as satisfies the following: \\
If $c_2\geq A({\cal S})$ and if some rational multiple of an ample line
bundle $H$ is contained in ${\cal S}$, then $\gamma_H(c_2)$ is equal to
the restriction of the Donaldson invariant $q(c_2): \Sym^{d(c_2)}\,
H_2(X,\ZZ)\rightarrow \ZZ$ to $\Sym^{d(c_2)}\, \NS(X) $.
In particular $\gamma_H(c_2)$ is independent of an ample line bundle $H$
contained in $\QQ \cdot {\cal S}$.
\end{prop}
Hence from well-known argument about Donaldson polynomials, which
originated in differential geometry, $\gamma_H(c_2)$ is independent of
an ample line bundle $H$ when $p_g(X)>0$ and $X$, $c_2$, ${\cal S}$ and
$H$ are as in this proposition. 
We would like to observe this independence from an algebraic
geometric point of view.
For this reason we shall study $\mu_-(C)^{d(c_2)}-\mu_+(C)^{d(c_2)}$
for a nonsingular curve $C$ in $X$. We often shorten $d(c_2)$ to $d$.
%
%

Since \eqref{codimgeq2} both $\phi_-: \tilmodm \rightarrow \modm$ and
$\bar{\phi}_+: \tilmodm \rightarrow \modp$ are birational,
$\mu_+(C)^d- \mu_-(C)^d$ is equal to
$\bar{\phi}_+^* \mu_+(C)^d- \phi_-^* \mu_-(C)^d$.
For $f\in\NS(X)$ , we set $\lambda_f^C:=\langle f\cdot C/2\rangle$.
Then Lemma \ref{lem:differenceofmu} implies that
{\allowdisplaybreaks
\begin{gather}\label{rewritedifference1}
\begin{split}
 & \bar{\phi}_+^*\mu_+(C)^d - \phi_-^*\mu_-(C)^d  \\
= & \,\Bigl\{ \bar{\phi}_+^*\mu_+(C) - \phi_-^*\mu_-(C) \Bigr\} \cdot
    \sum_{k=0}^{d-1} \phi_-^* \mu_-(C)^k \pmb{\cdot} \bar{\phi}_+^*
    \mu_+(C)^{d-1-k} \\
= & \,\sum_{k=0}^{d-1} \left[ \phi_-^* \mu_-(C)^k \pmb{\cdot} \bar{\phi}_+^*
\mu_+(C)^{d-1-k} \pmb{\cdot} -\sum_{\ff} \lambda_f^C E^{\ff} 
\right]_{\tilmodm},
\end{split}
\end{gather}
where $[\quad]_{\tilmodm}$ stands for the multiplication of line
bundles is calculated on $\tilmodm$.
By \cite[Page 297, Proposition 4]{Kl:ampleness}, \eqref{rewritedifference1}
is equal to
\begin{equation}\label{rewritedifference2}
\begin{split}
\sum_{\ff\in A^+(a)} &-\lambda_f^C\,
\sum_{k=0}^{d-1} \left[ \phi_-^* \mu_-(C)\big| _{E^{\ff}}^k \pmb{\cdot} 
\Bigl\{ \phi_-^*
\mu_-(C) -\sum_{\ff}  \lambda_f^C E^{\ff} \Bigr\}\Big|_{E^{\ff}}^{d-1-k}
\right]_{E^{\ff}}  \\
=\sum_{\ff} & -\lambda_f^C \sum_{k=0}^{d-1} \left[
\phi_-^* \mu_-(C)\big|_{E^{\ff}}^k \pmb{\cdot } \Bigl\{ \phi_-^*\mu_-(C)-
\lambda_f^C E^{\ff} \Bigr\} \Big|_{E^{\ff}}^{d-1-k} \right]_{E^{\ff}}, 
\end{split}
\end{equation}
since $E^{\ff}\cap E^{\ff'}=\emptyset$ if $\ff$ and $\ff'$ are different
by \ref{lem:setforcenter}. 
$\phi_-$ and $\bar{\phi}_+$ induce a commutative diagram
\begin{equation*}
\xymatrix{
E^{\ff} \ar[d]_{\phi_-} \ar[r]_{\bar{\phi}_+} & P^{-\ff} \ar[d]^{\tau_+}\\
P^{\ff} \ar[r]^{\tau_-} & T. }
\end{equation*}
Let ${\cal F}_0$ and ${\cal G}_0$ be ${\cal O}_{X_T}$-modules defined in
\eqref{F0G0}.
By Proposition \ref{prop:PfisPext} there are tautological line bundles
${\cal O}^{\ff}_-(1)$ and ${\cal O}^{-\ff}_+(1)$ on, respectively,
$P^{\ff}=\PP( Ext_{X_T/T}^1({\cal F}_0, {\cal G}_0(K_X)))$ and
$P^{-\ff}= \PP( Ext_{X_T/T}^1({\cal G}_0, {\cal F}_0(K_X)))$.
\begin{lem}\label{lem:OeE}
A line bundle ${\cal O}_{E^{\ff}}(-E^{\ff})$ on $E^{\ff}$ is naturally
isomorphic to $\phi_-^*{\cal O}_-^{\ff}(1)+ \bar{\phi}_+^* {\cal O}_+^{-\ff}
(1)$.
\end{lem}
\begin{pf}
We shorten ${\cal O}_-^{\ff}(1)$ and ${\cal O}_+^{-\ff}(1)$ to, 
respectively, ${\cal O}_-(1)$ and ${\cal O}_+(1)$.
Let $D^{\ff}$ denote a closed subscheme $(\tilde{\pi}_-)^{-1}(E^{\ff})$
of $\tilQm$. Then
\begin{equation*}
\xymatrix{
Q^{\ff} \ar[d]_{\pi_-} & D^{\ff} \ar[d]_{\tilde{\pi}_-} \ar[l]^{\varphi_-} &
D^{\ff}\cap U_{\alp} \ar@{_{(}-}[l] \ar[r]_-{\bar{\varphi}^{\alp}_+} &
Q^{-\ff} \ar[d]^{\pi_+} \\
P^{\ff} \ar[r]_{\phi_-} & E^{\ff} \ar[rr]^{\bar{\phi}_+} && P^{-\ff} }
\end{equation*}
is commutative for $U_{\alp}\subset\tilQm$ and $\bar{\varphi}_+^{\alp}$
in Lemma \ref{lem:locallift}. By \eqref{FandF0}, we can rewrite the
exact sequence \eqref{inducedHNFforWp} to obtain
\begin{equation}\label{rewriteInducedHNFforWp}
0 \longrightarrow {\cal G}_0\otimes L_2(-D^{\ff}) \longrightarrow
{\cal W}_+|_{X_{D^{\ff}}} \longrightarrow {\cal F}_0 \otimes L_1
\longrightarrow 0
\end{equation}
on $X_{D^{\ff}}$. Next, let
\begin{equation}\label{HNFpAgain}
0 \longrightarrow {\cal G}' \longrightarrow {\cal U}_+|_{Q^{-\ff}}
\longrightarrow {\cal F}' \longrightarrow 0
\end{equation}
be the exact sequence \eqref{HNFp} on $X_{Q^{-\ff}}$. Similarly to
\eqref{FandF0}, there are isomorphisms
\begin{equation*}
{\cal F}'\simeq {\cal F}_0\otimes M_1 \text{ and }
{\cal G}'\simeq {\cal G}_0\otimes M_2
\end{equation*}
with some line bundles $M_1$ and $M_2$ on $Q^{-\ff}$. Analogously to
Proposition \ref{prop:PfisPext}, $M_2\otimes M_1^{\vee}$ is isomorphic
to $\pi^*_+ {\cal O}_+(1)$.
Thus we obtain an exact sequence
\begin{equation}\label{rewriteHNFpAgain}
0 \longrightarrow {\cal G}_0\otimes (\bar{\varphi}_+^{\alp})^* M_2
\longrightarrow (\bar{\varphi}_+^{\alp})^* {\cal U}_+|_{D^{\ff}\cap
U_{\alp}} \longrightarrow {\cal F}_0 \otimes (\bar{\varphi}_+^{\alp})^*
M_1 \longrightarrow 0
\end{equation}
on $X_{D^{\ff}}\cap U_{\alp}$, pulling \eqref{HNFpAgain} back by
$\id_X\times \bar{\varphi}_+^{\alp}: X_{D^{\ff}\cap U_{\alp}} \rightarrow
X_{Q^{-\ff}}$.
Connecting \eqref{rewriteInducedHNFforWp} and \eqref{rewriteHNFpAgain}
by the isomorphism $\Phi_+^{\alp}$ in Lemma \ref{lem:locallift}, we
get the following:
\begin{equation}\label{connectRewrites}
\xymatrix{
0 \ar[r] & {\cal G}_0\otimes L_2(-D^{\ff}) \ar[r] \ar@{.>}[d]^{r'_{\alp}} &
{\cal W}_+|_{D^{\ff}\cap U_{\alp}} \ar[r] \ar[d]^{\Phi_+^{\alp}} &
{\cal F}_0 \otimes L_1|_{U_{\alp}} \ar[r] \ar@{.>}[d]^{r_{\alp}} & 0 \\
0 \ar[r] & {\cal G}_0 \otimes (\bar{\varphi}_+^{\alp})^* M_2 \ar[r] &
(\bar{\varphi}_+^{\alp})^* {\cal U}_+|_{D^{\ff}\cap U_{\alp}} \ar[r] &
{\cal F}_0 \otimes (\bar{\varphi}_+^{\alp})^* M_1 \ar[r] & 0 }
\end{equation}
As observed in the proof of Lemma \ref{lem:inverseofPpisDm}, there is an
isomorphism
$r_{\alp}: {\cal F}_0\otimes L_1|_{U_{\alp}} \rightarrow {\cal F}_0 \otimes
(\bar{\varphi}_+^{\alp})^* M_1 $ which makes \eqref{connectRewrites}
commutative. This $r_{\alp}$ induces an isomorphism $r'_{\alp}$ in
\eqref{connectRewrites}.
Because both ${\cal F}_0$ and ${\cal G}_0$ are flat families of simple
sheaves, $r_{\alp}$ and $r'_{\alp}$ induce isomorphisms
\begin{equation*}
\Gamma( r_{\alp}): L_1|_{U_{\alp}} \rightarrow (\bar{\varphi}_+^{\alp})^*
M_1 \quad\text{ and }\quad
\Gamma( r'_{\alp}): L_2\otimes {\cal O}_{D^{\ff}}(-D^{\ff})|_{U_{\alp}}
\rightarrow (\bar{\varphi}_+^{\alp})^* M_2
\end{equation*}
of line bundles on $D^{\ff}\cap U_{\alp}$.
$\Gamma( r_{\alp})$ and $\Gamma( r'_{\alp})$ induce an isomorphism
\begin{multline*}
\Gamma(r_{\alp})^{-1}\cdot \Gamma(r'_{\alp}):
{\cal O}_{D^{\ff}}(-D^{\ff})|_{U_{\alp}}= 
\tilde{\pi}_-^* {\cal O}_{E^{\ff}}(-E^{\ff})|_{U_{\alp}} \rightarrow \\
(L_1\otimes L_2^{\vee})|_{U_{\alp}} \otimes (\bar{\varphi}_+^{\alp})^*
(M_2\otimes M_1^{\vee}) \simeq 
\tilde{\pi}_-^* \phi_-^* {\cal O}_-(1)|_{U_{\alp}} \otimes
(\bar{\varphi}_+^{\alp})^* \pi_+^* {\cal O}_+(1) \\
=\tilde{\pi}^*\bigl( \phi_-^*{\cal O}_-(1) + \bar{\phi}_+^* {\cal O}_+(1)
\bigr) |_{D^{\ff}\cap U_{\alp}} \qquad
\end{multline*}
By \eqref{PhibPhia} one can check that
$ \Gamma( r_{\alp})^{-1} \cdot \Gamma( r'_{\alp})=
  \Gamma( r_{\beta})^{-1} \cdot \Gamma(r'_{\alp})$,
and hence glue $\Gamma( r_{\alp})^{-1}\cdot \Gamma( r'_{\alp})$ to
obtain an isomorphism
\begin{equation*}
\tilde{\pi}_-^* {\cal O}_{E^{\ff}}(-E^{\ff}) \simeq
\tilde{\pi}_-^* ( \phi_-^* {\cal O}_-(1)+ \bar{\phi}_+^* {\cal O}_+(1))
\end{equation*}
of line bundles on $D^{\ff}$.
One can also check this is an isomorphism of $\bar{G}$-linearized line
bundles.
Then we complete the proof of this lemma in similar fashion to
the proof of Lemma \ref{lem:differenceofmu}.
\end{pf}\par
From \eqref{rewritedifference2} and Lemma \ref{lem:OeE} we obtain that
\begin{multline}\label{rewritedifference3}
\mu_-(C)^{d(c_2)} -\mu_+(C)^{d(c_2)}= 
\sum_{\ff\in A^+(a)} -\lambda_f^C \cdot\\
\sum_{k=0}^{d(c_2)-1} \left[ \phi_-^* \mu_-(C)\big|_{E^{\ff}}^k \pmb{\cdot}
 \Bigl\{ \phi_-^* \mu_-(C)\big|_{E^{\ff}} + \lambda_f^C ( {\cal O}_-^{\ff}(1)
+ {\cal O}_+^{-\ff}(1) ) \Bigr\}^{d(c_2)-1-k} \right]_{E^{\ff}}.
\end{multline}
In the following section, we shall in detail examine the right side of
this equation in some special case.
%
%
\section{The relation to the intersection theory of 
$\PP( {\cal A}_-)\times \PP({\cal A}_-^{\vee})$}\label{ss:PtimesP}
From now on, adding to \eqref{expecteddim} and \eqref{codimgeq2} we
assume that the irregularity $q(X)=0$ and that 
\begin{equation}\label{pg>0}
\text{some section $\kappa\in\Gamma(K_X)$ gives a
nonsingular curve ${\cal K}\subset X$}
\end{equation}
in view of Proposition \ref{prop:AGinterpretIntro}. (We can expect this
will be weakened to the condition $p_g(X)>0$; to do so, we
have to adjust the assumption in Proposition \ref{prop:PveetimesP}.)
Moreover we assume the following about $\ff=(f,m,n)\in A^+(a)$:
\begin{equation}\label{ExtConstant}
\dim\Ext^1_{X_t}({\cal F}_0\otimes k(t), {\cal G}_0\otimes k(t))
=L_+ \text{ and } \dim \Ext^1_{X_t}({\cal G}_0\otimes k(t), {\cal F}_0
\otimes k(t))=L_-
\end{equation}
are independent of $t\in T=\Pic^{[f/2]}(X)\times \Hilb^m(X)\times \Hilb^n(X)$,
where ${\cal F}_0$ and ${\cal G}_0$ are ${\cal O}_{X_T}$-modules defined in
\eqref{F0G0}.
This assumption \eqref{ExtConstant} holds good if, for example, $K_X$ is
numerically equivalent to zero, but is not weak at all in general.
Assuming \eqref{ExtConstant} we see that both
\begin{equation}\label{Apm}
{\cal A}_-= Ext^1_{X_T/T}({\cal F}_0, {\cal G}_0(K_X)) \,\text{ and }\,
{\cal A}_+= Ext^1_{X_T/T}({\cal G}_0, {\cal F}_0(K_X))
\end{equation}
are locally free ${\cal O}_T$-modules, and hence $P^{\pm\ff}=
\PP({\cal A}_{\mp})$ are projective bundles over a nonsingular scheme $T$.
Under these assumptions we would like to examine
\begin{equation}\label{differenceofmu:Ef}
\begin{split}
\sum_{k=0}^{d(c_2)-1} & \left[ \phi_-^* \mu_-(C)\big|_{E^{\ff}}^k \pmb{\cdot}
\bar{\phi}_+^* \mu_+(C)\big|_{E^{\ff}}^{d-1-k} \right]_{E^{\ff}}= \\
\sum_{k=0}^{d(c_2)-1} & \left[ \phi_-^* \mu_-(C)\big|_{E^{\ff}}^k 
\pmb{\cdot} \bigl\{ \phi_-^* \mu_-(C)+\lambda_f^C ( {\cal O}_-^{\ff}(1)+
{\cal O}_+^{-\ff}(1)) \bigr\} \big|_{E^{\ff}}^{d-1-k} \right], 
\end{split}
\end{equation}
which appeared in \eqref{rewritedifference3}.
We shorten ${\cal O}_-^{\ff}(1)$ and ${\cal O}_+^{-\ff}$ to, respectively,
${\cal O}_-(1)$ and ${\cal O}_+(1)$ for the time being.
Since $P^{\pm\ff}$ are projective bundles over $T$ there are line bundles
${\beta}_{\pm}$ on $T$ and integers $N_{\pm}$ such that
\begin{equation*}
\mu_-(C)|_{P^{\ff}} = \tau_-^* (\beta_-)+ {\cal O}_-(N_-) \,\text{ and }\,
\mu_+(C)|_{P^{-\ff}}= \tau_+^* (\beta_+)+ {\cal O}_+(N_+)
\end{equation*}
in $\Pic(P^{\pm\ff})$.
By Lemma \ref{lem:differenceofmu} and Lemma \ref{lem:OeE} we have
\begin{equation}\label{rewritedifferenceofmu}
  \phi_-^* \tau_-^*(\beta_- -\beta_+) +\phi_-^*{\cal O}_-(N_-)-
\bar{\phi}_+^* {\cal O}_+(N_+) 
=  -\lambda_f^C (\phi_-^* {\cal O}_-(1)+ \bar{\phi}_+^* {\cal O}_+(1))
\end{equation}
in $\Pic( E^{\ff})$.
Suppose that $N_+\neq \lambda_f^C = \langle f\cdot C/2\rangle$.
Then \eqref{rewritedifferenceofmu} implies that ${\cal O}_{E^{\ff}}$ is
$\phi_-$-ample since
${\cal O}_{E^{\ff}}(-E^{\ff})=\phi_-^*{\cal O}_-(1)+ \bar{\phi}_+^*
{\cal O}_+(1)$ is $\phi_-$-ample. By \cite[II.5.1]{EGA:II}, the proper morphism
$\phi_-: E^{\ff}\rightarrow P^{\ff}$ should be finite if ${\cal O}_{E^{\ff}}$
is $\phi_-$-ample. This contradicts \eqref{codimgeq2}.
Therefore as divisors on $E^{\ff}$ we have 
\begin{equation*}
\phi_-^* \mu_-(C)|_{E^{\ff}}= \phi_-^*\{ \tau_-^*\beta +
{\cal O}_-(-\lambda_f^C) \} \quad\text{ and }\quad
\bar{\phi}_+^* \mu_+(C)|_{E^{\ff}}= \bar{\phi}_+^* \{ \tau_+\beta +
{\cal O}_+(\lambda_f^C)\}
\end{equation*}
with $\beta=\beta_-\in\Pic(T)$.
Hence one can check that \eqref{differenceofmu:Ef} is equal to
\begin{equation}\label{differenceofmu2}
\begin{split}
\sum_{k=0}^{d-1} & \left[ \phi_-^*( \beta+{\cal O}_-(-\lambda_f^C))^k
\pmb{\cdot} \bar{\phi}_+^*( \beta+{\cal O}_+(\lambda_f^C))^{d-1-k} 
\right]_{E^{\ff}}= \\
\sum_{t=0}^{d-1} & \sideset{_{d-1}}{_t}C\cdot (\lambda_f^C)^{d-1-t}
\sum_{s=0}^{d-1-t} \left[ \beta^t \pmb{\cdot} {\cal O}_+(1)^s \pmb{\cdot}
{\cal O}_-(-1)^{d-1-t-s} \right]_{E^{\ff}}
\end{split}
\end{equation}
by using the equation
$\sum_{l=0}^t\, \sideset{_{s+l}}{_l}C \pmb{\cdot}
\sideset{_{d-1-s-l}}{_{s-l}}C = \sideset{_{d-1}}{_t}C$. 

Let $E_0^{\ff}, \ldots , E_n^{\ff}$ be the reductions of all irreducible
components of $E^{\ff}$, and let $F_0^{\ff}, \ldots, F_n^{\ff}\subset
P^{\ff}\sideset{}{_T}\times P^{-\ff}$ be their image schemes by
$\phi_-\sideset{}{_T}\times
\bar{\phi}_+: E^{\ff} \rightarrow P^{\ff} \sideset{}{_T}\times P^{-\ff}$.
\cite[Section 1]{Kl:ampleness} implies that
$\sum_{s=0}^{d-1-t} \left[ \beta^t \pmb{\cdot} {\cal O}_+(1)^s \pmb{\cdot}
{\cal O}_-(-1)^{d-1-t-s} \right]_{E^{\ff}}$
in \eqref{differenceofmu2} is equal to
\begin{equation}\label{differenceofmu3}
\sum_{i=0}^n \deg_i \sum_{s=0}^{d-1-t}
\left[ \beta^t \pmb{\cdot} {\cal O}_+(1)^s \pmb{\cdot} 
{\cal O}_-(-1)^{d-1-t-s} \right]_{F_i^{\ff}}
\end{equation}
with some rational number $\deg_i$.
We shorten $F_i^{\ff}$ to $F^{\ff}$ for the time being. We fix some
integer $M$, and divide (the right side of) \eqref{differenceofmu3} into
\begin{multline}\label{diffofmu:divide}
\sum_{s=0}^M + \sum_{s=M+1}^{d-1-t} =
\left[ \beta^t \pmb{\cdot} {\cal O}_-(-1)^{d-1-t-M} \pmb{\cdot}
\sum_{s=0}^M \bigl( {\cal O}_+(1)^s \pmb{\cdot} {\cal O}_-(-1)^{M-s} \bigr)
\right]_{F^{\ff}} \\
 + \left[ \beta^t \pmb{\cdot} {\cal O}_+(1)^{M+1} \pmb{\cdot}
\sum_{s=0}^{d-2-t-M} \bigl( {\cal O}_+(1)^s \pmb{\cdot} 
{\cal O}_-(-1)^{d-2-t-s-M} \bigr) \right]_{F^{\ff}}.
\end{multline}
\eqref{diffofmu:divide} is related to the intersection theory on
$P^{\ff}\sideset{}{_T}\times P^{-\ff}= 
\PP({\cal A}_-) \sideset{}{_T}\times \PP({\cal A}_+)$ since $F^{\ff}$
is its closed subscheme.
In this section we would like to reduce the problem of computing
\eqref{diffofmu:divide} to the intersection theory on
$\PP({\cal A}_-) \sideset{}{_T}\times \PP({\cal A}_-^{\vee})$ and
$\PP({\cal A}_+^{\vee}) \sideset{}{_T}\times \PP({\cal A}_+)$ by
choosing $M$ suitably. The reason why we would like to do so will be
explained in the next section.
It is possible to connect $P^{\ff} \sideset{}{_T}\times P^{-\ff}$ with
$\PP({\cal A}_-) \sideset{}{_T}\times \PP({\cal A}_-^{\vee})$ because
$p_g(X)>0$. 

Since $T$ is projective over $\CC$, there is a line bundle $\beta_0$ on
$T$ such that coherent ${\cal O}_T$-modules
${\cal A}_-\otimes \beta_0$,
${\cal A}_-\otimes 2\beta_0$ and 
${\cal A}_-\otimes (\beta+ \beta_0)$ 
are generated by their global sections. Because
$\beta=\{ {\cal O}_-(1)+\beta+\beta_0\} -\{ {\cal O}_-(1)+\beta_0\}$
and
${\cal O}_-(1)= 2\{ {\cal O}_-(1)+\beta_0\}- \{ {\cal O}_-(1)+2\beta_0\}$,
one can express $\beta^s \pmb{\cdot} {\cal O}_-(-1)^{d-1-t-M}$ in
\eqref{diffofmu:divide} as
\[\sum_{i=1}^I N_i \prod_{j=1}^{d-1-M} ({\cal O}_-(1)+ L_j^i)\]
with integers $N_i$ and line bundles $L_j^i$ on $T$ such that
\begin{equation}\label{genbyglsects}
\tau_{-\, *}({\cal O}_-(1)+L_j^i)=
{\cal A}_-\otimes L_j^i\;
\text{is generated by its global sections.}
\end{equation}
Hence, in order to understand the first half of \eqref{diffofmu:divide},
let us examine
\begin{equation}\label{differenceofmu4}
\left[ \prod_{j=1}^{d-1-M} ({\cal O}_-(1)+L_j) \pmb{\cdot}
\sum_{s=0}^M \bigl({\cal O}_+(1)^s \pmb{\cdot} {\cal O}_-(-1)^{M-s} \bigr)
\right]_{F^{\ff}},
\end{equation}
where $L_j\in\Pic(T)$ satisfies \eqref{genbyglsects}. We shall denote the
natural projections by $p_{\mp}: F^{\ff} \hookrightarrow P^{\ff}
\sideset{}{_T}\times P^{-\ff} \rightarrow P^{\pm\ff}$.
\eqref{differenceofmu4} clearly is zero if $d-1-M> \dim p_-(F^{\ff})$,
and so we can assume that $d-1-M\leq \dim p_-(F^{\ff})$.
Then one can find nonzero global sections
$\lambda_j\in \Gamma( P^{\ff}, {\cal O}_-(1)\otimes L_j)=
\Gamma(T, {\cal A}_-\otimes L_j)$
such that
$\dim ( F^{\ff} \cap p_-^{-1}(\Lambda_1\cap \cdots \cap \Lambda_j))=
\dim F^{\ff}-j$, where $\Lambda_j \subset P^{\ff}$ is the effective Cartier
divisor of $P^{\ff}$ corresponding to $\lambda_j$.
These $\lambda_j$ induce a homomorphism
\begin{equation}\label{otimeslambda}
\bigoplus_j \,\otimes\lambda_j : L_1^{-1} \oplus \cdots \oplus
L_{d-1-M}^{-1}\rightarrow {\cal A}_-.
\end{equation}
$\Lambda_1\cap \cdots \cap\Lambda_{d-1-M}\subset P^{\ff}$ is just a
closed subscheme
$\PP( \Cok (\bigoplus_j\, \otimes \lambda_j)) \subset \PP( A_-)$.
By a general property of intersection number
\cite[Page 297, Proposition 4]{Kl:ampleness},
\eqref{differenceofmu4} is equal to
\begin{equation}\label{diffofmu:Lambda}
\sum_{s=0}^M\, \left[ {\cal O}_+(1)^s \pmb{\cdot} {\cal O}_-(-1)^{M-s}
\right]_{p_-^{-1}(\Lambda_1\cap \cdots \cap \Lambda_{d-1-M})}.
\end{equation}
On the other hand, $\kappa\in\Gamma(K_X)$ in \eqref{pg>0} induces a
homomorphism
\begin{equation}\label{kappam}
\otimes\kappa_-: {\cal A}_+^{\vee}=Ext^1_{X_T/T}({\cal F}_0, {\cal G}_0)
\rightarrow {\cal A}_-= Ext^1_{X_T/T}({\cal F}_0, {\cal G}_0(K_X))
\end{equation}
by virtue of Proposition \ref{prop:SerreDuality}. We define $l_-$ by
$l_-= \rk(\Cok(\otimes\kappa_-))$ and prove the following proposition.
\begin{prop}\label{prop:PveetimesP}
If $d-1-M \geq l_-+ \dim T$, then we can choose $\lambda_j$ so that
\begin{equation}\label{ApveeToCok}
 {\cal A}_+^{\vee}=Ext^1_{X_T/T}({\cal F}_0, {\cal G}_0)
\overset{\otimes\kappa}{\rightarrow} {\cal A}_-=
Ext^1_{X_T/T}({\cal F}_0 \otimes {\cal G}_0(K_X)) \twoheadrightarrow
\Cok ( \bigoplus_j \, \otimes \lambda_j) 
\end{equation}
is surjective. In particular,
$p_-^{-1} (\Lambda_1\cap \cdots \Lambda_{d-1-M})$ can be regarded as
a closed subscheme of $\PP({\cal A}_+^{\vee}) \sideset{}{_T}\times
P^{-\ff}= \PP({\cal A}_+^{\vee}) \sideset{}{_T}\times \PP({\cal A}_+)$.
\end{prop}
\begin{pf}
Suppose that the following lemma is valid:
\begin{lem}\label{lem:Ti}
Define a closed subscheme $T_i$ of $T$ by
\begin{equation*}
T_i=\left\{
\begin{array}{c|c}
t\in T & \rk\Cok\bigl\{ \otimes\kappa : \Ext^1_{X_t}({\cal F}_0\otimes k(t),
{\cal G}_0\otimes k(t)) \rightarrow \\
       & \quad\Ext^1_{X_t}( {\cal F}_0\otimes k(t), {\cal G}_0(K_X)
\otimes k(t)) \bigr\} \geq l_-+i
\end{array} \right\}.
\end{equation*}
Then $\operatorname{codim}(T_i,T)\geq i$ for all $i\geq 0$.
\end{lem}
Then the dimension of a closed subscheme $\PP(\Cok\, (\otimes\kappa))$ of
$\PP({\cal A}_-)= P^{\ff}$ is less than $l_-+\dim T$ since relative Ext
sheaves ${\cal A}_-$ and ${\cal A}_+^{\vee}$ are compatible with base
change by the assumption \eqref{ExtConstant}.
Thus if $d-1-M\geq l_-+\dim T$, then one can choose $\lambda_j$ suitably
so that
$\Lambda_1\cap \cdots \cap\Lambda_{d-1-M} \cap\PP(\Cok\, (\otimes\kappa))
=\emptyset$ in $P^{\ff}$, or
$L_1^{-1}\oplus \cdots \oplus L_{d-1-M}^{-1} \overset{\oplus\otimes \lambda_j}
{\longrightarrow} {\cal A}_-= Ext^1_{X_T/T}({\cal F}_0, {\cal G}_0(K_X))
\twoheadrightarrow \Cok (\otimes\kappa)$ is surjective.
Hence also ${\cal A}_+^{\vee} \overset{\otimes\kappa}{\longrightarrow}
{\cal A}_- \twoheadrightarrow \Cok( \oplus\, \otimes\lambda_j)$
is surjective, and so the proof of Proposition \ref{prop:PveetimesP}
is completed. 

To prove Lemma \ref{lem:Ti} let us observe good properties of $\Hilb(X)$.
${\cal F}_0\otimes k(t)$ and ${\cal G}_0\otimes k(t)$ are isomorphic to,
respectively, ${\cal O}(L)\otimes I_{Z_1}$ and ${\cal O}(c_1-L)\otimes
I_{Z_2}$ for some divisor $L$ on $X_t$ and codimension-two closed
subschemes $Z_1$ and $Z_2$ in $X_t$.
The long exact sequence of Ext sheaves associated with a short exact
sequence
\begin{multline*}
0 \longrightarrow {\cal O}(c_1-L)= {\cal G}_0\otimes k(t)
\overset{\otimes\kappa}{\longrightarrow}
{\cal O}(c_1-L+K_X)\otimes I_{Z_2}= \\
{\cal G}_0(K_X)\otimes k(t)
\longrightarrow {\cal O}(c_1-L+K_X)|_{{\cal K}} \otimes I_{Z_2}
\longrightarrow 0
\end{multline*}
tells us that
\begin{multline*}
\rk \Cok\bigl\{ \kappa: \Ext^1_{X_t}({\cal F}_0\otimes k(t),
{\cal G}_0\otimes k(t)) \rightarrow \Ext^1_{X_t}({\cal F}_0\otimes k(t),
{\cal G}_0(K_X)\otimes k(t)) \bigr\}= \\
L_- -L_+-\hom_X(I_{Z_1}, {\cal O}(c_1-2L+K_X)\otimes I_{Z_2}) 
+\hom_X(I_{Z_1}, {\cal O}(c_1-2L+K_X)|_{{\cal K}}\otimes I_{Z_2}),
\end{multline*}
where $L_{\pm}$ are those of \eqref{ExtConstant}.
Since 
\[\dim\Ext_{X_t}^1({\cal G}_0 \otimes k(t), {\cal F}_0 \otimes k(t))
= \dim\Ext^1_{X_t}(I_{Z_1}, {\cal O}(c_1-2L+K_X)\otimes I_{Z_2}) \]
is independent of  $t\in T$, $\hom_X(I_{Z_1}, {\cal O}(c_1-2L+K_X)
\otimes I_{Z_2})$ is independent of $t\in T$.
Moreover, if $t\in T$ is so general that $Z_1\cap {\cal K}=Z_2\cap\KK=
\emptyset$, then $\hom_X(I_{Z_1}, {\cal O}(c_1-2L+K_X)|_{\KK}\otimes
I_{Z_2})$ is equal to $h^0({\cal O}(c_1-2L+K_X)|_{\KK})$, which is
independent of $t\in T$ since $q(X)=0$.
Therefore one can show that
\begin{multline*}
\rk \Cok\bigl\{ \otimes\kappa : \Ext^1_{X_t}({\cal F}_0\otimes k(t),
{\cal G}_0 \otimes k(t))\rightarrow
\Ext^1_{X_t}({\cal F}_0\otimes k(t), {\cal G}_0(K_X)\otimes k(t))
\bigr\} \\
-l_- =\hom_X(I_{Z_1}, {\cal O}(c_1-2L+K_X)|_{\KK}\otimes I_{Z_2})-
h^0 ({\cal O}(c_1-2L+K_X)|_{\KK}).
\end{multline*}
Now we divide Artinian schemes $Z_1$ and $Z_2$ into $Z_1= W_1 \coprod
T_1$ and $Z_2= W_2\coprod T_2$ so that, set-theoretically, $W_1=
W_2=Z_1\cap Z_2\cap\KK$.
\begin{clm}\label{claim:cok-lm}
\begin{multline*}
\hom_X(I_{Z_1}, {\cal O}(c_1-2L+K_X)|_{\KK}\otimes I_{Z_2})-
h^0({\cal O}(c_1-2L+K_X)|_{\KK}) \leq \\
l(Z_1\cap\KK)+ l(Z_2\cap\KK)+ \hom_X({\cal O}_{W_2}, {\cal O}_{W_1})
+\hom_X({\cal O}_{W_1}, \IIm (\otimes\kappa: {\cal O}_{W_2} \rightarrow
{\cal O}_{W_2})).
\end{multline*}
\end{clm}
\begin{pf}
From the long exact sequence of Tor sheaves, one derives two exact
sequences
\begin{equation}\label{F2}
 0 \longrightarrow F_2 \longrightarrow {\cal O}_{Z_2} 
\overset{\otimes\kappa}{\longrightarrow} {\cal O}_{Z_2}
\longrightarrow {\cal O}_{Z_2\cap\KK} \longrightarrow 0
\end{equation}
and
\[ 0 \longrightarrow F_2 \longrightarrow I_{Z_2}|_{\KK} \longrightarrow
L_2=\Ker( {\cal O}_{\KK} \twoheadrightarrow {\cal O}_{\KK\cap Z_2})
\longrightarrow 0. \]
Hence one can show that
{\allowdisplaybreaks
\begin{align*}
 &\hom_X(I_{Z_1}, {\cal O}(c_1-2L+K_X)|_{\KK} \otimes I_{Z_2})-
h^0( {\cal O}(c_1-2L+K_X)|_{\KK}) \\
\leq & \hom_X(I_{Z_1}, {\cal O}(c_1-2L+K_X)|_{\KK})+ \hom_X(I_{Z_1}, F_2)
-h^0({\cal O}(c_1-2L+K_X)|_{\KK}) \\
= & \bigl[ \chi(I_{Z_1}, {\cal O}(c_1-2L+K_X))- \chi(I_{Z_1}, {\cal O}(c_1-
2L)) + \ext^1_X( I_{Z_1}, {\cal O}(c_1-2L+K_X)|_{\KK}) \bigr] \\
  & \quad-\bigl[\chi({\cal O}(c_1-2L+K_X)) -\chi({\cal O}(c_1-2L))
    + h^1({\cal O}(c_1-2L+K_X)|_{\KK}) \bigr] \\
  & \quad +\chi(I_{Z_1}, F_2)+ \ext^1_X(I_{Z_1}, F_2) \\
= & \ext^1_X(I_{Z_1}, {\cal O}(c_1-2L+K_X)|_{\KK})+l(F_2)+
\ext^1_X(I_{Z_1},F_2) -h^1({\cal O}(c_1-2L+K_X)|_{\KK}) \\
= & \ext^1_X(I_{Z_1}, {\cal O}(c_1-2L+K_X)|_{\KK})- h^1({\cal O}(c_1-2L+K_X)
|_{\KK}) \\
  & \quad +l(Z_2\cap\KK)+ \ext^1_X(I_{Z_1}, F_2)
\end{align*}}
by the Riemann-Roch theorem and \eqref{F2}.
If we define $F_1$ by an exact sequence
\[ 0 \longrightarrow F_1 \longrightarrow {\cal O}_{Z_1} 
\overset{\otimes\kappa}{\longrightarrow}
{\cal O}_{Z_1} \longrightarrow {\cal O}_{Z_1\cap\KK} \longrightarrow 0, \]
then we have that
{\allowdisplaybreaks
\begin{align*}
 & \ext^1_X( I_{Z_1}, {\cal O}(c_1-2L+K_X)|_{\KK})-
h^1( {\cal O}(c_1-2L+K_X)|_{\KK}) \\
\leq & \ext^2_X({\cal O}_{Z_1}, {\cal O}(c_1-2L+K_X)|_{\KK})=
\hom_X({\cal O}(c_1-2L)|_{\KK}, {\cal O}_{Z_1}) = 
\hom_X({\cal O}_{\KK},{\cal O}_{Z_1}) \\
= & \hom_X({\cal O}_{\KK}, F_1) \leq  l(F_1)=l(Z_1\cap\KK)
\end{align*}}
For $W_2\subset Z_2$ mentioned above, there is an exact sequence
\[ 0 \longrightarrow G_2 \longrightarrow {\cal O}_{W_2}
\overset{\otimes\kappa}{\longrightarrow} {\cal O}_{W_2} \longrightarrow
{\cal O}_{W_2\cap\KK} \longrightarrow 0. \]
$\ext^1_X(I_{Z_1}, F_2)= \ext^2_X({\cal O}_{Z_1}, F_2)=
\hom_X(F_2, {\cal O}_{Z_1})$ naturally equals
$\hom_X(G_2, {\cal O}_{W_1})$.
$G_2$ induces an exact sequence
\[ 0 \longrightarrow G_2 \longrightarrow {\cal O}_{W_2}
\longrightarrow \IIm (\otimes\kappa) \longrightarrow 0. \]
This sequence implies that
{\allowdisplaybreaks
\begin{align*}
 & \hom_X(G_2, {\cal O}_{W_1}) \leq \\
 & \hom_X({\cal O}_{W_2}, {\cal O}_{W_1})- \hom_X(\IIm(\otimes\kappa),
{\cal O}_{W_1})+ \ext_X^1(\IIm(\otimes\kappa), {\cal O}_{W_1}) \\
= & -\chi(\IIm(\otimes\kappa), {\cal O}_{W_1})+ \ext^2_X(\IIm(\otimes\kappa),
{\cal O}_{W_1})+ \hom_X({\cal O}_{W_2}, {\cal O}_{W_1}) \\
= & \hom_X({\cal O}_{W_1}, \IIm(\otimes\kappa))+
\hom_X({\cal O}_{W_2}, {\cal O}_{W_1}).
\end{align*}}
Hence we conclude the proof of this claim.
\end{pf}\par
For nonnegative integers $p$, $q$ and $r$, 
\begin{equation*}
W_{pqr}^{mn}= W_{pqr}= \left\{
\begin{array}{l|l}
 (Z_1,Z_2)\in & l(Z_1\cap\KK)=p,\, l(Z_2\cap\KK)=q, \\
 \Hilb^m(X)\times & \hom({\cal O}_{W_2}, {\cal O}_{W_1})+
 \hom({\cal O}_{W_1}, \IIm(\otimes\kappa : \\
 \Hilb^n(X) & {\cal O}_{W_2} \rightarrow {\cal O}_{W_2}))=r 
\end{array} \right\}
\end{equation*}
is a locally-closed subscheme of $\Hilb^m(X)\times\Hilb^n(X)$.
By the claim above, the proof of Lemma \ref{lem:Ti} is completed if we
prove that
\begin{equation}\label{Wpqr}
\operatorname{codim}( W_{pqr}^{mn}, \Hilb^m(X)\times \Hilb^n(X))\geq
p+q+r.
\end{equation}
Let $\Hilb^m(X,x)$ denote $\Hilb^m(\operatorname{Spec}( {\cal O}_{X,x}))$
for a closed point $x\in X$, and let $Z_p^m \subset \Hilb^m(X)$ be a
locally closed subscheme
$\bigl\{ z\in\Hilb^m(X) \bigm| l(Z\cap\KK)=p \bigr\} $
for $p\in\NN$.
\begin{clm}\label{clm:Wpqrx}
If we prove that
\begin{multline}\label{Wpqrx}
\operatorname{codim} \bigl(W_{pqr}^{mn} \cap [ \Hilb^m(X,x)\times 
\Hilb^n(X,x)],\,
\Hilb^m(X,x)\times\Hilb^n(X,x) \bigr) \\
\geq p+q+r+1
\end{multline}
and that
\begin{equation}\label{Zp}
\operatorname{codim}(Z_p^m, \Hilb^m(X))\geq p,
\end{equation}
then \eqref{Wpqr} follows.
\end{clm}
\begin{pf}
The proof is by induction on $(m,n)$. Fix $(m,n)$ and suppose that
\eqref{Wpqr}
holds good for $(m',n')\neq (m,n)$ such that $m'\leq m$ and $n'\leq n$.
If either $m$ or $n$ is zero, then \eqref{Wpqr} for $(m,n)$ is immediate
from \eqref{Zp}.
Hence we assume that both $m$ and $n$ are positive. We divide the proof into
several cases. Let $(Z_1,Z_2)$ be a member of $W_{pqr}^{mn}\subset
\Hilb^m(X)\times \Hilb^n(X)$.

First, suppose that $\sharp\, \supp(Z_1)\geq 2$ and $\sharp\, \supp(Z_2)
\geq 2$. Let $m_1$, $m_2$, $n_1$ and $n_2$ be positive integers such that
$m_1+m_2=m$ and $n_1+n_2=n$.
If we define an open subset $U^{m_1}$ of $\Hilb^{m_1}(X)\times
\Hilb^{m_2}(X)$ by
\[ U^{m_1}= \left\{ (Z_1^{(1)}, Z_1^{(2)}) \bigm| Z_1^{(1)}\cap
Z_1^{(2)} =\emptyset \right\}, \]
then we can define a natural map
$\varphi^{m_1}: U^{m_1}\rightarrow \Hilb^m(X)$. Similarly we can define
$\varphi^{n_1}: U^{n_1}\rightarrow \Hilb^n(X)$.
Let $V^{m_1,n_1}$ be an open subset of $U^{m_1}\times U^{n_1}$
\[ \left\{ (Z_1^{(1)}, Z_1^{(2)}, Z_2^{(1)}, Z_2^{(2)}) \bigm|
   Z_1^{(2)}\cap Z_2^{(1)}= Z_1^{(1)}\cap Z_2^{(2)}= \emptyset \right\}. \]
$(Z_1,Z_2)$ is contained in
$\varphi^{m_1}\times \varphi^{n_1} (V^{m_1, n_1})$ for some $m_1$,
$n_1$. It's easy to prove that, in $\Hilb^{m_1}\times \Hilb^{n_1}\times
\Hilb^{m_2}\times \Hilb^{n_2}$,
\begin{equation}\label{inductionWpqr}
 (\varphi^{m_1}\times \varphi^{n_1})^{-1}( W_{pqr}^{mn})\cap V^{m_1, n_1}
\subset \underset{(p_i,q_i,r_i)}{\bigcup}
W_{p_1 q_1 r_1}^{m_1 n_1}\ \times W_{p_2 q_2 r_2}^{m_2 n_2}, 
\end{equation}
where $(p_i, q_i, r_i)$ runs over the set of all triples such that
$p_1+ p_2=p$, $q_1+ q_2=q$ and $r_1+r_2=r$.
The inductive hypothesis tells us that the dimension of the right side of
\eqref{inductionWpqr} is not exceeding $2(m+n)-(p+q+r)$, since $\dim\Hilb^n(X)
=2n$ by \cite{Fo:Hilb}.
Hence $\dim( W_{pqr}^{mn}\cap (\varphi^{m_i}\times \varphi^{n_i})
(V^{m_i,n_i})\leq 2(m+n)-(p+q+r)$.

Unless $\sharp\, \supp(Z_1)\geq 2$ and $\sharp\, \supp(Z_2)\geq 2$,
it should hold either
$\sharp\, \supp(Z_1)=1$ and $\sharp\, \supp(Z_2)\geq 2$,
$\sharp\, \supp(Z_1)\geq 2$ and $\sharp\, \supp(Z_2)=1$,
$\supp(Z_1)=\supp(Z_2)=\{ x \}$ or
$\supp(Z_1)\cap \supp(Z_2)=\emptyset$.
In all cases one can verify that $(Z_1,Z_2)$ is contained in a subscheme
whose dimension does not exceed $2(m+n)-(p+q+r)$, similarly to the case
where $\sharp\,\supp(Z_1)\geq 2$ and $\sharp\,\supp(Z_2)\geq 2$.
%
%
%
%
\end{pf}\par
\begin{clm}
Let us denote $Z_p^m \cap \Hilb^m(X,x)$ by $Z^m_p(x)$. If we prove that
\begin{equation}\label{Zpx}
\operatorname{codim}( Z^m_p(x), \Hilb^m(X,x))\geq p-1,
\end{equation}
then \eqref{Wpqrx} and \eqref{Zp} follow.
\end{clm}
\begin{pf}
We can prove \eqref{Zp} by using \eqref{Zpx} in a similar fashion to the
proof of Claim \ref{clm:Wpqrx}.
Shorten $W_{pqr}^{mn}\cap [\Hilb^m(X,x)\times \Hilb^n(X,x) ]$ to
$W_{pqr}^{mn}(x)$. If $(Z_1,Z_2)\in W_{pqr}^{mn}(x)$, then
\begin{align*}
r= & \,\hom_X({\cal O}_{Z_2}, {\cal O}_{Z_1})+ \hom_X({\cal O}_{Z_1},
\IIm( \otimes\kappa: {\cal O}_{Z_2}\rightarrow {\cal O}_{Z_2})) \\
\leq & \, l(Z_1)+l(\IIm(\otimes\kappa))= l(Z_1)+l(Z_2)-l(Z_2\cap\KK)=m+n-q.
\end{align*}
Hence if $Z_{pqr}^{mn}(x)\neq\emptyset$, then \eqref{Zpx} means that
{\allowdisplaybreaks
\begin{align*}
 & 2(m+n)-(p+q+r)\geq 2(m+n)-(p+q+m+n-q) \\
= & m+n-p =m-1-(p-1)+(n-1)+1 \\
\geq & \dim(Z_p^m(x)\times \Hilb^n(X,x))+1 \geq  \dim(W_{pqr}^{mn}(x))+1
\end{align*}}
since $\dim\Hilb^m(X,x)=m+1$ by \cite{Bc:Hilb}.
Thus \eqref{Wpqrx} follows.
\end{pf}\par
Now we prove the following claim, which completes the proof of Proposition
\ref{prop:PveetimesP} because of the claim above.
\begin{clm}\label{clm:Wqix}
For an integer $i\geq 2$ and a closed point $x\in X$, we define a locally
closed subscheme $W_{qi}^m(x)$ of $Z_q^m(x)\subset \Hilb^m(X)$ by
\[ W^m_{qi}(x)= \bigl\{ z\in Z_q^m(x) \bigm| \dim_{\CC} (I_Z\otimes k(t))=i
\bigr\}. \]
Then it holds that
\begin{equation}\label{inductiveZq}
\dim Z_q^m(x)\leq m-q \qquad (1\leq q \leq m),
\end{equation}
and that
\begin{equation}\label{Wqix}
\dim W_{qi}^m(x)\leq m-q+2-i \qquad (2\leq i, \, 1\leq q\leq m).
\end{equation}
\end{clm}
\begin{pf}
It suffices to prove this in case where $x\in\KK$.
The proof is by induction on $m$. It's easy to prove this claim for
$m=1$. Fix $m$ and suppose that this claim is valid for all $m'\leq m$.
Referring to \cite{ES:intersection}, we here recall the incidence subvariety
$H_{m,m+1}$ of $\Hilb^m(X)\times \Hilb^{m+1}(X)$:
\[ H_{m,m+1}=\left\{ (Z_1, Z_2) \in\Hilb^m(X) \times\Hilb^{m+1}(X) \bigm|
   Z_1\subset Z_2 \right\}. \]
Let $f: H_{m,m+1} \rightarrow \Hilb^m(X)$ and $g: H_{m,m+1}\rightarrow
\Hilb^{m+1}(X)$ be the projections.
There is a natural morphism $q: H_{m,m+1} \rightarrow X$
sending $(Z_1,Z_2)$ to the unique point where $Z_1$ and $Z_2$ differ.
They give a (birational) morphism
$\phi=(f,q): H_{m,m+1} \rightarrow \Hilb^m(X)\times X$.
By \cite[Section 3]{ES:intersection} it holds that
\begin{equation}\label{fibPhi}
\dim \phi^{-1}(Z_1, y)= \dim_{\CC} (I_{Z_1}\otimes k(y))-1
\end{equation}
for $(Z_1, y)\in \Hilb^m(X)\times X$, and that if $(g,q)^{-1}(Z_2,y)\neq
\emptyset$ then
\begin{equation}\label{fibGQ}
\dim (g,q)^{-1}(Z_2, y)= \dim_{\CC} (I_{Z_2}\otimes k(y))-2
\end{equation}
for $(Z_2,y)\in \Hilb^{m+1}(X)\times X$.

First let us show \eqref{inductiveZq} for $m+1$.
Suppose that $q\leq m$. Then for any $Z_2\in Z_q^{m+1}(x)$ one can
find $Z_1\in Z_q^m(x)$ such that $(Z_1,Z_2)\in H_{m,m+1}$.
%
%
Thus $Z_q^{m+1}(x)\subset g\bigl(\phi^{-1} (Z_q^m(x)\times \{ x\} )\bigr)$.
$Z_q^m(x)$ clearly is equal to $\sideset{}{_{i\geq 2}}\bigcup W_{qi}^m(x)$,
and so
$\dim Z_q^{m+1}(x) \leq \sideset{}{_{i\geq 2}}\max\, \dim \phi^{-1}
(W_{qi}^m(x)\times \{ x\} )$.
The inductive hypothesis \eqref{Wqix} and \eqref{fibPhi} imply that
\begin{multline}\label{m-q+1}
 \dim \phi^{-1}(W_{qi}^m(x) \times \{ x\} ) \\
 \leq   \dim W_{qi}^m(x)+i-1 \leq m-q+2-i+i-1=m-q+1. 
\end{multline}
Now we claim that $\dim Z_{m+1}^{m+1}(x)=0$. Indeed, if $Z_2\in
Z_{m+1}^{m+1}(x)$, then ${\cal O}_{Z_2}$ is isomorphic to
${\cal O}_{Z_2\cap\KK}$, which is equal to ${\cal O}_{\KK}/ 
m_{\KK,x}^{m+1}$ since $\KK$ is a nonsingular curve.
Therefore \eqref{inductiveZq} is valid for $m+1$.

Next let us show \eqref{Wqix} for $m+1$. If $q=m+1$ or $i=2$, then
\eqref{Wqix} results form \eqref{inductiveZq}. So suppose that
$q\leq m$ and $i\geq 3$.
If $(Z_1, Z_2)\in H_{m, m+1}$ satisfies $Z_2\in W_{qi}^{m+1}(x)$, then
$Z_1\in \Hilb^m(X,x)$, $l(Z_1\cap \KK)=q-1$ or $q$, and $\dim_{\CC} (I_{Z_1}
\otimes k(x))=i-1$, $i$, or $i+1$. Hence
\begin{equation}\label{pullbackWqix}
g^{-1} \bigl( W_{qi}^{m+1}(x) \bigr) \subset 
 \bigcup_{j=i-1}^{i+1} \phi^{-1} \left( W_{q-1, j}^m(x)\times \{ x\} 
\right) 
\cup \bigcup_{j=i-1}^{i+1} \phi^{-1} \left( W_{q,j}^m(x)\times
\{ x\} \right).
\end{equation}
If $Z_1\in Z_{q-1}^m(x)$ and $Z_2\in Z_q^{m+1}(x)$ satisfy $Z_1\subset Z_2$,
then $I_{Z_2}$ is equal to
\[ \Ker\bigl( I_{Z_1} \twoheadrightarrow I_{Z_1}|_{\KK} \twoheadrightarrow
I_{Z_1\cap\KK}= m_{\KK,x}^{q-1} \twoheadrightarrow m_{\KK,x}^{q-1}/
m_{\KK,x}^q \simeq \CC \bigr) \]
since $\KK$ is nonsingular, where $m_{\KK,x}$ is the ideal sheaf of
$x\in\KK$.
Consequently the inductive hypothesis \eqref{Wqix} implies that
\begin{multline}\label{pullbackWqix:estimate}
\dim \phi^{-1} \bigl( W_{q-1, j}^m (x)\times \{ x\}\bigr)  \cap
     g^{-1} \bigl( W_{q,i}^{m+1}(x) \bigr) \leq \\
\dim W_{q-1,j}^m(x)\times \{ x\} \leq m-q+3-j \leq m-q+1
\end{multline}
since $j\geq i-1\geq 2$.
\eqref{m-q+1}, \eqref{pullbackWqix} and \eqref{pullbackWqix:estimate}
mean that $\dim\bigl( W_{qi}^{m+1}(x) \bigr) \leq m-q+1$.
By \eqref{fibGQ} we have
\[ \dim W_{qi}^{m+1}(x) \leq m-q+1-(i-2)= m+1-q+2-i. \]
Therefore we have proved \eqref{Wqix}.
\end{pf}\par
Claim \ref{clm:Wqix} concludes the proof of Proposition \ref{prop:PveetimesP}.
\end{pf}\par
Therefore \eqref{diffofmu:Lambda}, which is the first half of
\eqref{diffofmu:divide}, is related to the intersection theory on
$\PP({\cal A}_+^{\vee}) \sideset{}{_T}\times \PP({\cal A}_+)$ if
$d-1-M\geq l_-+\dim T$.
%
%
\section{The relation to incidence varieties}\label{ss:incidence}
To understand \eqref{diffofmu:Lambda} still
more, let us examine subschemes $F^{\ff}$ and $p_-^{-1}(\Lambda_1 \cap
\cdots \cap \Lambda_{d-1-M})$ of $P^{\ff} \sideset{}{_T}\times P^{-\ff}$.
We denote the reduction of $\varphi_-^{-1}(Q^{\ff})= \tilde{\varphi}_+^{-1}
(P^{-\ff})$ by $D^{\ff}$.
\begin{lem}\label{lem:sr=0}
Let $r: Hom_{X_{D^{\ff}}/ {D^{\ff}}} (\tilde{\cal G}, \tilde{\cal G}(K_X))
\rightarrow Ext^1_{X_{D^{\ff}}/ {D^{\ff}}} (\tilde{\cal G}, \tilde{\cal F}
(K_X))$ be a homomorphism induced by the restriction of \eqref{HNFtilQm}
to $X_{D^{\ff}}$. $($We here shorten $\tilde{\cal F}|_{X_{D^{\ff}}}$ to
$\tilde{\cal F}$, etc.$)$
The extension class of the third column of \eqref{Elemtransm2} gives an
element $s$ of
\begin{align*}
 & \Ext^1_{X_{D^{\ff}}} (\tilde{\cal F}, \tilde{\cal G}(-D_-))=
\Gamma (D^{\ff}, Ext^1_{X_{D^{\ff}}/ D^{\ff}}(\tilde{\cal F}, \tilde{\cal G}
(-D_-))) \\
\simeq \,& \Gamma (D^{\ff}, Ext^1_{X_{D^{\ff}}/ D^{\ff}}(\tilde{\cal G}(-D_-),
\tilde{\cal F}(K_X))^{\vee}) \\
= \,&\Hom_{D^{\ff}} (Ext^1_{X_{D^{\ff}}/ D^{\ff}}(\tilde{\cal G}, 
\tilde{\cal F}(K_X)), {\cal O}_{D^{\ff}}(-D_-)) 
\end{align*}
by virtue of Proposition \ref{prop:SerreDuality}.
Then $s\circ r: Hom_{X_{D^{\ff}}/ D^{\ff}}(\tilde{\cal G}, \tilde{\cal G}
(K_X)) \rightarrow {\cal O}_{D^{\ff}}(-D_-)$ is zero.
\end{lem}
\begin{pf}
We shall appeal to some obstruction theory.
For a closed point $t$ of $D^{\ff}$ the third column of \eqref{Elemtransm2}
induces an exact sequence
\begin{equation}\label{inducedHNFforWpFibr}
0 \longrightarrow \tilde{\cal G}_{k(t)}(-D_-) \longrightarrow
{\cal W}_+|_{X_{k(t)}} \longrightarrow \tilde{\cal F}_{k(t)}
\longrightarrow 0
\end{equation}
on $X_{k(t)}$. As observed in the proof of Lemma \ref{lem:nontrivial}
and lemma \ref{lem:BisWp}, the extension class $\sigma$ of
\eqref{inducedHNFforWpFibr} in
$\Ext^1_{X_{k(t)}}( \tilde{\cal F}_{k(t)}, \tilde{\cal G}_{k(t)}(-D_-))$
is the obstruction to extend a morphism
\[ \Spec(A')= \Spec({\cal O}_{\tilQm}/ {\cal O}(-D_-)+ \tilde{m}_t^{l+1})
\longrightarrow D_- \overset{\varphi_-}{\longrightarrow} V_- \]
to a morphism
$\Spec(A)= \Spec({\cal O}_{\tilQm}/ \tilde{m}_t^{l+1}) \rightarrow V_-$,
where $l$ is the integer in Lemma \ref{lem:nontrivial}.
Next, let
\[ r'_t: \Ext^1_{X}( \tilde{\cal F}_{k(t)}, \tilde{\cal G}_{k(t)}
(-D_-)) \longrightarrow \Ext^2_{X}( \tilde{\cal G}_{k(t)},
\tilde{\cal G}_{k(t)}(-D_-)) \]
be an homomorphism induced by \eqref{HNFtilQm}.
Then $r'_t(\sigma)$ is the obstruction to extend
$\tilde{\cal G} \sideset{}{_{D_-}}\otimes {\cal O}_{A'} \in \Coh(X_{A'})$
to an $A$-flat family of simple sheaves on $X_A$ by
\cite[Section 2.A]{HL:text}.
Moreover, the trace map 
\begin{equation}\label{trace}
 \operatorname{tr}: \Ext^2_{X}( \tilde{\cal G}_{k(t)},
\tilde{\cal G}_{k(t)}(-D_-)) \longrightarrow
H^2( {\cal O}_{X}(-D_-)) 
\end{equation}
sends $r'_t(\sigma)$ to the obstruction to extend a line bundle
$\det ( \tilde{\cal G} \sideset{}{_{D_-}}\otimes {\cal O}_{A'})$
on $X_{A'}$ to a line bundle on $X_A$ by \cite[Theorem 4.5.3]{HL:text}.
Now $\Pic(X)$ is smooth over $\CC$, and the trace map \eqref{trace} is
isomorphic since $\rk(\tilde{\cal G}_{k(t)})=1$.
Therefore $r'_t(\sigma)=0$.
%
Remark that
\begin{equation*}
\xymatrix{
\Ext^1_X( \tilde{\cal F}_{k(t)}, \tilde{\cal G}_{k(t)}(-D_-))
\ar[d]^{\Theta} \ar[r]_{r'_t} & 
\Ext^2_X( \tilde{\cal G}_{k(t)}, \tilde{\cal G}_{k(t)}(-D_-))
\ar[d]_{\Theta} \\
\Ext^1_X(\tilde{\cal G}_{k(t)}(-D_-), \tilde{\cal F}_{k(t)}(K_X))^{\vee}
\ar[r]^{r_t^{\vee}} &
\Hom_X( \tilde{\cal G}_{k(t)}(-D_-), \tilde{\cal G}_{k(t)}(K_X))^{\vee} }
\end{equation*}
is commutative, where $\Theta$ is an isomorphism induced by the Serre
duality \eqref{SerreDuality}, and
\[ r_t: \Hom ( \tilde{\cal G}_{k(t)}(-D_-), \tilde{\cal G}_{k(t)}(K_X))
\longrightarrow \Ext^1_X( \tilde{\cal G}_{k(t)}(-D_-), \tilde{\cal F}_{k(t)}
(K_X)) \]
is defined similarly to $r$ in this lemma.
One can verify that
\begin{multline*}
 \qquad Hom_{X_{D^{\ff}}/ D^{\ff}}( \tilde{\cal G}, \tilde{\cal G}(K_X))\otimes
k(t) \overset{\operatorname{can}}{\longrightarrow} \\
\Hom_{X_{k(t)}}( \tilde{\cal G}_{k(t)}, \tilde{\cal G}_{k(t)}(K_X))
\overset{r_t^{\vee}\circ \Theta(\sigma)}{\longrightarrow}
{\cal O}_{\tilQm}(-D_-)\otimes k(t) \qquad
\end{multline*}
is equal to $(s\circ r)\otimes k(t)$. Hence $(s\circ r)\otimes k(t)=0$
for every closed point $t\in D^{\ff}$, which implies $s\circ r=0$
since $D^{\ff}$ is reduced.
\end{pf}\par
An exact sequence \eqref{ExtonPextm} on $X_{P^{\ff}}$ induces a
homomorphism
\begin{equation}\label{rP}
r_P: Hom_{X_{P^{\ff}}/ P^{\ff}}( {\cal G}_0, {\cal G}_0(K_X)) \rightarrow
Ext^1_{X_{P^{\ff}}/ P^{\ff}}( {\cal G}_0, {\cal F}_0 \otimes {\cal O}_-(1)
(K_X)).
\end{equation}
\begin{lem}\label{restrictE}
The image scheme of
$\phi_- \sideset{}{_T}\times \bar{\phi}_+: (E^{\ff})_{\operatorname{red}}
 \rightarrow P^{\ff} \sideset{}{_T}\times P^{-\ff}$ is contained in 
a subscheme 
\[ \PP( \Cok (r_P)) \subset \PP( Ext^1_{X_{P^{\ff}}/ P^{\ff}}
( {\cal G}_0, {\cal F}_0 \otimes {\cal O}_-(1)(K_X)))= P^{\ff} 
\sideset{}{_T}\times P^{-\ff} \] 
defined in \eqref{rP}.
\end{lem}
\begin{pf}
We shorten $(E^{\ff})_{\operatorname{red}}$ to $E_r$ in this proof.
There is a natural exact sequence
\begin{equation}\label{ExtonPextp}
0 \longrightarrow {\cal G}_0 \otimes {\cal O}_+(1) \longrightarrow
 {\cal V}_+ \longrightarrow {\cal F}_0 \longrightarrow 0 
\end{equation}
on $X^{-\ff}$ similarly to \eqref{ExtonPextm}.
Pulling back \eqref{ExtonPextm} and \eqref{ExtonPextp} by, respectively,
$\phi_-$ and $\bar{\phi}_+$, we have two exact sequences
\begin{align}
0 \longrightarrow {\cal F}_0\otimes {\cal O}_-(1) \longrightarrow &
\tilde{\cal V}_- \longrightarrow {\cal G}_0 
\longrightarrow 0 \label{ExtonPextm:pullbk} \\
0 \longrightarrow {\cal G}_0 \otimes {\cal O}_+(1) \longrightarrow &
\tilde{\cal V}_+ \longrightarrow {\cal F}_0 
\longrightarrow 0 \label{ExtonPextp:pullbk}
\end{align} 
on $X_{E_r}$. They induce two homomorphisms
\begin{align*}
r_E: & Hom_{X_{E_r}/ E_r} ({\cal G}_0, {\cal G}_0(K_X)) \longrightarrow 
Ext^1_{X_{E_r}/ E_r} ({\cal G}_0, {\cal F}_0\otimes {\cal O}_-(1)(K_X)) 
\quad\text{and}\\
s_E: & Ext^1_{X_{E_r}/ E_r} ({\cal G}_0, {\cal F}_0\otimes {\cal O}_-(1)(K_X))
\longrightarrow  {\cal O}_-(1)\otimes {\cal O}_+(1).
\end{align*}
We pull them back by
$\tilde{\pi}_-: D^{\ff} \rightarrow E_r$. Then
\begin{equation*}
\xymatrix{
\tilde{\pi}_-^* Hom_{X_{E_r}/ E_r} ({\cal G}_0, {\cal G}_0(K_X))
\ar[r]_-{\tilde{\pi}_-^* (s_E\circ r_E)} \ar[d]^{f_1} &
\tilde{\pi}_-^* ({\cal O}_-(1) \otimes {\cal O}_+(1)) \ar[d]_{f_2} \\
Hom_{X_{D^{\ff}}/ D^{\ff}}( \tilde{\cal G}, \tilde{\cal G}(K_X))
\ar[r]^-{s\circ r} & {\cal O}_{D^{\ff}}(-D_-) }
\end{equation*}
is commutative, where $f_1$ is a natural homomorphism and $f_2$ is the
isomorphism in Lemma \ref{lem:OeE}. One can prove this by recollecting
the way to construct $\bar{\phi}_+$ and the proof of Proposition
\ref{prop:PfisPext}. Therefore $\tilde{\pi}_-^* (s_E\circ r_E)=0$
by Lemma \ref{lem:sr=0}.

On the other hand, $\tilde{\pi}_-: (\tilde{\pi}_-)^{-1}
(E^{\ff}_{\operatorname{red}}) \rightarrow E^{\ff}_{\operatorname{red}}=
E_r$ is a principal $\bar{G}$-bundle since $E^{\ff}\subset \tilmodms$.
Thereby $(\tilde{\pi}_-)^{-1}( E_{\operatorname{red}}^{\ff})$
is reduced, and hence $(\tilde{\pi}_-)^{-1} 
(E^{\ff})_{\operatorname{red}}= D^{\ff}$. Accordingly $\tilde{\pi}_-: D^{\ff}
\rightarrow E_r$ is faithfully-flat, and so $\tilde{\pi}_-^* (s_E\circ r_E)=0$
implies $s_E\circ r_E=0$.
In fact, $s_E$ gives a morphism $\PP(s_E): E_r \rightarrow 
\PP(Ext^1_{X_{E_r}/ E_r}( {\cal G}_0, {\cal F}_0(K_X))= 
E_r \sideset{}{_T}\times P^{-\ff}$ and $\PP(s_E)$ is equal to $\id
\sideset{}{_T}\times \bar{\phi}_+$ because of its definition.
Thus $s_E\circ r_E=0$ implies this lemma.
\end{pf}\par
Here we remark that $F^{\ff}$ also is contained in $\PP(\Cok (r_P))\subset
P^{\ff} \sideset{}{_T}\times P^{-\ff}$ by virtue of its definition
and the lemma above.

Let us proceed to study a closed subscheme $p_-^{-1}(\Lambda_1\cap \cdots
\cap \Lambda_{d-1-M})$ of $F^{\ff}$ in \eqref{diffofmu:Lambda}.
We assume that $q(X)>0$, \eqref{expecteddim}, \eqref{codimgeq2}, 
\eqref{pg>0} and \eqref{ExtConstant}.
\begin{defn}
Dualizing a canonical quotient ${\cal A}_+ \sideset{}{_T}\otimes 
{\cal O}_{P^{-\ff}} \twoheadrightarrow {\cal O}_+(1)$, we have an exact
sequence
\begin{equation}\label{Cokp}
0 \longrightarrow {\cal O}_+(-1) \longrightarrow {\cal A}_+^{\vee}
\sideset{}{_T}\otimes {\cal O}_{P^{-\ff}} \longrightarrow \Cok_+
\longrightarrow 0 
\end{equation}
on $P^{-\ff}$. The {\it incidence subvariety} $\DD^{-\ff}$ of
$\PP( {\cal A}_+^{\vee} \sideset{}{_T}\otimes {\cal O}_{P^{-\ff}})=
\PP( {\cal A}_+^{\vee}) \sideset{}{_T}\times \PP( {\cal A}_+)$ is
a closed subscheme $\PP( \Cok_+)$.
\end{defn}
\begin{lem}\label{lem:containedPcok}
Suppose that the homomorphism \eqref{ApveeToCok} is surjective.
Then a closed subscheme $p_-^{-1} (\Lambda_1\cap \cdots \Lambda_{d-1-M})$
of $\PP( Ext^1_{X_T/T}( {\cal F}_0, {\cal G}_0)) \sideset{}{_T}\times
P^{-\ff}$ is contained in the incidence variety $\DD^{-\ff}$.
\end{lem}
\begin{pf}
In the proof we shorten $\Lambda_1\cap \cdots \cap \Lambda_{d-1-M}$ to
$\Lambda$ for simplicity.
From the assumption we have the following commutative diagram of
${\cal O}_T$-modules:
\begin{equation}\label{ApveeAmToCok}
\xymatrix{
{\cal A}_+^{\vee}= Ext^1_{X_T/T}( {\cal F}_0, {\cal G}_0)
\ar@{->>}[r] \ar[d]_{\otimes\kappa} & \Cok \bigl( \oplus \otimes \lambda_j)\\
{\cal A}_-= Ext^1_{X_T/T}( {\cal F}_0, {\cal G}_0(K_X))
\ar@{->>}[ur] & }
\end{equation}
This induces two closed immersions:
\begin{equation}\label{TwoImmersions}
\xymatrix{
\PP( {\cal A}_+^{\vee}) & \PP\bigl( \Cok(\oplus\otimes\lambda_j) \bigr)=
\Lambda \ar@{_{(}->}[l]^-{i^{\vee}} \ar@{_{(}->}[dl]_i \\
\PP( {\cal A}_-)= P^{\ff} & }
\end{equation}
Now we can find isomorphisms
\[j^{\vee}:   i^{\vee\, *}{\cal O}_+^{\vee}(1) \longrightarrow 
 i^* {\cal O}_-(1) \quad\text{ and}\quad
j:  i^* {\cal O}_-(1) \longrightarrow {\cal O}_{\Lambda}(1) \]
such that
\begin{equation}\label{ThreeQuotients}
\xymatrix{
{\cal A}_+^{\vee} \sideset{}{_T}\otimes {\cal O}_{\Lambda}
\ar[r]_{\otimes\kappa}  \ar@{->>}[d]_{(i^{\vee})^*\lambda_+^{\vee}}
& {\cal A}_- \sideset{}{_T}\otimes {\cal O}_{\Lambda} \ar@{->>}[r]
\ar@{->>}[d]_{i^*(\lambda_-)} &
\Cok (\oplus\otimes\lambda_j) \sideset{}{_T}\otimes {\cal O}_{\Lambda}
\ar@{->>}[d]_{\lambda_{\Lambda}} \\
(i^{\vee})^* {\cal O}^{\vee}_+(1) \ar[r]_{j^{\vee}}  &
i^* {\cal O}_-(1) \ar[r]_j & {\cal O}_{\Lambda}(1) }
\end{equation}
is commutative, where $\lambda_-$, $\lambda_+^{\vee}$ and $\lambda_{\Lambda}$
are the natural surjections on $P^{\ff}$, $\PP({\cal A}_+^{\vee})$ and
$\Lambda$, respectively.
Pull back this diagram by $\phi_-: E_{\Lambda}:= \phi_-^{-1}(\Lambda)
\rightarrow \Lambda$, which is a restriction of $\phi_-: 
E^{\ff}_{\operatorname{red}} \rightarrow P^{\ff}$.
The following diagram on $E_{\Lambda}$ is commutative:
\begin{equation}\label{incidenceHomo}
\xymatrix{
{\cal O}_+(-1)|_{E_{\Lambda}} \ar[r]_-{\operatorname{can}} \ar[dd]_l &
Hom_{X_{E_{\Lambda}}/ E_{\Lambda}} ( {\cal F}_0 \otimes {\cal O}_+(1),
{\cal F}_0 ) \ar[d]_{\otimes\kappa} \\
 & Hom_{X_{E_{\Lambda}}/ E_{\Lambda}}( {\cal F}_0 \otimes {\cal O}_+(1),
{\cal F}_0(K_X)) \ar[d]_{r'_E} \\
{\cal A}_+^{\vee} \sideset{}{_T}\otimes {\cal O}_{E_{\Lambda}} 
\ar[r]^-{\otimes\kappa} \ar@{->>}[d]_{ (i^{\vee}\circ \phi_-)^*
\lambda_+^{\vee}} & {\cal A}_- \sideset{}{_T}\otimes 
{\cal O}_{E_{\Lambda}} = Ext^1_{X_{E_{\Lambda}}/ E_{\Lambda}}
({\cal F}_0, {\cal G}_0(K_x)) \ar@{->>}[d]_{ (i\circ\phi_-)^* \lambda_-} \\
{\cal O}_+^{\vee}(1) \ar[r]^{\phi_-^* (j^{\vee})} & {\cal O}_-(1)}
\end{equation}
Here, $l$ is the pull back of ${\cal O}_+(-1) \rightarrow {\cal A}_+^{\vee}
\sideset{}{_T}\otimes {\cal O}_{P^{-\ff}}$ in \eqref{Cokp} by
$E_{\Lambda} \hookrightarrow E^{\ff} \rightarrow P^{-\ff}$,
$r'_E$ is defined by using \eqref{ExtonPextp:pullbk} similarly to
$r_P$ \eqref{rP}, and the lower diagram is obtained from the left side
of \eqref{incidenceHomo}.
Then, $(i\circ\phi_-)^* \lambda_-$ in \eqref{incidenceHomo} coincides
with the homomorphism
$s'_E: Ext^1_{X_{E_{\Lambda}}/ E_{\Lambda}}( {\cal F}_0, {\cal G}_0(K_X))
\twoheadrightarrow {\cal O}_-(1)$ defined by using \eqref{ExtonPextm:pullbk}
similarly to $s_E$ in the proof of Lemma \ref{lem:containedPcok}.
Hence one can verify that $(i\circ\phi_-)^* \lambda_- \circ r'_E=0$ in
the same way as the proof of Lemma \ref{lem:containedPcok}, which implies
$(i^{\vee}\circ\phi_-)^* \lambda_+^{\vee}\circ l=0$ in
\eqref{incidenceHomo}.
This and \eqref{ThreeQuotients} imply that
\[ {\cal O}_+(-1)|_{E_{\Lambda}} \overset{l}{\longrightarrow} 
{\cal A}_+^{\vee}\sideset{}{_T}\otimes {\cal O}_{E_{\Lambda}}
\twoheadrightarrow \Cok(\oplus\otimes\lambda_j) \sideset{}{_T}\otimes
{\cal O}_{E_{\Lambda}} \overset{\phi_-^*(\lambda_{\Lambda})}{\longrightarrow}
{\cal O}_{\Lambda}(1) \]
is the zero map. By this we can conclude the proof.
\end{pf}\par
For the time being we suppose a homomorphism \eqref{ApveeToCok} is
surjective. Moreover, we assume that
$\dim F^{\ff}= \dim E^{\ff}= d-1$ since \eqref{differenceofmu3} is zero
unless this holds good.
Then, by the lemma above a subvariety $p_-^{-1}(\Lambda_1 \cap \cdots \cap
\Lambda_{d-1-M})$ of $\PP({\cal A}_+^{\vee}) \sideset{}{_T}\times
\PP({\cal A}_+)$ gives an algebraic cycle $\omega\in A^r(\DD^{-\ff})$
of the incidence variety $\DD^{-\ff}$ with
$r= \operatorname{codim}( p^{-1}_-(\Lambda_1 \cap \cdots \cap 
\Lambda_{d-1-M}), \DD^{-\ff})$.
$\DD^{-\ff}$ is nonsingular, so we can use the intersection
theory of $\DD^{\ff}$.
Because ${\cal O}_-(1)|_{\Lambda_1 \cap \cdots \cap \Lambda_{d-1-M}}=
{\cal O}_+^{\vee}(1)|_{\Lambda_1 \cap \cdots \cap \Lambda_{d-1-M}}$ as
mentioned in \eqref{ThreeQuotients}, one can verify that
\eqref{diffofmu:Lambda} is equal to
\begin{equation}\label{ToIntTheory}
\sum_{t=0}^M \deg \left( c_1({\cal O}_+(1))^t \pmb{\cdot} 
c_1( {\cal O}_+^{\vee}(-1))^{M-t} \pmb{\cdot} \omega \right)_{\DD^{-\ff}},
\end{equation}
where $(\,)_{\DD^{\ff}}$ designates the multiplication in the Chow ring
$A(\DD^{-\ff})$. We shall omit $c_1(\, )$ from now on.
%
%
Moreover, one can write $\omega\in A^r(\DD^{-\ff})$ as
$\omega= \sum_{j=0}^{L_+-2} b_j \, {\cal O}_+^{\vee}(1)^j$
with some $b_j \in A^{r-j}(P^{-\ff})=A^{r-j}(\PP({\cal A}_+))$ because
the sheaf $\Cok_+$ in \eqref{Cokp} is a vector bundle on $P^{-\ff}$
whose rank is $L_+-1= \rk {\cal A}_+^{\vee}-1$.
By \eqref{Cokp}, $(-1)^M$ times \eqref{ToIntTheory} is equal to
{\allowdisplaybreaks
\begin{gather}\label{IntTheory1}
\begin{split}
 & \deg \sum_{j=0}^{L_+-2} \left( b_j \pmb{\cdot} {\cal O}_+^{\vee}(1)^j 
\pmb{\cdot}
\sum_{t=0}^M {\cal O}_+^{\vee}(1)^t \pmb{\cdot} {\cal O}_+(-1)^{M-t} 
\right)_{\DD^{-\ff}} \\
 = & \deg \sum_{j=0}^{L_+-2} \left( b_j \pmb{\cdot} \sum_{t=0}^{M+j}
{\cal O}_+^{\vee}(1)^t \pmb{\cdot} {\cal O}_+(-1)^{M+j-t} 
\right)_{\DD^{-\ff}} \\
   & -\deg \sum_{j=0}^{L_+-2} \left( b_j \pmb{\cdot} \sum_{t=0}^{j-1} 
  {\cal O}_+^{\vee}(1)^t \pmb{\cdot} {\cal O}_+(-1)^{M+j-t} 
\right)_{\DD^{-\ff}} \\
= & \deg \sum_{j=0}^{L_+-2} \left( b_j \pmb{\cdot} \sum_{t=0}^{M+j} 
{\cal O}_{\PP(\Cok_+)}(1)^t \pmb{\cdot} {\cal O}_+(-1)^{M+j-t} 
\right)_{\PP(\Cok_+)} \\
  & -\deg \sum_{j=0}^{L_+-2} \left( b_j \pmb{\cdot} 
{\cal O}_+(-1)^{M-1} \pmb{\cdot} \sum_{t=0}^{j-1} 
{\cal O}_{\PP(\Cok_+)} (1)^t \pmb{\cdot} {\cal O}_+(-1)^{j-1-t} 
\right)_{\PP(\Cok_+)} \\
= & \deg \sum_{j=0}^{L_+-2} \left( b_j \pmb{\cdot} \sum_{t=0}^{M-j}
s_{t-(L_+-2)} ( (\Cok_+)^{\vee}) \pmb{\cdot} {\cal O}_+(-1)^{M+j-t}
 \right)_{P^{-\ff}} \\
  & -\deg \sum_{j=0}^{L_+-2} \left( b_j \pmb{\cdot} 
{\cal O}_+(-1)^{M-1} \pmb{\cdot} \sum_{t=0}^{j-1} s_{t-(L_+-2)}
( (\Cok_+)^{\vee}) \pmb{\cdot} {\cal O}_+(-1)^{j-1-t} 
\right)_{P^{-\ff}}. 
\end{split}
\end{gather}}
Here $s_l((\Cok_+)^{\vee})\in A^l(P^{-\ff})$ is the Segre class
of a vector bundle $(\Cok_+)^{\vee}$ on $P^{-\ff}$, which is explained in
\cite[Section 3.1]{Fu:intersection}.

In general, the Chern polynomial $c_t({\cal V})= \sum_{j=0}^{\infty}
c_j({\cal V})\, t^j$ of a vector bundle ${\cal V}$ satisfies that
$c_t({\cal V})^{-1}= \sum_{j=0}^{\infty} s_j({\cal V})\, t^j$
as power serieses. Thus the dual of \eqref{Cokp} tells us that
\[ s_t(( \Cok_+)^{\vee}) \cdot (1+ \sum_{j>0} {\cal O}_+(-1)^j\,
t^j)= s_t({\cal A}_+ \sideset{}{_T}\otimes {\cal O}_{P^{-\ff}}). \]
In addition, $s_j({\cal V})=0$ if $j<0$. 
We see that \eqref{IntTheory1} is
equal to the degree of
{\allowdisplaybreaks
\begin{gather}\label{IntTheory2}
\begin{split}
 &\sum_{j=0}^{L_+-2} \left( b_j \pmb{\cdot} s_{M+j-(L_+-2)}({\cal A}_+
\sideset{}{_T}\otimes {\cal O}_{P^{-\ff}}) \right)_{P^{-\ff}} \\
 & \quad-\sum_{j=0}^{L_+-2} \left( b_j \pmb{\cdot} {\cal O}_+(-1)^{M+1} 
\pmb{\cdot} s_{j-1-(L_+-2)}({\cal A}_+ \sideset{}{_T}\otimes 
{\cal O}_{P^{-\ff}} ) 
\right)_{P^{-\ff}} \\
= & \sum_{j=0}^{L_+-2} \left( b_j \pmb{\cdot} s_{M+j-(L_+-2)}({\cal A}_+
\sideset{}{_T}\otimes {\cal O}_{P^{-\ff}} )\right)_{P^{-\ff}}, 
\end{split}
\end{gather}}
taking into account that $j-1-(L_+-2)< 0$. 
Since ${\cal A}_+$ is a vector bundle on $T$, 
$s_{M+j-(L_+-2)}( {\cal A}_+ \sideset{}{_T}\otimes
{\cal O}_{P^{-\ff}})=0$ provided that $M+j-(L_+-2)\geq M-(L_+-2)>
\dim T$. Therefore we obtain the following proposition as a result
of \eqref{diffofmu:divide}, Proposition \ref{prop:PveetimesP},
\eqref{ToIntTheory}, \eqref{IntTheory2}, etc.
\begin{prop}\label{prop:FirstTermIsZero}
If $d-1-M \geq l_-+ \dim T$ and $M-(L_+-2)> \dim T$, then the first
term $\sum_{s=0}^M$ in \eqref{diffofmu:divide} is zero.
\end{prop}
Furthermore, the second term $\sum_{s=M+1}^{d-1-t}$ of
\eqref{diffofmu:divide},
\begin{equation}\label{SecondTermIsZero}
\left[ \beta^t \pmb{\cdot} {\cal O}_+(1)^{M+1} \pmb{\cdot} \sum_{s=0}^{d-2-t-M}
{\cal O}_+(1)^s \pmb{\cdot} {\cal O}_-(-1)^{d-2-t-s-m} \right]_{F^{\ff}},
\end{equation}
clearly is zero if $M+1> \dim P^{-\ff}= L_+-1+ \dim T$.
\begin{prop}\label{prop:ContributionIsZero}
The contribution of $E^{\ff}$ to $\mu_-(C)^{d(c_2)} -\mu_+(C)^{d(c_2)}$, 
that is \eqref{differenceofmu:Ef}, is equal to zero if $d\geq L_+-l_- +
2\dim T$.
\end{prop}
\begin{pf}
By Proposition \ref{prop:FirstTermIsZero} and \eqref{SecondTermIsZero},
\eqref{diffofmu:divide} is equal to zero if $M\leq d-1-l_- -\dim T$
and $M> L_+-2+\dim T$. One can find such an integer $M$ if
$d-1-l_- -\dim T > L_+ -2+\dim T$.
\end{pf}\par
As observed before Claim \ref{claim:cok-lm},
{\allowdisplaybreaks
\begin{align*}
l_- = & L_- -L_+ -\hom_X(I_{Z_1}, {\cal O}(c_1-2L+K_X)\otimes I_{Z_2})
 + h^0( {\cal O}(c_1-2L+K_X)|_{\KK}) \\
 = & -\chi(I_{Z_1}, {\cal O}(c_1-2L+K_X) \otimes I_{Z_2}) -L_+ +
 h^0( {\cal O}(c_1-2L+K_X)|_{\KK}),
\end{align*}}pp
where $L$ and $Z_i$ satisfies that $([2L-c_1], l(Z_1), l(Z_2))=\ff=
(f,m,n)$.
From the Riemann-Roch theorem and Clifford's theorem
\cite[Theorem IV.5.4]{Ha:text}, we deduce that
{\allowdisplaybreaks
\begin{gather}\label{lm+Lp}
\begin{align*}
l_- +L_+=  &-f \cdot (f-K_X)/2 +(m+n) -\chi(\Ox) +h^0( {\cal O}(c_1-
2L+K_X)|_{\KK}) \\
 \leq & -f\cdot(f-K_X)/2 +(m+n) -\chi(\Ox)  \\
      & \qquad +\max( -K_X\cdot f, (K_X-f)\cdot K_X/2 +1, 0).
\end{align*}
\end{gather}}
On the other hand, $\dim T= \dim( \Pic(X)\times \Hilb^m(X)\times 
\Hilb^n(X))=2(m+n)$ and $m+n=c_2+(f^2-c_1^2)/4$ since $(f,m,n)\in A^+(a)$.
Therefore one can verify that
\begin{equation}\label{EstimateTheCondition}
d-(L_+ +l_- +2\dim T) \geq -c_2-(3/4) f^2 -2\chi(\Ox)+
\min(\pm K_X\cdot f/2, -K_X^2/2 -1).
\end{equation}
%
%
Now fix a compact subset ${\cal S}$ in $\Amp(X)$. Then one can find a constant
$d_0({\cal S})$ depending on ${\cal S}$ such that
$|f\cdot K_X| \leq d_0({\cal S})\cdot
\sqrt{ -f^2}$ if $W^f \cap {\cal S}\neq\emptyset$, as shown in the proof of
\cite[Lemma 2.1]{Qi:birational}.
Hence one can find constants $d_1({\cal S})$ and $d_2({\cal S})$
depending on ${\cal S}$ such
that if $-f^2 > (4/3) c_2 +d_1({\cal S})\sqrt{ c_2} +d_2({\cal S})$, then
\eqref{lm+Lp} is greater than zero.

Therefore we arrive at Proposition \ref{thm:MainIntro} in Introduction,
which is the observation of Proposition \ref{prop:AGinterpretIntro}
in algebro-geometric view.
\begin{rem}\label{rem:K3}
Suppose that $X$ is $K3$ surface and that assumptions \eqref{expecteddim}
and \eqref{codimgeq2} hold good for $(0,c_2)$. Then \eqref{pg>0} and
\eqref{ExtConstant} are always valid, and furthermore, the homomorphism
\eqref{ApveeToCok} is always surjective. (It is not necessary to assume
that $d-1-M\geq l_- +\dim T$.) Thus one can prove
$\gamma_{H_-}(c_2)=\gamma_{H_+}(c_2)$.

\end{rem}
\providecommand{\bysame}{\leavevmode\hbox to3em{\hrulefill}\thinspace}
\providecommand{\MR}{\relax\ifhmode\unskip\space\fi MR }
\providecommand{\MRhref}[2]{%
  \href{http://www.ams.org/mathscinet-getitem?mr=#1}{#2}
}
\providecommand{\href}[2]{#2}

\end{document}